\documentclass[12pt,a4paper]{amsart}

\usepackage{euscript,amsfonts,amssymb,amsmath,amscd,amsthm,enumerate,hyperref}

\usepackage{comment}

\usepackage{mathrsfs}

\usepackage{graphicx}

\usepackage{color}

\usepackage[margin=1.1in]{geometry}

\usepackage{bbm}

\newtheorem{thm}{Theorem}[section]
\newtheorem{lemma}[thm]{Lemma}
\newtheorem{prop}[thm]{Proposition}
\newtheorem{cor}[thm]{Corollary}
\newtheorem{asmp}[thm]{Assumption}
\newtheorem{definition}[thm]{Definition}

\newtheorem{rmk}[thm]{Remark}

\newcommand{\del}{{\delta}}
\newcommand{\eps}{{\epsilon}}
\newcommand{\eq}{\begin{equation}}
\newcommand{\en}{\end{equation}}
\newcommand{\eqs}{\begin{equation*}}
\newcommand{\ens}{\end{equation*}}
\newcommand{\rr}{\mathbb{R}}
\newcommand{\pp}{\mathbb{P}}
\newcommand{\qq}{\mathbb{Q}}
\newcommand{\nn}{\mathbb{N}}
\newcommand{\zz}{\mathbb{Z}}

\newcommand{\ev}{\mathbb E}  
\newcommand{\oS}{\overline{\mathcal S}}
\newcommand{\oSe}{\oS}
\newcommand{\DxS}{\oS_x}

\newcommand{\CC}{\mathcal C}

\newcommand{\cA}{\mathcal A}
\newcommand{\tcA}{\widetilde{\mathcal A}}
\newcommand{\ieta}{{\iota_\star}}
\newcommand{\wt}{\widetilde}
\newcommand{\Rg}{\red{R^{(\gamma)}}}
\newcommand{\Rro}{\red{R^{(\rho)}}}
\newcommand{\Rxi}{\red{R^{(\xi)}}}
\newcommand{\Rgx}{\red{R_x^{(\gamma)}}}
\newcommand{\Rgs}{\red{R^{(\gamma_{\star})}}}
\newcommand{\Rwg}{\red{R^{(\widehat{\gamma})}}}
\newcommand{\Rge}{\red{R^{(\gamma^{\epsilon})}}}
\newcommand{\Rs}{R_\star}
\newcommand{\Mstar}{M_1^{(0)}(\rr)}
\newcommand{\Miota}{M_1^{(\iota)}(\rr)}
\newcommand{\Miotas}{M_1^{(\iota-1)}(\rr)}
\newcommand{\Mieta}{M_1^{(\ieta)}(\rr)}
\newcommand{\rre}{\rr_T}
\newcommand{\FF}{\mathcal F}
\newcommand{\GG}{\mathcal G}
\newcommand{\GGc}{{\mathcal G}^\star}
\newcommand{\RR}{\mathcal R}
\newcommand{\LL}{\mathcal L}
\newcommand{\HH}{\mathcal H}
\newcommand{\EE}{\mathcal E}
\def\v{\textsf{v}}
\newcommand{\ess}{\mathrm{ess}\!\!\!\!\!}
\newcommand{\ep}{\hfill $\Box$}
\newcommand{\abbr}[1]{{\sc\lowercase{#1}}}
\newcommand{\red}[1]{{
#1}}

\numberwithin{equation}{section}

\title[Interacting diffusions]{Large deviations for diffusions interacting through their ranks} 

\author{\quad\;Amir Dembo}
\address{Departments of Statistics and Mathematics, Stanford University, Stanford, CA 94305}
\email{adembo@stat.stanford.edu}

\author{Mykhaylo Shkolnikov}
\address{Department of Statistics, University of California Berkeley, CA 94720-3860}
\email{mshkolni@gmail.com}

\author{S. R. Srinivasa Varadhan}
\address{Courant Institute, New York University, New York, NY 10012-1185}
\email{varadhan@cims.nyu.edu}

\author{Ofer Zeitouni}
\address{School of Mathematics, University of Minnesota, Minneapolis, MN 55455 
and Faculty of Mathematics, Weizmann Institute of Science, Rehovot, Israel 76100}
\email{zeitouni@umn.edu}

\thanks{We thank N. Krylov for suggesting the method applied in the proof of Lemma \ref{lem-l3rx} (which is a key ingredient in the proof of Prop \ref{L3prop}(A)) and for describing to one of us some of the details needed for the implementation of the method. Our research was partially supported by NSF grants \#DMS-1106627 (Dembo and Shkolnikov), \#DMS-1208334 (Varadhan), \#DMS-0804133 (Zeitouni), and the Israel Science Foundation
grant \#111/11 (Zeitouni).}

\begin{document}

\maketitle

\begin{abstract}
We prove a Large Deviations Principle (\abbr{LDP}) for systems of diffusions (particles) interacting through their ranks, when the number of particles tends to infinity. We show that the limiting particle density is given by the unique solution of the appropriate McKean-Vlasov equation and that the corresponding cumulative distribution function evolves according to a non-degenerate generalized porous medium equation with convection. The large deviations rate function is provided in explicit form. This is the first instance of a \abbr{LDP} for interacting diffusions, where the interaction occurs both through the drift and the diffusion coefficients and where the rate function can be given explicitly. In the course of the proof, we obtain new regularity results for tilted versions of such generalized porous medium equation.
\end{abstract}

\section{Introduction}\label{sec_intro}

Recently, systems of diffusion processes (particles) interacting through their ranks have received much attention. For a fixed number of particles $N\in\nn$, these are given by the unique weak solution of the stochastic differential system (\abbr{sds}),
\eq\label{rank_sde}
\mathrm{d}X_i(t)=\sum_{j=1}^N b_j\,\mathbf{1}_{\{X_i(t)=X_{(j)}(t)\}}\,\mathrm{d}t
+\sum_{j=1}^N \sigma_j\,\mathbf{1}_{\{X_i(t)=X_{(j)}(t)\}}\,\mathrm{d}W_i(t),\;\; i=1,\ldots,N,
\en
where $b_1,b_2,\ldots,b_N$ are arbitrary real constants, $\sigma_1,\sigma_2,\ldots,\sigma_N$ are arbitrary positive constants, $W_1,W_2,\ldots,W_N$ are independent standard Brownian motions and $X_{(1)}(t)\leq X_{(2)}(t)\leq\ldots\leq X_{(N)}(t)$ are the ordered particles at time $t$. In this paper we study the behavior of the solution to \eqref{rank_sde} as $N$ becomes large in the regime when $|b_{j+1}-b_j|+|\sigma_{j+1}^2-\sigma_j^2|=O(N^{-1})$ for all $j$ (see Assumption 
\ref{main_ass} 
below for the details). In other words, the drift and the diffusion coefficients of a particle change slowly as it changes its rank in the particle system. 

The existence and uniqueness of the weak solution to \eqref{rank_sde} was shown in the work \cite{bp}, which was motivated by questions in filtering theory. The system 
\eqref{rank_sde} has also reappeared in the context of stochastic portfolio theory under the name \textit{first-order market model} (see the book \cite{fe} and the survey article \cite{fk}). In the latter context, our choice of the regime $|b_{j+1}-b_j|+|\sigma_{j+1}^2-\sigma_j^2|=O(N^{-1})$ agrees with the economic intuition that 
a small change of a company's rank cannot lead to a large jump in the growth rate and the volatility coefficient of its market capitalization. Due to its central role in the analysis of capital distributions in financial markets and long-term portfolio performance therein, as well as its intriguing mathematical features, the ergodicity and sample path properties of this model have undergone a detailed analysis in the case that the number of particles is fixed (see \cite{cp,cp2,ik,iks,ipbkf}). Moreover, concentration properties of the solution to \eqref{rank_sde} for large values of $N$ have been studied in \cite{ps}, 
and an analogous infinite particle system has been constructed and analyzed in \cite{pp}.

\smallskip
In \cite{sh} it was observed that the \abbr{SDS} \eqref{rank_sde} can be rewritten as
\eq\label{sde}
\mathrm{d}X_i(t)=b(F_{\rho^N(t)}(X_i(t)))\,\mathrm{d}t + \sigma(F_{\rho^N(t)}(X_i(t)))\,\mathrm{d}W_i(t),\quad i=1,2,\ldots,N,
\en
where $\rho^N(t)=\frac{1}{N}\sum_{i=1}^N \delta_{X_i(t)}$ is the \textit{empirical measure} of the particle system at time $t$, $F_{\rho^N(t)}$ is the corresponding cumulative distribution function (\abbr{cdf}), and $b:\,[0,1]\rightarrow\rr$, $\sigma:\,[0,1]\rightarrow(0,\infty)$ are functions satisfying $b\big(\frac{j}{N}\big)=b_j$, $\sigma\big(\frac{j}{N}\big)=\sigma_j$ for all $j=1,2,\ldots,N$. The representation \eqref{sde} allows us to view the particle system \eqref{rank_sde} as a system of diffusion processes interacting through their mean-field.
It gives rise to questions on the large $N$ behavior of the empirical measure in \eqref{sde}, where the mathematical challenge is due to the discontinuity of the diffusion coefficients in \eqref{sde}.
A 
\red{Law of Large Numbers (\abbr{lln})}
for $\rho^N(t)$ is obtained in \cite{sh} 
\red{for non-decreasing $i \mapsto X_i(0)$, with $\{X_i(0)-X_1(0)\}$
chosen according to the stationary distribution of} 
the process of spacings between consecutive ordered particles in \eqref{sde} 
(in particular, assuming that $b$ \red{and $\sigma$ are} such that a stationary distribution exists). In this case, it was shown that the limiting particle measure \red{path $t \mapsto \gamma(t)(\cdot)$ 
satisfies} the \textit{McKean-Vlasov equation}
\eq
\int_\rr f\,\mathrm{d}\gamma(t)-\int_\rr f\,\mathrm{d}\gamma(0)
=\int_0^t \red{\mathrm{d}s} \, 
\int_\rr \big[b(F_{\gamma(s)}(\cdot))f'+\frac{1}{2} \sigma(F_{\gamma(s)}(\cdot))^2 f''\big]\,\mathrm{d}\gamma(s)  
\en
for all Schwartz functions \red{$x \mapsto f(x)$ and any} $t\ge 0$
\red{(hereafter $\mathrm{d}\gamma(s)$ is short hand for integration 
with respect to the probability measure $\gamma(s)(\cdot)$ on $\rr$, 
at the fixed time $s$)}.
Further, the corresponding \red{\abbr{cdf}-path 
$\Rg(t,x):=F_{\gamma(t)}(x)$
evolves} according to the
\textit{
non-degenerate generalized porous medium equation with convection} (see \cite{Va} and the references therein for a thorough treatment),
\eq\label{PME}
R_t=(\Sigma(R))_{xx}-(\Gamma(R))_x\,,
\en
where $\Sigma(r)=\int_0^r \frac{1}{2}\sigma^2(u)\,\mathrm{d}u$ and $\Gamma(r)=\int_0^r b(u)\,\mathrm{d}u$.

\smallskip
In this paper, we establish a Large Deviations Principle (\abbr{LDP}) for the sequence 
$\{\rho^N$, $N\in\nn\}$ of paths $\{\rho^N(t), t\in[0,T]\}$ of empirical measures, where $T>0$ is arbitrary, but fixed throughout. Among other things, this \abbr{LDP} implies the \abbr{lln} 
for $\{\rho^N,\;N\in\nn\}$. Such \abbr{lln}, 
and the \abbr{LDP} upper bound are shown under 
the mild regularity Asmp. \ref{main_ass}
on the functions $b$, $\sigma$ and the initial 
empirical measures $\{\rho^N(0),\;N\in\nn\}$ 
(dispensing of 
the stationarity assumption on the process of spacing, which plays a crucial role
in \cite{sh}). Our next definition is useful for stating Asmp. \ref{main_ass}.
\begin{definition} 
Let $M_1(\rr)$ denote the space of Borel probability measures on $\rr$, 
endowed with the L\'{e}vy distance metric
\eqs
d_L(\alpha_1,\alpha_2):=\inf\{\epsilon>0|
\alpha_1(O)\leq \alpha_2(O_\epsilon)+\epsilon, 
\alpha_2(O)\leq \alpha_1(O_\epsilon)+\epsilon\;
\forall\,{\rm open}\,O\subset\rr
\}
\ens
(where $O_\epsilon$ stands for the $\epsilon$-neighborhood of $O$ in $\rr$).
Then, for $\iota \in [0,1]$ let 
$\Miota$ denote the subset of all $\mu \in M_1(\rr)$ such that 
$\int_\rr |x|^{1+\iota} \,\mathrm{d}\mu<\infty$ and
$\frac{\mathrm{d} \mu}{\mathrm{d} x} \in L^{q_0}(\rr)$ for some $q_0>1$.
\end{definition}

\begin{asmp}\label{main_ass}
\red{$\sigma_j=\sigma(\frac{j}{N})$ and $b_j=b(\frac{j}{N})$, where:} 
\newline 
(a) The function $A:=\frac{1}{2}\sigma^2$ is uniformly bounded below 
by some $\underline{a}>0$, and $A'$ is bounded on $[0,1]$.
\newline 
(b) The function 
$b(\cdot)$ is Lipschitz continuous on $[0,1]$.
\newline
(c) As $N\rightarrow\infty$
the deterministic initial empirical measures $\{\rho^N(0)\}$ converge weakly 
to some $\rho_0 \in \Mieta$, $\ieta \in (0,1]$,
and $\sup_{N\in\nn} \int_\rr |x|^{1+\ieta}\,\mathrm{d}\rho^N(0)<\infty$. 
\newline
(d) The function $A'$ is further Lipschitz continuous on $[0,1]$.
\end{asmp}

Throughout the paper we 
let $\CC=C([0,T],M_1(\rr))$ stand for the space of 
continuous functions from $[0,T]$ to $M_1(\rr)$ endowed with the metric
\eq\label{eq:dL-sup-def}
d(\gamma_1(\cdot),\gamma_2(\cdot))=\sup_{t\in[0,T]} d_L(\gamma_1(t),\gamma_2(t))\,.
\en
Further, throughout we adopt the convention that $\frac{0}{0}=0$,
and identify each $\gamma\in \CC$ with the \red{corresponding
\abbr{cdf}-path $\Rg (t,x):=F_{\gamma(t)}(x)$}. The following spaces
of functions then play a major role in our \abbr{LDP} on $\CC$.

\begin{definition}\label{dfn:cA-FF-J}
Let $\oSe$ denote the space of functions $g$ on $\rre := [0,T] \times \rr$
which are infinitely differentiable
and such that, for all $t\in[0,T]$, $g(t,\cdot)$ is a Schwartz function on $\rr$.
Next, let
\begin{align*}
\tcA := \{  \gamma\in \CC : &  
\; \Rg \in C_b(\rr_T), \; \; \gamma(0) \in \Mstar, \\
&
t\mapsto\int_\rr g(t,\cdot)\,\mathrm{d}\gamma(t) \;\; \text{ abs. cont. on }
[0,T] 
\;\; \forall g\in\oS \},
\end{align*}
and for each $\iota \in (0,1]$, $\mu \in \Mstar$, its subset
\eq\label{eq:eta-mom}
\cA_{\iota,\mu} := \{ \gamma \in \tcA : \; \gamma(0)=\mu, \;\; 
\int_{\rr_T} |x|^{1+\iota}\,\mathrm{d}\gamma(t)\,\mathrm{d}t <\infty
\} \,.
\en
Further, we use 
\begin{align}
\FF_q := \Big\{R = \Rg
: &\gamma \in \CC,
\nonumber \\ & 
R_t,\,R_{xx}\in L^{q}(\rr_T),\;R_x\in L^3(\rr_T),\;
\frac{R_t^2}{R_x},\frac{R_{xx}^2}{R_x}\in L^1(\rr_T)\Big\} 
\label{eq:thm-13}
\end{align}
with $\FF:=\FF_{3/2}$ and
\eq\label{eq:J-def}
\widetilde{J} (\gamma):=\left\{\begin{array}{ll}
    \frac{1}{4}\Big\|\frac{R_t-(A(R)R_x)_x+b(R)R_x}{(A(R)R_x)^{1/2}}\Big\|_{L^2(\rr_T)}^2, & \;\textrm{if}\;\;\gamma\in \tcA, \;
     R = \Rg \in \FF \\
    \infty\;, & \;\textrm{otherwise},
   \end{array}
  \right.
\en 
with $J_{\iota,\mu}(\gamma)$ 
\red{defined as in \eqref{eq:J-def} except for replacing there
$\tcA$ by the smaller $\cA_{\iota,\mu}$. 
} 
\end{definition}
Our main result then reads as follows.
\begin{thm}\label{main_thm}
Under Asmp. \ref{main_ass} with $\ieta=1$, the sequence 
$\{\rho^N$, $N\in\nn\}$
satisfies the \abbr{LDP} on $\CC$ with scale $N$
and good rate function \red{$J_{1,\rho_0}(\cdot)$} 
of Definition \ref{dfn:cA-FF-J}.
\end{thm}

\red{
\begin{rmk}\label{rmk:ieta-small} 
As shown in Prop. \ref{local_ubd_exp_tight} (and Corollary \ref{rate_fct_prop}),
the exponential tightness of $\{\rho^N\}$ and the 
\abbr{LDP} upper bound of Theorem \ref{main_thm} 
with rate $J_{\ieta,\rho_0}(\cdot)$ 
apply for any value of $\ieta>0$ in Asmp. \ref{main_ass}(c) 
and do not require part (d) of Asmp. \ref{main_ass}. 
\end{rmk}
In view of the preceding remark we have} the following \abbr{lln}. 
\begin{cor}\label{LLN}
Under Asmp. \ref{main_ass}(a)-(c) \red{for some $\ieta \in (0,1]$}, 
the sequence $\{\rho^N$, $N\in\nn\}$ 
converges almost surely 
to the unique path \red{$\gamma_{\star} \in \cA_{\ieta,\rho_0}$ for which 
$\Rgs \in \FF$} 
is a generalized solution of the Cauchy problem
\eq\label{PMEcauchy}
R_t=(A(R)R_x)_x-b(R)R_x,\quad R(0,\cdot)=F_{\rho_0}(\cdot)\,.
\en 
\end{cor}
\noindent\textbf{Proof.} \red{Setting $J:=J_{\ieta,\rho_0}$,
recall Remark \ref{rmk:ieta-small} that}
the exponentially tight $\{\rho^N\}$ satisfy 
the \abbr{LDP} upper bound in the metric space $(\CC,d)$ 
with some rate function $I_{\ieta,\rho_0} (\cdot) \ge J(\cdot)$.
%
\red{Necessarily, $I_{\ieta,\rho_0}$ has compact, non-empty, level sets and
in particular 
$J^{-1}(0) = \{ \gamma \in \CC : J(\gamma)=0 \}$ is non-empty (and pre-compact).
Considering the \abbr{LDP} upper bound for the complement 
of any finite $\delta$-cover of $J^{-1}(0)$, we further deduce by the first
Borel-Cantelli lemma that a.s.} 
$d(\rho^N,J^{-1}(0)) \to 0$. From the explicit formula \eqref{eq:J-def}  
we know that $J^{-1}(0) \subseteq \red{\cA_{\ieta,\rho_0}}$ 
and further to each $\gamma \in J^{-1}(0)$ 
corresponds $\Rg \in \FF$ which is a non-negative continuous bounded generalized 
solution of the problem \eqref{PMEcauchy}, in the sense of \cite[Definition 4]{gi}.
Recall \cite[Theorem 4]{gi} that such generalized solution is unique.
Consequently, $J^{-1}(0)=\{\gamma_{\star}\}$ is a single point to which 
$\rho^N$ converges a.s. when $N \to \infty$.
\ep

\smallskip
In \cite{dg}, the authors prove a \abbr{LDP} for systems of diffusions
with the same constant diffusion coefficient, where further 
each drift coefficient is a continuous function of the value of the diffusion and the 
empirical measure of the whole system. In this context
the (local) \abbr{LDP} is established by a clever application of Girsanov's theorem, 
which allows one to move from the system of interacting diffusions to the corresponding 
system of independent diffusions (on the event that the path of empirical measures 
is near a deterministic path of probability measures). Such 
an approach is not viable in our case because 
of the interaction through the diffusion coefficients in \eqref{sde}. Moreover, the discontinuity of the drift and the diffusion coefficients presents an additional challenge. Even on the level of the \abbr{lln} as in Corollary \ref{LLN}, previous works had to assume that there is no interaction through the diffusions coefficients (see 
\cite{bt1,bt2,jo} and the references therein), or be restricted to special initial conditions (see \cite{sh}). We overcome these challenges, but remark that our analysis relies on the particular form of the drift and the diffusion coefficients in \eqref{sde}. 

A crucial part of the proof of Theorem \ref{main_thm} is devoted to the 
study of generalized solutions to porous medium equations with tilt:
\eq\label{PDEtilt}
R_t=(A(R)R_x)_x+h\,A(R)\,R_x.
\en
The following regularity result, which we need in the proof of Theorem \ref{main_thm}, is also of independent interest.
\begin{thm}\label{main_pde_thm}
Let $R \in C_b(\rr_T)$ be 
such that, for every $t\in[0,T]$, the function $R(t,\cdot)$ is a \abbr{cdf} of a probability measure $\gamma(t)$. Suppose that $R$ is a generalized solution to \eqref{PDEtilt} with initial condition 
$R(0,\cdot)=F_{\mu}(\cdot)$, 
where $A(\cdot)$ satisfies Asmp. \ref{main_ass}(a),
$\mu \in \Mstar$ 
and $h$ is a function on $\rr_T$ such that 
\eq\label{eq:h-l2}
\int_{\rr_T} h^2(t,x)\,\mathrm{d}\gamma(t) \,\mathrm{d}t<\infty.
\en
If, in addition, $\gamma(\cdot)$ 
satisfies the moment condition (\ref{eq:eta-mom}) 
for some $\iota > 0$, then 
$\gamma \in \cA_{\iota,\mu}$ and
$R \in \FF_{q}$ 
of \eqref{eq:thm-13}, 
for all $\frac{6}{5}\leq q\leq \frac{3}{2}$. 
\end{thm}

\section{Outline: proofs of Theorems \ref{main_thm} and \ref{main_pde_thm}}\label{sec_outl}

In this section we establish Theorems \ref{main_thm} and \ref{main_pde_thm} as consequences of Propositions
\ref{L3prop} and \ref{approx_prop_main}--\ref{local_lbd}. The latter are in turn proved in the following five sections, in the order in which they are stated here. We note in passing that the proofs of Propositions \ref{L3prop} and \ref{approx_prop_main} are of analytic nature (relying for proving Prop. \ref{L3prop}
on results from \cite{KR,KR2,lsu} about 
parabolic equations with non-smooth coefficients),
whereas those of Propositions \ref{local_ubd_exp_tight} and \ref{local_lbd} are mostly probabilistic, involving tools from large deviations theory and stochastic analysis. More precisely, the local large deviations upper bound of Prop. \ref{local_ubd_exp_tight} is established by integrating suitable test functions against $\rho^N$, proving a Freidlin-Wentzell type local large deviations upper bound for the resulting processes and optimizing over such test functions; and the local large deviations lower bound of Prop. \ref{local_lbd} is shown via a tilting argument which relies on an appropriate Girsanov change of measure.  

We proceed with few notations and definitions
that are used throughout this paper. First, 
{we write 
$m$ for the Lebesgue measure on $\rr_T=[0,T]\times\rr$,}
$(\alpha,f)$ for $\int_\rr f\,\mathrm{d}\alpha$, any $f$ in the space of continuous 
bounded functions $C_b(\rr)$ and any $\alpha\in M_1(\rr)$, 
with $(\alpha,f)(s) = \int_\rr f(s,\cdot) \, \mathrm{d}\alpha(s,\cdot)$, in case
of $s \mapsto f(s,\cdot) \in C_b(\rr)$ and $s \mapsto \alpha(s,\cdot) \in M_1(\rr)$
(or more generally, whenever $f(s,\cdot)$ is integrable with respect to 
$\alpha(s,\cdot)$). We further let $\DxS = \{g_x : 
g \in \oS \}$ denote the space of spatial derivatives 
of test functions from $\oS$.
\begin{definition}
Setting  
\eq\label{eq:rgamma-def}
{\RR}^\gamma g=g_t+b(\Rg)g_x+A(\Rg)g_{xx} \,,
\en
{with} the functional
\eq\label{eq:phi-def}
\Phi_{\gamma} (t,g)
= (\gamma,g)(t)-(\gamma,g)(0)-\int_0^t (\gamma,{\RR}^\gamma g)(s)\,\mathrm{d}s
\en
on $\oS$, $\Phi_{\gamma} (g):=\Phi_{\gamma} (T,g)${,} and inner product
\eq\label{eq:ip-def}
(f,g)_\gamma=\int_{\rr_T} f_x\,g_x\,A(\Rg)\,\mathrm{d}\gamma(t)\,\mathrm{d}t 
\en
on $\oS$, we consider the (rate) functions on $\CC$ {given by}
\eq\label{eq:tilde-I-def}
I_{\iota,\mu}(\gamma)=\left\{\begin{array}{ll}
\sup_{g\in\oS} \big[ \Phi_\gamma (g) - (g,g)_\gamma \big], & \gamma\in \cA_{\iota,\mu}\\
\infty,& \mathrm{otherwise}.
\end{array}\right. 
\en
\end{definition}
Theorem \ref{main_pde_thm} is 
a direct consequence of the following proposition
(which is also key to proving Theorem \ref{main_thm}).
\begin{prop}\label{L3prop}
Suppose Asmp. \ref{main_ass}(a) holds and the function 
$b(\cdot)$ is uniformly bounded.
If $I_{\iota,\mu}(\gamma)<\infty$
for some $\iota > 0$ and $\mu \in \Mstar$, then $\Rg \in \FF_q$ 
for all $\frac{6}{5}\leq q\leq \frac{3}{2}$. Namely, $\Rg=R$ such that: 
\begin{enumerate}[(A)]
\item $R_x\in L^3(\rr_T)$,
\item $R_t,\,R_{xx}\in L^q(\rr_T)$ for all $\frac{6}{5}\leq q\leq \frac{3}{2}$,
\item $\int_{\rr_T} \frac{R_{xx}^2}{R_x}\,\mathrm{d}m<\infty$, $\int_{\rr_T} \frac{R_t^2}{R_x}\,\mathrm{d}m<\infty$.
\end{enumerate}
\end{prop}

\noindent\textbf{Proof of Theorem \ref{main_pde_thm}.} 
After integration by parts in space, we see 
that having $R=\Rg$ as
generalized solution of \eqref{PDEtilt}, is 
equivalent to   
\eq\label{eq:PDE-dual}
\Phi_\gamma (t,g) 
= - \int_0^t (\gamma,b(R) g_x + h A(R) g_x)(s) 
\mathrm{d}s  \,,
\en
for any $g \in \oS$ and $t \in [0,T]$. In particular, upon 
comparing \eqref{eq:phi-def} and \eqref{eq:PDE-dual}, we deduce that 
$t \mapsto (\gamma,g)(t)$ is absolutely continuous for any
$g \in \oS$, and thus 
$\gamma \in \cA_{\iota,\mu}$ (for 
$R \in C_b(\rr_T)$, and 
the moment and initial conditions have all been assumed in
Theorem \ref{main_pde_thm}).  
Further, taking here $b \equiv 0$ without loss of generality, 
we see that for $\gamma$ as in Theorem \ref{main_pde_thm},
\begin{align}
I_{\iota,\mu}(\gamma) 
&=
\sup_{f\in\DxS} \Big[\int_{\rr_T} \big(-h\,A(R)f-A(R)f^2\big)\,\mathrm{d}\gamma(t)\,\mathrm{d}t\Big]
\nonumber \\
&\le 
\sup_{f\in L^2(\rr_T,\,\mathrm{d}\gamma(t)\mathrm{d}t)}\Big[\int_{\rr_T} \big(-h\,A(R)f-A(R)f^2\big)\,\mathrm{d}\gamma(t)\,\mathrm{d}t\Big].
\nonumber 
\end{align}
The latter supremum is attained for $f=-\frac{1}{2} h$ and its value is finite due to our assumption \eqref{eq:h-l2}.
Consequently, in this case $I_{\iota,\mu}(\gamma)<\infty$ and by Prop. \ref{L3prop} such $R=\Rg$ satisfies
the regularity properties \eqref{eq:thm-13} for
all $\frac{6}{5} \le q \le \frac{3}{2}$, as claimed.  
\ep

We start the proof of Theorem \ref{main_thm} by establishing the following corollary of
Prop. \ref{L3prop}.  
\begin{cor}\label{rate_fct_prop}
Suppose $A(\cdot)$ satisfies Asmp. \ref{main_ass}(a), 
$\mu \in \Mstar$ and $b(\cdot)$ is uniformly bounded.
If $I_{\iota,\mu}(\gamma)<\infty$ then 
$J_{\iota,\mu} (\gamma)=I_{\iota,\mu}(\gamma)$. 
In particular, 
$J_{\iota,\mu}(\gamma)\leq I_{\iota,\mu}(\gamma)$ 
for all $\gamma\in \CC$.
\end{cor}

\noindent\textbf{Proof.}
Fixing $\gamma\in \cA_{\iota,\mu}$  with 
\[
I_{\iota,\mu} 
(\gamma) = \sup_{g\in\oS} \big[\Phi_\gamma (g)-(g,g)_\gamma\big] <\infty\,,
\]
consider the Hilbert space $\mathbb{H}$ 
given by \red{identifying and completing 
$\oS$ under the semi-norm}
corresponding to the inner product 
$(\cdot,\cdot)_{\gamma}$ of (\ref{eq:ip-def}). 
\red{By scaling, the linear functional $\Phi_\gamma(\cdot)$ 
is bounded on $\oS$ 
by $2 \sqrt{I_{\iota,\mu}(\gamma)}$ times
this semi-norm, and with
$A(\cdot)$ uniformly bounded below, 
if
$(g,g)_{\gamma}=0$ then by
\eqref{eq:rgamma-def}-\eqref{eq:ip-def} also 
$\Phi_\gamma(g)=0$.
Hence, there exists a unique bounded linear functional $\overline{\Phi}_\gamma$ 
on $\mathbb{H}$ which} coincides with $\Phi_\gamma$ on $\oS$.
Now, by the Riesz representation theorem, there is a unique element $\tilde{h}\in\mathbb{H}$, which satisfies 
$\overline{\Phi}_\gamma (g)=(\tilde{h},g)_\gamma$ for all $g\in\mathbb{H}$. Combining this with the fact that $\oS$ is by definition dense in $\mathbb{H}$, we obtain that
\eq\label{Ifromh}
I_{\iota,\mu} (\gamma)=\sup_{g\in\mathbb{H}} 
\big[\overline{\Phi}_\gamma (g)-(g,g)_\gamma\big]
=\sup_{g\in\mathbb{H}} \big[(\tilde{h},g)_\gamma-(g,g)_\gamma\big]=\frac{1}{4}(\tilde{h},\tilde{h})_\gamma.
\en
Furthermore, by the definition of $\tilde{h}$ and
$\Phi_\gamma$, we have that $\tilde{h}_x \in L^2(\rr_T,
\mathrm{d}\gamma(t)\,\mathrm{d}t)$ satisfies
\eq\label{3.4}
\Phi_\gamma(t,g) = \int_0^t 
(\gamma,A(\Rg) \tilde{h}_x \, g_x) (s) \,\mathrm{d}s
\en
for $t=T$ and any $g\in\oS$. 
In particular, considering Schwartz functions $g$ 
supported on $\rr_t$ we have that \eqref{3.4} 
applies also for any $t \in [0,T]$. Comparing 
this with \eqref{eq:PDE-dual} we deduce that
$R=\Rg$ is a generalized solution of the \abbr{PDE}
\eqref{PDEtilt} for 
\eq\label{htildeh}
h=-\tilde{h}_x-\frac{b(R)}{A(R)}\,. 
\en
By the assumed boundedness of $\frac{b}{A}$, clearly 
$h \in L^2(\rr_T,\mathrm{d}\gamma(t)\,\mathrm{d}t)$. By 
Theorem \ref{main_pde_thm}, this implies in turn that $R_t$, $R_{xx}$ and the
$L^1(\rr_T)$ density $R_x$ are 
elements of $L^{3/2}(\rr_T)$ and, moreover,
the functions $R_t R_x^{-1/2}$, $R_x^{3/2}$ and $R_{xx} R_x^{-1/2}$ 
are elements of $L^2(\rr_T)$. 
Thus, the identity
\eq\label{whatish}
h (A(R) R_x)^{1/2} =\frac{R_t-(A(R)R_x)_x}{(A(R)R_x)^{1/2}}
\en 
holds in $L^2(\rr_T)$. Finally, putting \eqref{whatish}, \eqref{htildeh} and \eqref{Ifromh} together, we end up with
\eqs
I_{\iota,\mu} 
(\gamma)=\frac{1}{4} \|\widetilde{h}_x (A(\Rg)\Rgx)^{1/2} \|_{L^2(\rr_T)} 
= J_{\iota,\mu}(\gamma)
\ens
of \eqref{eq:J-def}, as claimed. 
\ep

Corollary \ref{rate_fct_prop}
allows us to replace the function 
$I_{1,\rho_0}(\cdot)$ of the large deviations
upper bound for Theorem \ref{main_thm}
(see \eqref{ld-local-ubd}), by $J_{1,\rho_0}(\cdot)$ 
of the corresponding lower bound
(see Prop. \ref{local_lbd}). 
The task of proving such lower bound is further simplified 
thanks to the next proposition, for which we first introduce  
some relevant notations.
\begin{definition}\label{new-def-eta}
Let $\GG$ denote the subset 
of $\{ \gamma \in \CC : \widetilde{J} (\gamma) < \infty \}$ 
for which $R := \Rg \in C^{\infty}_b(\rre)$, $R_x$ is 
strictly positive and 
\eqref{PDEtilt} holds point-wise for some
$h \in C_b(\rre)$ with $x \mapsto h(t,x)$ uniformly Lipschitz continuous on $\rre$.
\end{definition}

\begin{prop}\label{approx_prop_main}
Suppose Asmp. \ref{main_ass}(a),(b) and (d) hold, and that 
$J_{\red{1},\mu}(\gamma)<\infty$ for 
some $\mu \in M_1^{(1)}(\rr)$. Then, {there exist}
$\gamma^{\ell,\eps} = (1-\ell^{-1}) \gamma^\eps + \ell^{-1} \widehat{\mu}^\eps \in \GG$
for some $\widehat{\mu}^\eps \in \Mstar$ and $\gamma^\eps \in \CC$ such that 
\begin{align}\label{eq:iota-ubdd}
\sup_{\eps} \int_\rr |x|^{2} \mathrm{d}\gamma^\eps(0) &< \infty \,,\\
\label{eq:d-eta-conv}
\lim_{\eps \to 0} \limsup_{\ell \to \infty} \;d(\gamma^{\ell,\eps},\gamma) &= 0\,,\\
\label{rate_control}
\limsup_{\eps \to 0} \limsup_{\ell \to \infty} 
\; {\widetilde J} \, (\gamma^{\ell,\eps}) 
&= J_{\red{1},\mu} (\gamma) < \infty\,.
\end{align}
\end{prop}

We proceed to state our basic local large deviations bounds.
\begin{prop}\label{local_ubd_exp_tight}
With $B(\gamma,\delta)$ denoting the open ball of radius $\delta>0$ centered at 
arbitrary $\gamma \in \CC$, under Asmp. \ref{main_ass}(a)-(c)
we have the local large deviations upper bound
\eq\label{ld-local-ubd}
\lim_{\delta\downarrow0}
\limsup_{N\rightarrow\infty}\frac{1}{N}
\log \pp(\rho^N\in B(\gamma,\delta))
\leq -I_{\red{\ieta,\rho_0}} (\gamma) \,.
\en
Moreover, the sequence $\{\rho^N$, $N\in\nn\}$ is exponentially tight 
in the sense that for any $M<\infty$ there exists a compact set 
$K_M\subset \CC$, for which
\eq\label{eq:ex-tight}
\limsup_{N\rightarrow\infty}\frac{1}{N}\log \pp(\rho^N\notin K_M)\leq -M.
\en
\end{prop}

\begin{prop}\label{local_lbd} 
Under Asmp. \ref{main_ass}(a)-(c) with $\ieta=1$ and 
$\{\gamma^{\ell,\eps}\} \subset \GG$ of Prop. \ref{approx_prop_main}: 
\newline
(a) Asmp.~\ref{main_ass}(c) holds for $\rho_0^{\ell,\eps}=\gamma^{\ell,\eps}(0)$,
$\ieta=0$ and 
deterministic initial empirical measures $\{\rho^{N,\ell,\eps}(0), N \in \nn\}$ 
such that 
\begin{equation}\label{eq:coupling}
\limsup_{\eps \to 0} \limsup_{\ell \to \infty}
\limsup_{N \to \infty} \, d(\rho^N,\rho^{N,\ell,\eps}) = 0 \,.
\end{equation}
(b) The corresponding local large deviations lower bound
\eq\label{eq:ldp-loc-lbd}
\liminf_{N\rightarrow\infty}\frac{1}{N}\log \pp(\rho^{N,\ell,\eps}\in 
B(\gamma^{\ell,\eps},\delta)) \geq - {\widetilde J} (\gamma^{\ell,\eps})
\en
holds for any $\ell$, $\eps,\delta>0$. 
\end{prop}

\noindent\textbf{Proof of Theorem \ref{main_thm}.} From \cite[Theorem 4.1.11]{dz} and \cite[Lemma 1.2.18]{dz} we conclude that Theorem \ref{main_thm} follows once we show
that for \red{$J=J_{1,\rho_0}$} and any $\gamma\in  \CC$,
\begin{align}
\lim_{\delta\downarrow0}\limsup_{N\rightarrow\infty}\frac{1}{N}\log \pp(\rho^N\in B(\gamma,\delta))&\leq -J(\gamma)\,,\;\;\;
\label{UBD1}\\
\lim_{\delta\downarrow0}\liminf_{N\rightarrow\infty}\frac{1}{N}\log \pp(\rho^N\in B(\gamma,\delta))&\geq -J(\gamma)\,, \label{LBD1}
\end{align}  
and that the sequence $\{\rho^N$, $N\in\nn\}$ is exponentially tight in the sense of (\ref{eq:ex-tight}).
To this end, note that $\{\rho^N$, $N\in\nn\}$ is exponentially tight 
(by Prop.\ref{local_ubd_exp_tight}), and 
the upper bound \eqref{UBD1} is a consequence of Prop. \ref{local_ubd_exp_tight} and Corollary \ref{rate_fct_prop}. Next, it clearly suffices to establish the lower bound \eqref{LBD1} when  
$J(\gamma)<\infty$. To this end, 
fixing such $\gamma$ and 
$\gamma^{\ell,\eps} \in \GG$ as in Prop. \ref{approx_prop_main}
\red{(where $\mu=\rho_0$)},
we know from part (a) of Prop. \ref{local_lbd} and \eqref{eq:d-eta-conv}
that, 
for any $\delta>0$, $\eps \le \eps_0(\del)$, $\ell \ge \ell_0(\eps,\del)$,  
if  
$\rho^{N,\ell,\eps} \in B(\gamma^{\ell,\eps},\delta)${,} then 
$\rho^N \in B(\gamma,3\delta)$ for all $N$ large enough. 
From part (b) of Prop. \ref{local_lbd} we deduce that then  
\eqs
\liminf_{N\rightarrow\infty}\frac{1}{N}\log \pp(\rho^N\in B(\gamma,3\delta))\geq - 
{\widetilde J} (\gamma^{\ell,\eps}) \,,
\ens
so using \eqref{rate_control} we get \eqref{LBD1} upon
taking $\ell \to \infty$, then $\eps \downarrow 0${,} and finally
$\delta \downarrow 0$.
\ep

\section{Proof of Proposition \ref{L3prop} (A) and (B)}\label{sec_L3}

Throughout this section, the function $A(\cdot)$ satisfies Asmp. \ref{main_ass}(a),
and the function $b(\cdot)$ is uniformly bounded. Then, fixing $\iota > 0$, 
$\mu \in M_1(\rr)${,} and $R=\Rg$
for $\gamma \in \CC$ with $I_{\iota,\mu}(\gamma)<\infty$,
we prove Prop. \ref{L3prop} (A) in a series of
three lemmas, starting in Lemma \ref{lem-l2rx} 
by showing that $R_x \in L^{3/2}(\rr_T)$, 
which we improve in Lemma \ref{lem-lprx} to $L^p$ 
estimates on $R_x$ for all $\frac{3}{2} \leq p<3$. 
Finally, Lemma \ref{lem-l3rx} establishes the uniform boundedness of the 
corresponding norms when $\mu \in \Mstar$,
resulting with $R_x\in L^3(\rr_T)$. 
\begin{lemma}\label{lem-l2rx}
If $R=\Rg$ and $I_{\iota,\mu}(\gamma)=I<\infty\;$ for some $\mu \in M_1(\rr)$ then:
\begin{enumerate}[(a)]
\item The function $a:=A(R)$ is uniformly continuous on $\rr_T$.
\item The measure $\mathrm{d}\gamma(t)\,\mathrm{d}t$ on $\rr_T$ has a density 
with respect to the Lebesgue measure on $\rr_T$, whose $L^2$ norm restricted to any 
strip $S_{n,r}:=[0,T]\times[n-\frac{r}{4},n+\frac{r}{4}]$ is bounded by a constant 
$C(T,I,r)<\infty$ (independent of $n\in\zz$). 
In particular, this density is locally square integrable, so that the 
weak derivative $R_x$ exists as a function in $L^2_{loc}(\rr_T)$. 
\item The weak derivative $R_x$ is an element of 
$L^{3/2}(\rr_T)$.
\end{enumerate}
\end{lemma}

\noindent\textbf{Proof.} 
(a). Recall that $A(\cdot)$ is assumed to be Lipschitz, 
so it suffices to show uniform continuity
of $R=\Rg$ on $\rr_T$. We further
assumed that $\gamma \in \cA_{\iota,\mu}$, hence
$R \in C_b(\rr_T)$ is uniformly continuous 
on compact sets. In addition, the continuity of 
$t\mapsto\gamma(t)$ with respect to the topology of
weak convergence in $M_1(\rr)$ implies that the 
image $\{\gamma(t)\}_{t\in[0,T]}$ 
of the compact $[0,T]$ is compact in $M_1(\rr)$, 
hence by Prokhorov's Theorem, uniformly tight. 
Consequently, for every $\alpha>0$, there 
exists finite $M=M_\alpha$ such that
\eq\label{eq:Runif-tight}
\sup_{t\in[0,T],x\geq M}\max(R(t,-x),1-R(t,x)) 
< \alpha \,,
\en
extending the uniform continuity of $R$ to all of 
$\rr_T$. 

\smallskip
\noindent (b). Setting hereafter $S_n=S_{n,2}$,  
it suffices to show that the inequality
\eq\label{L2locbound}
\int_{S_n} \psi(t,x)\,\mathrm{d}\gamma(t)\,\mathrm{d}t\leq C(T,I)\|\psi\|_2\,,
\en
holds for some constant $C(T,I)<\infty$, 
all $n\in\zz${,} and any continuous
$\psi: S_n \rightarrow\rr_+$. Indeed, this implies 
the existence of the Radon-Nikodym derivative 
$R_x=\mathrm{d}\gamma(t)\mathrm{d}t/\mathrm{d}m$
on $\rr_T$, whose $L^2(S_n)$-norm is bounded by 
$C(T,I)$ (by the same argument we used 
en-route to \eqref{Ifromh}). Turning to prove 
\eqref{L2locbound}, by 
definition of $I_{\iota,\mu}(\gamma)$ and the identity
$\inf_{\lambda>0} \{ \lambda y^2 + \lambda^{-1} z^2 \}
= 2 |y z|$, we have that for
$C_1=2 (T I \|A\|_\infty)^{1/2}$ finite,
$a(t,x):=A(R(t,x))${,} and any $g\in\oS$,
\eq\label{bd0}
C_1 \|g_x\|_{\infty} \ge 
2\,I^{1/2} \left(\int_{\rr_T} (g_x)^2 a(t,x)\,\mathrm{d}\gamma(t)\,\mathrm{d}t\right)^{1/2} 
\ge \Phi_\gamma (g) \,.
\en 
Considering first $n=0$, we use \cite[Theorem 2]{KR}
(taking there $d=1$ and $f(t,x):=\psi(T-t,x)$). It provides
a universal finite constant $C_2$ and space-time 
kernels $k^\epsilon (t,x)=\epsilon^{-2} \xi(t/\epsilon) 
\xi(x/\epsilon)$, $\epsilon>0$, for some infinitely differentiable, 
probability density $\xi(\cdot)$ of compact support, 
with the following property:
\newline 
For any continuous $\psi:S_0 \mapsto \rr_+$ 
there exists bounded measurable $z:\,\rr_T \mapsto (-\infty,0]$ 
(depending on $\psi$), non-decreasing in $t$ and supported 
on larger strip $S_{0,4}$, \textit{within which} it is 
convex in $x$, such that for any $\epsilon>0$ small 
enough, the smooth functions $\psi^\epsilon=\psi \ast k^\epsilon$ 
and $z^\epsilon=z\ast k^\epsilon$ satisfy 
on the intermediate strip $S_{0,3}$ the inequalities
\begin{align}
\forall c\geq0:\quad &c^{1/2} \psi^\epsilon\leq C_2 \left(z^\epsilon_t+c\,z^\epsilon_{xx}\right), \label{bd1}\\
&\frac{1}{2}|z^\epsilon_x| \leq -z^\epsilon \leq C_2 \,\|\psi\|_2\,. \label{bd2}
\end{align}
In the preceding, the compact support of $z$ 
is specified in the proof of \cite[Theorem 2]{KR} 
after \cite[display (29)]{KR}.

\smallskip
Picking a $[0,1]$-valued truncation function $\phi \in\oS$, supported on $S_{0,3}$ 
with $\phi \equiv 1$ on $S_{0,2}$, we note that 
for each $\epsilon >0$ as above, the non-negative 
function $g= -z^{\epsilon} \phi$ is in $\oS$, supported on $S_{0,3}$ such that by \eqref{bd2} 
both $\|g\|_\infty$ and $\|g_x\|_{\infty}$ are 
bounded by $C_3 \|\psi\|_2$ (for the universal 
constant $C_3 = (2+\|\phi_x\|_{\infty}) C_2$).  
Consequently, applying the bound \eqref{bd0} for 
such choice of $g$, we deduce that
$$
C_1 C_3 \|\psi\|_2 
\ge \Phi_\gamma (g) 
\ge - \int_{S_{0,3}} (g_t + a g_{xx}) 
\mathrm{d}\gamma(t) \mathrm{d}t
- T \|b(R)\|_{\infty} \|g_x\|_{\infty} 
- \|g(0,\cdot)\|_{\infty} \,.
$$
Next with $\phi_t$, $\phi_x$ and $\phi_{xx}$ uniformly bounded
by some universal constant $C_4$, it follows from
Leibniz's rule 
and \eqref{bd2} that 
$$
|(g_t + a g_{xx}) + (z^\epsilon_t + a z^\epsilon_{xx}) \phi|
= |z^\epsilon (\phi_t + a \phi_{xx}) + 2 a \phi_x z^\epsilon_x|  
\le C_2 C_4 (1 + 5 \|a\|_\infty) \|\psi\|_2 \,,
$$
out of which we deduce by simple algebra that 
\eq\label{bd99}
C_5 \|\psi\|_2 \ge 
\int_{S_{0,3}} (z^{\epsilon}_t + a z^{\epsilon}_{xx}) 
\phi \mathrm{d}\gamma(t) \mathrm{d}t 
\en
for finite $C_5$ depending only on 
$T$, $\|b\|_{\infty}$, $\|A\|_{\infty}$ and the 
constants $C_i$, $i=1,\ldots,4$. With $z$ 
non-decreasing in $t$, so are 
$z^\epsilon = z \ast k^\epsilon$.
Further, as $\xi(\cdot)$ has compact support, 
the convexity of $z$ in $x$ within $S_{0,4}$ implies 
that $z^\epsilon$ is convex in $x$ on $S_{0,3}$, 
provided $\epsilon > 0$ is small enough. 
Thus, 
both $z^\epsilon_t$ and $z^\epsilon_{xx}$ are non-negative,
so considering \eqref{bd1} for the strictly positive $c=\underline{a}$
of Asmp. \ref{main_ass}(a) and recalling 
that $\phi \equiv 1$ on $S_0=S_{0,2}$ (with $\phi \ge 0$ 
elsewhere), we have from \eqref{bd99} that
\[
\begin{split}
C_2 C_5 \|\psi\|_2 \ge 
C_2 \int_{S_0} (z^{\epsilon}_t + a z^{\epsilon}_{xx}) 
\phi \mathrm{d}\gamma(t)\,\mathrm{d}t  
\ge 
\int_{S_0} C_2 (z^{\epsilon}_t + c z^{\epsilon}_{xx}) 
\mathrm{d}\gamma(t)\,\mathrm{d}t \\
\ge c^{1/2} \int_{S_0} \psi^\epsilon \,\mathrm{d}\gamma(t)\,\mathrm{d}t \,.
\end{split}
\]
With $\psi \in C_b(S_0)$, clearly  
$\psi^\epsilon \to \psi$ uniformly on $S_0$
when $\epsilon \downarrow 0$, leading to
\eqref{L2locbound} for $n=0$ and the
universal finite constant 
$C(T,I)=C_2 C_5 c^{-1/2}$ which depends only 
on $T$, $I${,} and the functions $b(\cdot)$ 
and $A(\cdot)$ on $[0,1]$.  

\smallskip
To extend \eqref{L2locbound} to other values of $n \in \zz$, let 
$\mu_n(\cdot) := \mu(\cdot+n)$ and note that the path 
$\gamma_n(\cdot,\cdot)=\gamma(\cdot,\cdot+n) \in
\cA_{\iota,\mu_n}$ if and only if $\gamma \in \cA_{\iota,\mu}$. Setting next
$g_n(\cdot,\cdot)=g(\cdot,\cdot+n)$ we clearly 
have that $(\gamma_n,{\RR}^{\gamma_n} g_n)(s)$
is independent of $n$ for each $s \in [0,T]$. 
Hence,  $I_{\iota,\mu_n}(\gamma_n) = I$  
and{,} since to any non-negative $\psi \in C_b(S_n)$ 
corresponds non-negative $\psi_n(\cdot,\cdot)=\psi(\cdot,\cdot+n) \in C_b(S_0)$, by the preceding proof:
\[
\int_{S_n} \psi(t,x)\,\mathrm{d}\gamma(t)\,\mathrm{d}t = \int_{S_0} \psi_n(t,x)\,\,\mathrm{d}\gamma_n(t)\,\mathrm{d}t 
\leq C(T,I)\|\psi_n\|_2 = C(T,I)\|\psi\|_2 
\]
(as $C$ is independent of the initial measure). This completes the proof of part (b).

\smallskip
\noindent (c). Fixing $n\in\zz$ we already know 
that $R_x \in L^p(S_n)$ for $p=1,2$. Further, 
upon applying Cauchy-Schwarz inequality with respect 
to Lebesgue measure on $S_n$, we have that
$$
\int_{S_n} R_x^{3/2}\,\mathrm{d}m \leq 
\Big(\int_{S_n} R_x\,\mathrm{d}m\Big)^{1/2}
\Big(\int_{S_n} R_x^2\,\mathrm{d}m\Big)^{1/2}\,.
$$
We have shown in part (b) that the 
right-most term is bounded uniformly in $n$,
so our claim that the left-side is summable over
$n \in \zz$ follows from finiteness of 
$\sum_{|n| \ge 1} \kappa_n$, with 
$\kappa_n := \Big(\int_{S_n} R_x\,\mathrm{d}m\Big)^{1/2}$.
Next{,} taking $\iota > 0$ as in the given moment
condition \eqref{eq:eta-mom}, we get by
Cauchy-Schwarz that for 
$c_1 := \sum_{|n| \ge 1} |n|^{-(1+\iota)}$ and 
$c_2 := \sup \{ |n/y|^{1+\iota} : y \in S_n, |n| \ge 1 \}$ 
finite, 
\begin{align*}
\big(\sum_{|n| \ge 1} \kappa_n \big)^2  
\le c_1 \sum_{|n| \ge 1} |n|^{1+\iota} \kappa_n^2  
&\leq c_1 c_2 \sum_{|n| \ge 1} 
\int_{S_n} |y|^{1+\iota} R_x(t,y) \,\mathrm{d}m(t,y) \\
&\le c_1 c_2 \int_{\rr_T} |y|^{1+\iota} \mathrm{d}\gamma(t)\,\mathrm{d}t < \infty
\end{align*}
(see \eqref{eq:eta-mom} for the right-most inequality).
\ep

\begin{lemma}\label{lem-lprx}
If $R=\Rg$
and $I_{\iota,\mu}(\gamma)<\infty$, then  
$R_x$ exists as element of $L^p(\rr_T)$ 
for all $\frac{3}{2}\leq p<3$. 
\end{lemma}

\noindent\textbf{Proof.} We fix 
$\frac{3}{2}\leq p<3$ and set $\frac{3}{2}<q\leq3$ 
such that $\frac{1}{p}+\frac{1}{q}=1$. For any 
given non-negative $f\in\oS$, 
we consider the backward Cauchy problem
\eq\label{eq-cauchy}
u_t+a(t,x)\,u_{xx}+f(t,x)\,u=0, \quad u(T,\cdot)=1
\en
with $a(t,x):=A(R(t,x))$. While $a(t,x) \in C_b(\rr_T)$ 
is not differentiable in $x$, recall part (a) of Lemma 
\ref{lem-l2rx} that it is uniformly continuous and
bounded away from zero. Considering \cite[Theorem 2.1]{KR2}
for $w:=u-1$, such $a(\cdot,\cdot)$ and 
any $f \in L^r(\rr_T)$, we see that \eqref{eq-cauchy} 
has a weak solution, for which $w=u-1$ is an element of $W^{1,2}_r(\rr_T)$. 
Further, here  $f\in L^q(\rr_T)\cap L^2(\rr_T) \cap L^3(\rr_T)$ from which 
it follows that the norm bounds of \cite[Theorem 4.1]{KR2} can be refined to an 
estimate on the norms in $L^q(\rr_T)\cap L^2(\rr_T)\cap L^3(\rr_T)$. From the
latter estimate, we conclude upon applying the method of continuity 
(see e.g. \cite[section III.2]{Li}), that for $f \in \oS$
the problem \eqref{eq-cauchy} has a weak solution, for which $u-1$ is an 
element of $W^{1,2}_q(\rr_T)\cap W^{1,2}_2(\rr_T)\cap W^{1,2}_3(\rr_T)$. Next, let 
$Y(\cdot)$ denotes the canonical process having the law $P$ 
of a diffusion with generator $a(t,x)\,\frac{\mathrm{d}^2}{\mathrm{d}x^2}$ 
(which exists thanks to uniform continuity and boundedness from above and below of $a$).
Applying It\^o's formula in the form of the last identity in \cite[chapter 10, Theorem 1]{kr} (again using the boundedness from above and below of $a$), we obtain the stochastic representation of such weak solution of \eqref{eq-cauchy},
\eq
u(t,x)=\ev^P\Big[\exp\Big(\int_t^T f(s,Y(s))\,\mathrm{d}s\Big)\Big|Y(t)=x\Big]\,.
\en 
The latter 
representation implies that $u\geq 1$, and by Portenko's Lemma (see \cite[inequality (6)]{Po}, note that it only relies on the standard heat kernel estimate for the diffusion with law $P$ and that $q>\frac{3}{2}$ throughout this proof), there is a non-decreasing 
function $G_q:\,[0,\infty)\rightarrow [1,\infty)$ depending only on $\frac{3}{2}<q\leq3$ (and \textit{not on $f$}), 
such that
$$
1 \leq u(t,x) \leq G_q \big(\|f\|_{L^q(\rr_T)}\big),\quad (t,x)\in\rr_T.
$$
Next, observe that the non-negative $\v:=\log u$ inherits the bound
\eq\label{vfnormbound}
0 \le \v (t,x) \leq 
G_q \big(\|f\|_{L^q(\rr_T)}\big),\quad (t,x)\in\rr_T \,.
\en
Recall that $u_t,\,u_{xx} \in L^3(\rr_T)$,
and with $u \in W^{1,2}_2(\rr_T)$ further
$u_x \in L^6(\rr_T)$ 
(by the parabolic Sobolev inequality in the form of \cite[chapter II, Lemma 3.3]{lsu}),
which for positive $u$ bounded away from zero, imply also that 
$\v_t,\,\v_{xx}\in L^3(\rr_T)$ and $\v_x\in L^6(\rr_T)$ (with $\v_{xx}+\v_x^2=u_{xx}/u$).
From this it follows in turn that 
$m$-a.e. on $\rr_T$ such {$\v$} 
and its generalized derivatives,
satisfy the backward equation
\eq\label{vPDEL3}
\v_t+a(t,x)\,\v_{xx}+a(t,x)\,\v^2_x+f(t,x)=0, \quad \v(T,\cdot)=0 \,.
\en 
Recall that while proving Corollary \ref{rate_fct_prop} we found in 
\eqref{htildeh}
a function $h$ whose
$L^2(\rr_T, a R_x \mathrm{d}m)$-norm 
is bounded by $C_6 := 2 I^{1/2} + \|b\|_\infty (T/\underline{a})^{1/2}$
(due to \eqref{Ifromh} and the assumed
bounds on $b(\cdot)$ and $A(\cdot) \ge \underline{a}$), {so} that{,}
for any $g \in \oS$, 
\eq\label{g-chain}
(\gamma,g)(T)-(\gamma,g)(0) = 
\int_{\rr_T} (g_t + a g_{xx} - a g_x h) R_x \mathrm{d} m 
\en
(consider \eqref{3.4}). With $R_x \in L^{3/2}(\rr_T)$ (by part (c) of Lemma
\ref{lem-l2rx}), and H\"older's inequality
$$
\int_{\rr_T} |\varphi| R_x \mathrm{d} m
\le 
\Big( \int R_x^{3/2} \mathrm{d}m\Big)^{2/3}
\Big( \int |\varphi|^3 \mathrm{d} m\Big)^{1/3} \,,
$$
we deduce that $\v_t$, $a \v_{xx}$ and 
$a \v_x^2$ are integrable with respect 
to $R_x \mathrm{d}m$ (as is $a h^2$).
Thus, with $\v$ bounded, we get from 
\eqref{g-chain} upon approximating
$\v$ in suitable (mixed)-norm by 
functions $g_k \in \oS$ 
that
$$
(\gamma,{\v})(T)
-(\gamma,{\v})(0)
=\int_{\rr_T} (\v_t+a \v_{xx}-a \v_x h)\,R_x\,\mathrm{d}m \,.
$$
The latter identity, in combination with 
\eqref{vPDEL3}, 
results with 
\begin{align*}
\int_{\rr_T} f R_x \mathrm{d} m
&= - \int_{\rr_T} (\v_x h + \v^2_x) a \,R_x\,\mathrm{d}m
+ (\gamma,{\v})(0) \\
&\leq \frac{1}{4} \int_{\rr_T} h^2 a \,R_x\,\mathrm{d}m + 
\sup_{(t,x) \in \rr_T} {\v}(t,x) 
\le \frac{1}{4} C_6^2 +  
G_q \big(\|f\|_{L^q(\rr_T)}\big) \,, 
\end{align*}
where the first inequality is merely the
non-negativity of 
$L^2(\rr_T, a R_x \mathrm{d}m)$-norm
of $\v_x+h/2$, and the second follows 
from \eqref{vfnormbound}. With 
$C_6$ independent of $f$, we have that
the linear functional 
$f \mapsto \int_{\rr_T}
f R_x \mathrm{d} m$ is bounded
on $\oS$ with 
respect to the $L^q(\rr_T)$-norm, 
hence $R_x\in L^p(\rr_T)$ for $1/p=1-1/q$,
as claimed. \ep

\begin{lemma}\label{lem-l3rx}
Suppose $R=\Rg$
and $I_{\iota,\mu}(\gamma)=I<\infty$ for
some $\mu \in \Mstar$. Then
$R_x \in L^3(\rr_T)$.
\end{lemma} 

\noindent\textbf{Proof.}
Recall {from} Lemma \ref{lem-lprx} that 
$R_x \in L^p(\rr_T)$ for any $p \in [3/2,3)$. Hence,
$R_x \in L^3(\rr_T)$, provided 
$p \mapsto \|R_x\|_{L^p(\rr_T)}$ is uniformly bounded 
over such $p$, as we prove here.

\smallskip
\noindent\textbf{Step 1.} Recall that $a(t,x)=A(R)$ and
$\int_{\rr_T} (g_t R_x + g_{xt} R) \mathrm{d} m = 0$ 
(integration by parts in $x$). When substituted in 
\eqref{g-chain}, these yield that
the path of \abbr{cdf}-s 
$R=\Rg$
is a generalized solution of 
the Cauchy problem
\eq\label{R_cauchy}
R_t-(A(R)R_x)_x=h\,A(R)R_x\,,\quad R(0,\cdot)=F_{\mu}
\en
with respect to the collection $\DxS := \{g_x :  g \in \oS\}$ of test functions, for 
some function $h$ whose $L^2(\rr_T,aR_x \mathrm{d} m)$-norm
is bounded by finite $C_6=C_6(T,I)$. That is,
\eq\label{gen-sol}
\int_{\rr} (g_x R) (T,x) \mathrm{d} x 
- \int_{\rr} (g_x R) (0,x) \mathrm{d} x =
\int_{\rr_T} (g_{xt} R - g_{xx} a R_x + g_x h a R_x) 
\mathrm{d} m \,,
\en
for any $g \in \oS$. 

Next, applying for $f=A(R) R_x$ the 
H\"older's inequality in the form of 
\eq\label{Holder_bnd} 
\int_{\rr_T} |h|^q\,f^q\,\mathrm{d}m \leq \Big( \int_{\rr_T} |h|^2\,f\,\mathrm{d}m \Big)^{q/2} \Big(\int_{\rr_T} f^p\,\mathrm{d}m \Big)^{q/(2p)}
\en
with $\frac{3}{2}\leq p<3$ and $\frac{q}{2}+\frac{q}{2p}=1$,
namely $q:=\frac{2p}{p+1}\in[6/5,3/2)$, we deduce that $h\,A(R)R_x\in L^q(\rr_T)$. Thus, by the regularity theory for 
the heat equation (see inequalities (3.1) and (3.2) in \cite[chapter IV]{lsu}, or \cite[Theorem 2.1]{KR2}), we have that the function
$$
V(t,x):=\int_0^t \int_\rr (h\,A(R)R_x)(s,y) \varphi(t-s,x-y)\,\mathrm{d}y\,\mathrm{d}s,
$$
obtained by convolving 
with the heat kernel $\varphi(t,x)=(4\pi t)^{-1/2} \exp(-x^2/4t)$, is a
generalized solution of the auxiliary Cauchy problem
\eq\label{V_cauchy}
V_t-V_{xx}=h\,A(R)R_x\,,\quad V(0,\cdot)=0
\en
(with respect to the collection $\DxS$). Further,
$$
\|V\|_{W^{1,2}_q(\rr_T)}\leq C_7 \|h\,A(R)R_x\|_{L^q(\rr_T)} 
$$   
where $C_7<\infty$ is a uniform constant (which in particular does not depend on $q$ as long as $q$ 
belongs to a compact interval). Thus{,} $V\in W^{1,2}_q(\rr_T)$ and{,}  
due to the parabolic Sobolev inequality (in the 
form of \cite[chapter II, Lemma 3.3]{lsu}), also
$V,\,V_x\in L^{r}(\rr_T)$ for any $q \le r \le p'$ where 
$p':=\big(\frac{1}{q}-\frac{1}{3}\big)^{-1}=\frac{6p}{3+p}$.  
Further, the constants in this parabolic Sobolev 
inequality can be chosen uniformly over $p$ in 
any given compact interval, hence for some uniform $C_8<\infty$ and 
\textit{all} $p \in [3/2,3)$,
\eq\label{bd-vx}
\|V_x\|_{L^{p'}(\rr_T)} \leq C_8 \|h\,A(R)R_x\|_{L^q(\rr_T)} 
\en 
where $q=\frac{2p}{p+1}$ and $p'=\frac{6p}{3+p}$.
Similarly, the function 
$$
Z(t,x):=\int_{\rr} F_{\mu}(y) \varphi(t,x-y)\,\mathrm{d}y,
$$
is a classical solution of the initial value problem
\eq\label{Z_cauchy}
Z_t-Z_{xx}=0\,,\quad Z(0,\cdot)=F_{\mu}\,.
\en
Since $F_\mu(y)=\int_{{-\infty}}^y \theta(z)\mathrm{d}z$ for
$\theta:=\frac{\mathrm{d}\mu}{\mathrm{d}x}$, clearly
$Z_x$ is given by the convolution 
(in space) of $\theta$ and the heat kernel 
$\varphi(t,\cdot)$.
Consequently, by Fubini's theorem and Young's inequality,
for any $p,q \ge 1$ such that $1/p+1/q=1/r+1$,
$$
\|Z_x\|_{L^r(\rr_T)}^r
= \int_0^T \|Z_x(t,\cdot)\|_{L^r(\rr)}^r \mathrm{d}t
\le \|\theta\|_{L^q(\rr)}^r 
\int_0^T \|\varphi(t,\cdot)\|_{L^p(\rr)}^r \mathrm{d} t \,.  
$$
By assumption $\mu \in \Mstar$, so 
$\theta\in L^1(\rr)\cap L^{q_0}(\rr)$ for some $q_0>1$,
and hence
the norms $\|\theta\|_{L^q(\rr)}$, $1\leq q\leq q_0$ are
uniformly bounded. Further,
$\|\varphi(t,\cdot)\|_{L^p(\rr)} \le t^{-(p-1)/(2p)}$ 
for all $t>0$ and $p \ge 1$. Considering $q=q_0 \wedge r${,} 
with $p \ge 1$ such that for $1 \le r \le 3$ 
the value of $r(p-1)/(2p)$ is bounded away from one,
we conclude that $\|Z_x\|_{L^r(\rr_T)}$, 
$1\leq r \leq 3$ are uniformly (in $r$) 
bounded by some finite $C_9=C_9(T,\|\theta\|_{q_0})$.

From \eqref{R_cauchy}, \eqref{V_cauchy}{,} and
\eqref{Z_cauchy} it follows by a direct computation 
that $U:=R-V-Z$ is a generalized solution of 
the Cauchy problem
\eq\label{U_cauchy}
U_t-(A(R)U_x)_x=\big[(A(R)-1)V_x+(A(R)-1)Z_x\big]_x\,,\quad U(0,\cdot)=0 \,,
\en
with respect to test functions in $\DxS$ (interpreted via
integration by parts in $t$ and $x$, similarly to what we
have done in \eqref{gen-sol}).

\noindent\textbf{Step 2.}
We proceed to obtain the uniform in $p$ bounds on
$L^p$-norms of $U_x$ (and thereby those of $R_x$)
out of the bounds we already have for $V_x$ and $Z_x$. 
To this end, recalling part (a) of Lemma \ref{lem-l2rx}, 
refining the norm estimates in 
\cite[Theorem 6.2]{KR2} to estimates in $L^p(\rr_T)\cap L^{p'}(\rr_T)$ with $p'=\frac{6p}{3+p}$ and applying the method of continuity (see e.g. \cite[section III.2]{Li}), 
we deduce the existence of a solution $\hat{U}$ of the problem \eqref{U_cauchy} in the space ${\HH}_p(T)\cap {\HH}_{p'}(T)$ defined in \cite{KR2}. In particular, $\hat{U}\in L^{p'}(\rr_T)$, $\hat{U}_x\in L^{p'}(\rr_T)$ and
\eq\label{bd-ux}
\|\hat{U}\|_{L^{p'}(\rr_T)}+\|\hat{U}_x\|_{L^{p'}(\rr_T)}\leq C_{10} \big(\|V_x\|_{L^{p'}(\rr_T)}+\|Z_x\|_{L^{p'}(\rr_T)}\big)
\en
with $p'=\frac{6p}{3+p}$ and where the finite constant $C_{10}=C_{10}(T,\|A\|_\infty,\underline{a})$ can be 
chosen uniformly for all $2 \leq p'<3$. 

We show in Step 3 below that
$U=\hat{U}$ Lebesgue a.e. on $\rr_T$. Thus,  
$U_x\in L^{p'}(\rr_T)$ and consequently 
also $R_x=U_x+V_x+Z_x\in L^{p'}(\rr_T)$. 
Hence, combining the norm bounds \eqref{bd-vx} and
\eqref{bd-ux} with H\"older's inequality 
\eqref{Holder_bnd}, yields the estimate
\eq\label{bd-rx}
\|R_x\|_{L^{p'}(\rr_T)} \leq C_{11}\,
\|R_x\|^{1/2}_{L^p(\rr_T)}+ C_{12}
\en
for the finite constants 
$C_{11}=(1+C_{10}) C_8 (C_6 \|A\|_\infty)^{1/2}$ and $C_{12}=(1+C_{10}) C_9$ 
(both independent of $\frac{3}{2}\leq p<3$). 
Since $p' \ge p > 1$, by Jensen's inequality 
(for the convex function $x^{(p'-1)/(p-1)}$ and 
the probability measure $T^{-1} R_x\,\mathrm{d}m$ 
on $\rr_T$), we deduce that for any $p \in [3/2,3)$ 
$$
\|R_x\|_{L^p(\rr_T)}^{r(p)} \le C_{13} 
\|R_x\|_{L^{p'}(\rr_T)} \,,
$$
where $r(p) =\frac{p(p'-1)}{p'(p-1)} \ge 1$ and 
$C_{13} =\max(1,T^{1/2})$
is finite.
Combining this with \eqref{bd-rx} we conclude that
$$
\|R_x\|_{L^p(\rr_T)} \leq \max \{ 1, 
C_{13} (C_{11} \,\|R_x\|^{1/2}_{L^p(\rr_T)}+C_{12}) \} \,.
$$
As explained before, having $\|R_x\|_{L^p(\rr_T)}$ bounded,
uniformly over $\frac{3}{2}\leq p<3${,} yields  
that $R_x\in L^3(\rr_T)$ (whose $L^3$-norm is bounded by 
some $C(T,I)$ finite).

\noindent\textbf{Step 3.} Turning to show that 
Lebesgue a.e. $U=\hat{U}$ on $\rr_T$, we start by 
verifying that $R - Z \in L^p(\rr_T)$. 
Indeed, as $|R-Z| \le 1$, it suffices to check this 
for $p=1$, in which case by the triangle inequality, 
\begin{align}\label{rz-decom}
\int_{\rr_T} |R-Z| \,\mathrm{d}m
&\le 
\int_{\rr^+_T} (1-R)\,\mathrm{d}m+\int_{\rr^+_T} (1-Z)\,\mathrm{d}m+\int_{\rr^-_T} R\,\mathrm{d}m
+\int_{\rr^-_T} Z\,\mathrm{d}m \nonumber \\
&= \int_{\rr_T} |x| R_x \,\mathrm{d}m+\int_{\rr_T} |x| Z_x \,\mathrm{d}m \,,   
\end{align} 
where $\rr_T^+=[0,T]\times\rr_+$, 
$\rr_T^-=[0,T]\times\rr_-$ and the last equality applies
since $R(t,\cdot)$ and $Z(t,\cdot)$ are \abbr{cdf}-s having densities $R_x$ and $Z_x$ with 
respect to Lebesgue measure on $\rr_T$.
Since $\gamma\in \cA_{\iota,\mu}$ the first term 
on the right-side of \eqref{rz-decom}      
is finite (see \eqref{eq:eta-mom}), 
whereas the second term amounts to 
$\int_0^T \mathbf{E}[|Y+W(2t)|] \mathrm{d}t$ 
for $Y$ of law $\mu$, independent of 
the standard Brownian motion $\{W\}$. By the 
triangle inequality, the latter term is at most 
$\frac{1}{3}\,(2T)^{3/2}+T\,\int_\rr |x|\,\mathrm{d}\mu${,}
hence finite (since $\mu \in \Mstar$). 

Now, with $R-Z \in L^p(\rr_T)$, and having seen already 
in Step 1 that $V \in L^p(\rr_T)$ (by the parabolic 
Sobolev inequality for $V \in W^{1,2}_q(\rr_T)$), 
we conclude that $U=R-Z-V \in L^p(\rr_T)$. Similarly,
$U_x\in L^p(\rr_T)${,} since $R_x\in L^p(\rr_T)$ (from 
Lemma \ref{lem-lprx}), and we have already established in 
Step 1 that $Z_x\in L^p(\rr_T)$ and $V_x\in L^p(\rr_T)$. 

\smallskip
\noindent
Recall that $h \in L^1(R_x \mathrm{d}m)${,} so $V \in L^1(\rr_T)${,}
and hence $U\in L^1(\rr_T)$ as well. Now, fixing 
$g\in\oS$ we let 
\eq\label{eq:f-def}
f(t,x)=-\int_x^\infty g(t,y)\,\mathrm{d}y,\quad(t,x)\in\rr_T 
\en
and claim that, for Lebesgue almost every $0\leq t_1<t_2\leq T$,
\eq\label{U_PDE_weak}
\begin{split}
\int_\rr f(t_2,x)\,U(t_2,x)\,\mathrm{d}x-\int_\rr f(t_1,x)\,U(t_1,x)\,\mathrm{d}x - \int_{t_1}^{t_2}\int_\rr U\,f_t\;\mathrm{d}m \quad\quad\quad\\
=-\int_{t_1}^{t_2}\int_\rr \big[A(R)\,U_x + (A(R)-1)(V_x+Z_x) \big]\,f_x\;\mathrm{d}m.      
\end{split}
\en
Indeed, from the weak formulation of the \abbr{PDE} in \eqref{U_cauchy} we have
\eqref{U_PDE_weak} when $f \in \DxS$. This in turn
extends to all $f$ as in \eqref{eq:f-def}, by the 
uniform joint approximation on compacts of the continuously differentiable 
and bounded $(f,f_t,f_x)$ by $(\widetilde{g},\widetilde{g}_t,\widetilde{g}_x)$ 
for some $\widetilde{g} \in \DxS$ (taking such compact 
large enough that the contribution to \eqref{U_PDE_weak} from its complement 
be as small as desired, which we can do since
$A(R)U_x+(A(R)-1)(V_x+Z_x) \in L^p(\rr_T)$ and  
$U \in L^1(\rr_T)$). Next, integrating by parts in space on the left side of \eqref{U_PDE_weak}, we conclude that the function $W(t,x)=\int_{-\infty}^x U(t,y)\,\mathrm{d}y$ is a generalized solution of the Cauchy problem
\eq
W_t=A(R)W_{xx}+(A(R)-1)(V_x+Z_x),\quad W(0,\cdot)=0
\en  
on $\rr_T$ (with respect to test functions from $\oS$). 
This, $W_{xx}=U_x\in L^p(\rr_T)$, $V_x\in L^p(\rr_T)${,} and $Z_x\in L^p(\rr_T)$ imply $W_t\in L^p(\rr_T)$. Thus, in view of \cite[norm estimate (6.1)]{KR2}, we have that 
$U\in{\HH}_p(\rr_T)$, 
so that $U=\hat{U}$ Lebesgue almost everywhere by \cite[Theorem 2.4]{KR2}. 
\ep

\smallskip
\noindent\textbf{Proof of Proposition \ref{L3prop} (B).} 
{From} the proof of Lemma \ref{lem-l3rx} we recall that
$R-Z,Z_x,R_x \in L^r(\rr_T)$ for all $r \in [1,3]$ and
$V \in W^{1,2}_q(\rr_T)$ for all $q \in [6/5,3/2]$
(where{,} since $R_x\in L^3(\rr_T)$, the argument leading 
to \eqref{bd-vx} now applies also for $p=3$ and $q=3/2$).
Consequently, $V_t,V_{xx} \in L^{q}(\rr_T)$ for any $q \in [6/5,3/2]$
and{,} by the parabolic Sobolev inequality (in the form of 
\cite[chapter II, Lemma 3.3]{lsu}), also $V,V_x \in L^{p'}(\rr_T)$
for any $p'=(\frac{1}{q}-\frac{1}{3})^{-1} \in [3/2,3]$. 
Since $U=R-Z-V${,} this in turn implies that $U,U_x \in L^{p}(\rr_T)$
for all $p \in [3/2,3]$. Further, $Z_t=Z_{xx} \in L^q(\rr_T)$ 
for all $q \in [6/5,3/2]$, hence it suffices to show that 
$U \in W_{q}^{1,2}(\rr_T)$, as then $U_t,U_{xx} \in L^q(\rr_T)$
implying the same for $R_t$ and $R_{xx}$.
To this end, rewriting  
\eqref{U_cauchy} we have that $U$ solves 
\eq\label{U_non_div}
U_t-A(R)U_{xx}=A'(R)R_x\,U_x+\big[(A(R)-1)V_x+(A(R)-1)Z_x\big]_x 
\en   
with respect to test functions in $\DxS$, starting 
at $U(0,\cdot)=0$. 
For any $q \in [6/5,3/2]$ there exist $r,p \in [3/2,3]$ such that 
$\frac{1}{p}+\frac{1}{r}=\frac{1}{q}${,} so with $A$ and $A'$ bounded,
by H\"older's inequality and the preceding integrability properties,
the right side in \eqref{U_non_div} belongs to $L^{q}(\rr_T)$ for 
all $q \in [6/5,3/2]$. Thus, fixing $q \in [6/5,3/2]$ we may apply 
\cite[Theorem 2.1]{KR2} to deduce that there is a function 
$\tilde{U}\in W^{1,2}_{q}(\rr_T) \cap W^{1,2}_{3/2}(\rr_T)$ which satisfies
\eqs
\tilde{U}_t-A(R)\tilde{U}_{xx}=A'(R)R_x\,U_x+\big[(A(R)-1)V_x+(A(R)-1)Z_x\big]_x 
\ens 
and $\tilde{U}(0,\cdot)=0$. In particular, 
$\tilde{U},\,\tilde{U}_x\in L^{3}(\rr_T)$. 
Now, we let $\phi_k$, $k\in\nn$ be a truncation sequence such that $\phi_k\in C^\infty(\rr)$, $0\leq\phi_k\leq 1$, $\phi_k\equiv1$ on $[-k,k]$, $\phi_k\equiv0$ on $(-\infty,-k-1]\cup[k+1,\infty)$ and $\max(|\phi_k'|,|\phi_k''|)\leq2$. Next, fixing $k\in\nn$ we set $\tilde{\Delta}=\phi_k(U-\tilde{U})$. Then, $\tilde{\Delta}$ is a 
generalized solution of the problem
\eq\label{tilde_Delta_PDE}
\tilde{\Delta}_t-(A(R)\tilde{\Delta}_x)_x+A'(R)R_x\tilde{\Delta}_x=\tilde{\psi}_k,\quad \Delta(0,\cdot)=0
\en
with respect to test functions in $\DxS$, where
\eqs
\tilde{\psi}_k=-A(R)\phi_k''(U-\tilde{U})-2A(R)\phi_k'(U-\tilde{U})_x
\ens
is in $L^3(\rr_T)$. Further, $\tilde{\psi}_k(t,x) = 0$ for all $x \in (-k,k)${,} 
so by our choice of $\phi_k$,   
\eq\label{psi_k-norm}
\lim_{k \to \infty} \|\tilde{\psi_k}\|_{L^3(\rr_T)} \, = 0 \,.
\en
Now, a careful reading of the proof of \cite[chapter III, Theorem 3.3]{lsu} shows that the solution of the problem \eqref{tilde_Delta_PDE} in the space $W^{0,1}_2(\rr_T)$ is unique and satisfies
\eq\label{W0,1bound}
\|\tilde{\Delta}\|_{W^{0,1}_2(\rr_T)}
\leq C_9 \Big(\int_0^T \Big(\int_\rr |\tilde{\psi_k}|^{q_1}\,\mathrm{d}x\Big)^{q_2/q_1}\,\mathrm{d}t\Big)^{1/q_2}
\en
for all $q_1\in[2,\infty]$, $q_2\in[2,4]$ 
with $\frac{1}{q_1}+\frac{2}{q_2}=1$, provided that
\eqs
\int_0^T \Big(\int_\rr |A'(R)R_x|^{q_1}\,\mathrm{d}x\Big)^{q_2/q_1}\,\mathrm{d}t<\infty.
\ens
We choose $q_1=q_2=3$, so that the latter condition is satisfied. In addition,
by \eqref{psi_k-norm} and \eqref{W0,1bound}, the norm $\|\tilde{\Delta}\|_{W^{0,1}_2(\rr_T)}$ tends to $0$ in the limit $k\rightarrow\infty${,} and we conclude 
that $U=\tilde{U} \in W^{1,2}_{q}(\rr_T)$, as claimed. 
\ep 

\section{Proof of Proposition \ref{L3prop} (C)}\label{sec_dir}

The proof of Prop. \ref{L3prop} (C) consists of four steps.
In {Step 1} we convert the variational formula 
$I_{\iota,\mu} (\gamma) < \infty$ into the formula \eqref{sup_psi_u} 
corresponding to a suitable one-dimensional reversible diffusion. 
{Step 2} then deduces the existence of a sufficiently regular 
solution to the corresponding backward Cauchy problem \eqref{PDE_psi_u}, 
which enables  
us to employ Dirichlet form calculus for establishing in {Step 3} 
the integrability of $R_{xx}^2/R_x$ \red{for $R=\Rg$}. 
From the latter we deduce in 
{Step 4} that $R_t^2/R_x$ is integrable.

\smallskip
\noindent
\textbf{Step 1.} The functions $a(t,x)=A(R)$ and $b(R)$ are uniformly bounded and 
our assumption that $I_{\iota,\mu} (\gamma) < \infty$ implies that 
$\mathrm{d}\gamma\,\mathrm{d}t$ has density
$R_x \in L^3(\rr_T)$ with respect to Lebesgue measure 
$\mathrm{d}m$ (by Prop. \ref{L3prop} (A)). 
These facts imply by multiple applications of H\"older's inequality that the 
functional $g \mapsto \Phi_\gamma (g)-(g,g)_{\gamma}$ in \eqref{eq:tilde-I-def}  
is continuous with respect to the norm
\eqs
\|g\|_\infty+\|g_t\|_{L^{3/2}(\rr_T)}+\|g_x\|_{L^{3/2}(\rr_T)}+\|g_x\|_{L^3(\rr_T)}+\|g_{xx}\|_{L^{3/2}(\rr_T)}.
\ens
Therefore, denoting by $\widetilde{W}^{1,2}_{3/2}(\rr_T)$ the sub-space 
of all $g\in C_b(\rr_T)$ for which
$g_t,\,g_{xx}\in L^{3/2}(\rr_T)$ and $g_x\in L^{3/2}(\rr_T)\cap L^3(\rr_T)$,
the assumption 
$I_{\iota,\mu} (\gamma)<\infty$ amounts to 
\eq\label{sup_0}
\sup_{g\in\widetilde{W}^{1,2}_{3/2}(\rr_T)} \Big[(\gamma,g)(T)-(\gamma,g)(0)
-\int_0^T (\gamma,{\RR}^\gamma g+a g_x^2)(t)\,\mathrm{d}t\Big]<\infty\,,
\en
for ${\RR}^\gamma$ of \eqref{eq:rgamma-def}. 
Now, fixing $\psi \in C^\infty(\rr)$ such that 
$\lim_{|x|\rightarrow\infty}\frac{\psi(x)}{|x|}=1$, $\|\psi'\|_\infty<\infty$ 
and $\alpha_0:=e^{-\psi(x)}\,\mathrm{d}x$ is a probability measure, we 
introduce the parabolic operator
\eq\label{eq:rpsi-def}
{\RR}^{\gamma,\psi}
=\frac{\partial}{\partial t}+e^{\psi(x)}\frac{\partial}{\partial x}a(t,x)e^{-\psi(x)}\frac{\partial}{\partial x}
={\RR}^\gamma+(a_x-a\psi'-b(R))\frac{\partial}{\partial x}
\en
and show next that \eqref{sup_0} implies that
\eq\label{sup_psi}
\sup_{g\in \widetilde{W}^{1,2}_{3/2}(\rr_T)} 
\Big[(\gamma,g)(T)-\log \big(\int_{\rr} e^{g(0,x)} \mathrm{d}\alpha_0\big)  
-\int_0^T (\gamma,{\RR}^{\gamma,\psi} g+a g_x^2)(t) \,\mathrm{d}t\Big]<\infty.
\en 
Indeed, by Cauchy-Schwarz we have the bound
$$
\kappa_3:=\int_{\rr_T} (b(R)-a_x+a\psi')g_x\,R_x\,\mathrm{d}m \leq C_2 \sqrt{(g,g)_\gamma}
$$
for the finite, positive
$$
C_2:=\underline{a}^{-1/2} 
\big(\|b\|_\infty T^{1/2}+\|A'\|_\infty\|R_x\|_{L^3(\rr_T)}^{3/2}\big)
+\|\sqrt{a} \psi'\|_\infty T^{1/2} \,.
$$
Further, the implication
\begin{align}
&\forall\;\lambda>0:\;\lambda\kappa_1\leq C_1+\lambda^2\kappa_2,\;\;\kappa_3\leq C_2\sqrt{\kappa_2} \nonumber \\
\Longrightarrow \quad &\forall\;\tilde{\lambda}>0:\;\tilde{\lambda}\kappa_1+\tilde{\lambda}\kappa_3\leq 2C_1+\frac{1}{2}C_2^2+\tilde{\lambda}^2\kappa_2
\label{eq:implication}
\end{align} 
holds for all $\kappa_1,\kappa_3 \in\rr$ and $\kappa_2$, $C_1$, $C_2$ positive. 
Therefore, scaling the test functions $g$ in \eqref{sup_0} and \eqref{sup_psi}
by $\lambda>0$ and $\tilde{\lambda}=2\lambda$, respectively,  
then considering \eqref{eq:implication} for 
$C_1=I_{\iota,\mu} (\gamma)$, $\kappa_1=\Phi_\gamma (g)$ and
$\kappa_2=(g,g)_\gamma$ proves that
the change from ${\RR}^\gamma$ 
to ${\RR}^{\gamma,\psi}$ in \eqref{sup_0} does not make the supremum infinite. 
Moreover, the change in value
due to  
the terms corresponding to the initial condition, is bounded by
$$
\sup_{g\in \widetilde{W}^{1,2}_{3/2}(\rr_T)}\Big[\int_\rr g(0,x)\mathrm{d}\mu
-\log\int_\rr e^{g(0,x)}\,\mathrm{d}\alpha_0 \Big]
\leq H(\mu|\alpha_0)
$$
(see for example \cite[Lemma 6.2.13]{dz}). The latter relative entropy is finite since
\eq
H(\mu|\alpha_0)=\int_\rr \log\frac{\mathrm{d}\mu}{\mathrm{d}\alpha_0}\,\mathrm{d}\mu
=\int_\rr \log\frac{\theta}{e^{-\psi}}\,\mathrm{d}\mu
\leq\int_\rr |\psi|\,\mathrm{d}\mu+\int_\rr \log\theta\,\mathrm{d}\mu \,,
\en
while $\int_\rr |x|\,\mathrm{d}\mu<\infty$ and $\theta \in L^{q_0}(\rr)$ 
for some $q_0>1$ (by definition, for $\mu \in \Mstar$).

Next, let $\EE \widetilde{W}^{1,2}_{3/2}(\rr_T)$ denote the 
collection of $u : \rr_T \mapsto (0,\infty)$ such that 
$\log u = g \in \widetilde{W}^{1,2}_{3/2}(\rr_T)$. It is easy to 
verify that $\EE \widetilde{W}^{1,2}_{3/2}(\rr_T)$ consists
of all positive $u \in \widetilde{W}^{1,2}_{3/2}(\rr_T)$ which are
bounded away from zero, in terms of which \eqref{sup_psi} becomes  
\eq\label{sup_psi_u}
\sup_{u\in \EE \widetilde{W}^{1,2}_{3/2}(\rr_T)} \Big[(\gamma,\log u)(T) - 
\log \int_\rr u(0,x)\,\mathrm{d}\alpha_0
-\int_0^T \big(\gamma,\frac{{\RR}^{\gamma,\psi} u}{u}\big)(t)\,\mathrm{d}t\Big]
<\infty \,.
\en

\smallskip
\noindent\textbf{Step 2.} We claim that to every 
$f\in\oS$ corresponds 
$u\in \EE\widetilde{W}^{1,2}_{3/2}(\rr_T)$ such that  
\eq\label{PDE_psi_u}
{\RR}^{\gamma,\psi} u - f\,u = 0,\quad u(T,\cdot)=1
\en
where all terms of \eqref{PDE_psi_u} are then in $L^{3/2}(\rr_T)$ 
(by definition 
$u_x,R_x \in L^r(\rr_T)$ for all $r \in [3/2,3]$, hence 
$(a_x-a\psi')u_x\in L^{3/2}(\rr_T)$ by the Cauchy-Schwarz inequality and boundedness of 
$A',A$ and $\psi'$). Clearly, having such solution for \eqref{PDE_psi_u} 
amounts to finding a solution $w$ of
\eqs
{\RR}^{\gamma,\psi} w - f\,w = f,\quad w(T,\cdot)=0
\ens
or{,} equivalently (see \eqref{eq:rpsi-def}), a solution of 
\eq\label{wPDE}
w_t+(aw_x)_x-\psi'\,a\,w_x- f\,w = f, \quad w(T,\cdot)=0
\en
such that $u=(w+1) \in \EE\widetilde{W}^{1,2}_{3/2}(\rr_T)$.
To this end, we
employ \cite[Theorem 6.2]{KR2} together with the method of continuity to first
get a generalized solution $w$ of \eqref{wPDE} in the space $W^{0,1}_6(\rr_T)\cap W^{0,1}_2(\rr_T)$. Indeed, the norm bound in \cite[inequality (6.3)]{KR2} extends to a norm bound for functions in ${\HH}^1_6\cap {\HH}^1_2$ with respect to the norms $\|\cdot\|_{\mathbb{H}^{-1}_6}+\|\cdot\|_{\mathbb{H}^{-1}_2}$ and $\|\cdot\|_{\mathbb{H}^1_6}+\|\cdot\|_{\mathbb{H}^1_2}$ defined in \cite{KR2} (one only needs to add the norm bounds \cite[inequality (6.3)]{KR2} for $p=6$ and $p=2$). Applying the method of continuity (see \cite[section III.2]{Li}), and relying on such a refined norm estimate to interpolate between the \abbr{PDE} \eqref{wPDE} and the corresponding \abbr{PDE} with a smooth coefficient $a$, we find a solution of \eqref{wPDE} in ${\HH}^1_6\cap {\HH}^1_2$, which in particular belongs to $W^{0,1}_6(\rr_T)\cap W^{0,1}_2(\rr_T)$. Moreover,
by H\"older's inequality $a_x w_x \in L^2(\rr_T) \cap L^{3/2}(\rr_T)$
(since $a_x = A'(R) R_x \in L^3(\rr_T)$ and $w_x \in L^r(\rr_T)$ for all $r \in [2,6]$). 
Hence, 
refining the norm bounds in \cite[Theorem 4.1]{KR2} to an estimate on the norms in $L^2(\rr_T)\cap L^{3/2}(\rr_T)$ and applying the method of continuity in a similar fashion, we also have a generalized solution $\hat{w}\in W^{1,2}_2(\rr_T)\cap W^{1,2}_{3/2}(\rr_T)$ of the equation
\eq
\hat{w}_t+a\,\hat{w}_{xx}-a\,\psi'\hat{w}_x - f\,\hat{w}=-a_x\,w_x+f,\quad \hat{w}(T,\cdot)=0\,.
\en
Proceeding to show that $w=\hat{w}$, we let $\{\phi_k, k\in\nn\}$ be the same
truncation sequence as in the proof of Prop. \ref{L3prop} (B), 
fix $k\in\nn${,} and set $\Delta:=\phi_k(\hat{w}-w)$. Then, $\Delta\in W^{0,1}_2(\rr_T)$ is a generalized solution of
$$
\Delta_t+(a\Delta_x)_x-(a_x+\psi'\,a) \Delta_x- f\,\Delta = \psi_k, \quad \Delta(T,\cdot)=0,
$$ 
where
$$
\psi_k=\phi_k''\,a(\hat{w}-w)+2\,\phi_k'\,a (\hat{w}-w)_x-\phi_k'\,\psi'\,a(\hat{w}-w).
$$
As in the derivation of \eqref{psi_k-norm}, since 
$\max(|\phi_k'|,|\phi_k''|) \leq 2{\cdot}\mathbf{1}_{|x| \in [k,k+1]}$, 
having $a$, $\psi'$ uniformly bounded and $(\hat{w}-w),(\hat{w}-w)_x\in L^2(\rr_T)$
implies that 
\eq\label{eq:psik-l2-norm}
\lim_{k\rightarrow\infty}\|\psi_k\|_{L^2(\rr_T)}=0\,.
\en
Hence, writing $\Delta=\phi_k\,\hat{w}-\phi_k\,w$, following the paragraph after the statement of \cite[chapter III, Theorem 3.3]{lsu}{,} and applying the energy inequality 
of \cite[chapter III, Theorem 2.1]{lsu}, {with} $r=q=2$ and $n=1$
(so that $\frac{2}{r}+\frac{n}{q} \le \frac{n+4}{2}$), we conclude that 
\eq\label{eq:Delta-l2-norm}
\|\Delta\|_{L^2(\rr_T)} \le C_3 
\Big(\int_0^T 
\Big(\int_\rr |\psi_k|^{q}\,\mathrm{d}x\Big)^{r/q}\,\mathrm{d}t\Big)^{1/r}
\rightarrow 0\,,  \quad \textrm{as} \quad  k\rightarrow\infty \,.
\en 
Therefore, {$m$}-a.e. $w=\hat{w}$ on $\rr_T$. All in all, we have found a solution $u$ 
to \eqref{PDE_psi_u} such that $w=u-1$ is an element of 
$W^{1,2}_2(\rr_T)\cap W^{1,2}_{3/2}(\rr_T)$, and hence also 
$u_x=w_x\in L^6(\rr_T)$ (by the parabolic Sobolev inequality 
in the form of \cite[chapter II, Lemma 3.3]{lsu} for $p=6$ and $q=2$).

It thus remains only to show that $u\in C_b(\rr_T)$ and that $u$ is 
bounded away from $0$ on $\rr_T$. To establish this 
we first apply \cite[chapter III, Theorem 5.2]{lsu} to find a generalized solution $\tilde{w}$ of \eqref{wPDE} in the subspace of $W^{0,1}_2(\rr_T)$ whose elements satisfy
\eq
\ess\,\sup_{t\in[0,T]}
 \int_\rr \tilde{w}(t,x)^2\,\mathrm{d}x+\int_{\rr_T} \tilde{w}_x^2\,\mathrm{d}m < \infty. 
\en
Next, we apply \cite[Theorem 10 (vi)]{Ar}, with the constant $\gamma>0$ there being arbitrarily small, to conclude that $\tilde{u}:=\tilde{w}+1$ has to be the unique 
generalized solution of \eqref{PDE_psi_u} in the sense of \cite{Ar}. It is thus 
given by
\eq
\tilde{u}(t,x)=\int_\rr \Gamma(t,x;T,y)\,\mathrm{d}y
\en
with $\Gamma$ denoting the weak 
fundamental solution of \eqref{PDE_psi_u}, defined as in \cite[Section 6]{Ar}. Now, \cite[Theorem C]{Ar} implies that $\tilde{u}$ is locally H\"older continuous in $(t,x)$, and hence continuous in $(t,x)$ on $[0,T)\times\rr$. Putting this together with \cite[Theorem 10 (vi)]{Ar}, we conclude that $\tilde{u}$ is continuous on the whole of $\rr_T$. Finally, we use the heat kernel estimates on $\Gamma$ 
from \cite[Theorem 7]{Ar} to conclude that $\tilde{u}$ has to be bounded between two positive constants. Therefore, all we need to show now is that $\tilde{u}=u$ or{,} equivalently{,} $\tilde{w}=w$. To this end, re-using the truncation sequence $\{\phi_k\}$, we fix 
$k\in\nn$ and set $\tilde{\Delta}:=\phi_k(\tilde{w}-w)$. Then 
$\tilde{\Delta}\in W^{0,1}_2(\rr_T)$ is a generalized solution of
\eqs
\tilde{\Delta}_t+(a\tilde{\Delta}_x)_x-\psi'\,a\,\tilde{\Delta}_x- f\,\tilde{\Delta} = \tilde{\psi}_k, \quad \tilde{\Delta}(T,\cdot)=0
\ens 
where
\eqs
\tilde{\psi}_k=\phi_k''\,a(\tilde{w}-w)+2\,\phi_k'\,a(\tilde{w}-w)_x+\phi_k'\,a_x\,(\tilde{w}-w)-\phi_k'\,\psi'\,a(\tilde{w}-w).
\ens
The $L^2(\rr_T)$-norm of $\phi_k' a_x(\tilde{w}-w)$ decays 
as $k \to \infty$ (because $a_x \in L^3(\rr_T)$ and 
$\|\tilde{w}-w\|_{L^6(S_{k,4} \cap S_{-k,4})}$ is uniformly bounded in $k$).
Thus, with $w,\tilde{w},w_x,\tilde{w}_x \in L^2(\rr_T)$ we deduce
as in the derivation of \eqref{eq:psik-l2-norm} that 
\eqs
\lim_{k\rightarrow\infty}\|\tilde{\psi}_k\|_{L^2(\rr_T)}=0 \,.
\ens
This implies, as in the derivation of \eqref{eq:Delta-l2-norm}, that
%
$\|\tilde{\Delta}\|_{L^2(\rr_T)}\rightarrow0$ as $k\rightarrow\infty$, 
and consequently that {$m$}-a.e. $w=\tilde{w}$ on $\rr_T$. 

\smallskip
\noindent\textbf{Step 3.} In view of \eqref{PDE_psi_u}, we get from \eqref{sup_psi_u} that 
\eqs
\infty > \sup_{f\in\oS} \Big[-\log \int_\rr u^{(f)}(0,x)\,\mathrm{d}\alpha_0-\int_0^T (\gamma,f)(t) \,\mathrm{d}t\Big]\,,
\ens
where $u^{(f)} \in \EE\widetilde{W}^{1,2}_{3/2}(\rr_T)$ satisfies \eqref{PDE_psi_u}.
We then deduce that{,} for some $C \in (0,\infty)$,  
\eq\label{supg}
\infty>\sup_{g\in\oS} \Big[ 
- \int_0^T (\gamma,g_x+C\,g^2)(t)\,\mathrm{d}t \Big]
=\sup_{g\in\oS} \int_{\rr_T} (g\,R_{xx}-C\,g^2\,R_x)\,\mathrm{d}m
\en
by showing that the $L^2(\alpha_0)$-norm of $u^{(f)}(0,\cdot)$ 
(and hence also $\log \int_\rr u^{(f)}(0,x)\,\mathrm{d}\alpha_0$),
is uniformly bounded over all $f=g_x+C\,g^2$ with $g\in\oS$.
To this end, recall {from} \eqref{eq:rpsi-def} that 
${\RR}^{\gamma,\psi} v = v_t + {\LL} v$ for 
$v \in \widetilde{W}^{1,2}_{3/2}(\rr_T)$ and the bounded linear
operator 
$$
{\LL} v = e^{\psi} (a e^{-\psi} v_x)_x : 
\widetilde{W}^{1,2}_{3/2}(\rr_T) \mapsto L^{3/2}(\rr_T) \,.
$$
We further set ${\LL}_t \phi = e^{\psi}(a(t,x) e^{-\psi} \phi')_x$ for
$\phi \in C^\infty_c(\rr)$. Then, considering positive $\phi_k \in \oS$ that 
converge to $u^{(f)}$ in $\widetilde{W}^{1,2}_{3/2}(\rr_T)$, we deduce from \eqref{PDE_psi_u}
that{,} 
for any $s \in [0,T]$,
\begin{align}\label{eq:uf-bd}
1
- \|u^{(f)}(s,\cdot)\|^2_{L^2(\alpha_0)}
&=2\int_{s}^{T}\mathrm{d}t \, \int_\rr u^{(f)}\,u^{(f)}_t\,\mathrm{d}\alpha_0
\nonumber \\
&=-2\int_{s}^{T}\mathrm{d} t \,
\int_\rr u^{(f)} ({\LL}\,u^{(f)}-f\,u^{(f)})\,\mathrm{d}\alpha_0 \nonumber \\ 
&\geq -2\int_{s}^{T}\lambda_f (t)\,\|u^{(f)}(t,\cdot)\|^2_{L^2(\alpha_0)}\,\mathrm{d}t
\end{align} 
where
$$
\lambda_f (t)
= \sup_{\phi\in C_c^\infty(\rr):\,\phi \geq0,\,\|\phi\|_{L^2(\alpha_0)}=1}\;\; 
\int_\rr \phi({\LL}_t \,\phi-f(t,\cdot) \,\phi)\,\mathrm{d}\alpha_0 \,.
$$
Recall that $\mathrm{d}\alpha_0 = e^{-\psi} \mathrm{d} x$, so 
integrating by parts and then taking $\phi \mapsto \sqrt{\phi}$ yields 
\begin{align*}
\lambda_f (t) 
&=\sup_{\phi\in C_c^\infty(\rr):\,\phi\geq0,\,\|\phi\|_{L^2(\alpha_0)}=1} \;\;
\Big[-\int_\rr a(t,\cdot) \,(\phi')^2\,\mathrm{d}\alpha_0-\int_\rr f(t,\cdot)\,\phi^2\,\mathrm{d}\alpha_0 \Big]\\
&=\sup_{\phi\in C_c^\infty(\rr):\,\phi\geq0,\,\|\phi\|_{L^1(\alpha_0)}=1}\;\; 
\Big[-\frac{1}{4}\int_\rr \frac{(\phi')^2}{\phi}\,a(t,\cdot)\,\mathrm{d}\alpha_0
-\int_\rr f(t,\cdot)\,\phi\,\mathrm{d}\alpha_0\Big] \,.
\end{align*}
Further, for any $g \in \oS$ and smooth $\phi$ such that 
$\phi \alpha_0 \in M_1(\rr)$, we get by integration by parts and
the Cauchy-Schwarz inequality that
\eq
\Big|\int_\rr g_x(t,\cdot)\,\phi\,\mathrm{d}\alpha_0\Big|=\Big|\int_\rr g(t,\cdot) (\phi'-
\phi \psi')\,\mathrm{d}\alpha_0\Big|\leq\sqrt{\upsilon_1\upsilon_2}+\|\psi'\|_\infty\sqrt{\upsilon_2}, 
\en 
{with} $\upsilon_1=\int_\rr \frac{(\phi')^2}{\phi}\,\mathrm{d}\alpha_0$ and
$\upsilon_2=\upsilon_2(t)=\int_\rr g^2(t,\cdot)\,\phi\,\mathrm{d}\alpha_0$. Hence, for 
$f=g_x + C g^2$ with $g \in \oS${,} we end up with
\eq\label{lambda_bnd}
\sup_{t \in [0,T]} \{\lambda_f (t)\}\leq\sup_{\upsilon_1,\upsilon_2\geq0} 
\Big[-\frac{1}{4} \underline{a}\upsilon_1+\sqrt{\upsilon_1\upsilon_2}+\|\psi'\|_\infty\sqrt{\upsilon_2}-C\upsilon_2\Big]\,.
\en
For $C>0$ large enough, the right side of \eqref{lambda_bnd} is 
finite, hence applying Gr\"onwall's lemma for 
$s \mapsto\|u^{(f)} (T-s,\cdot)\|^2_{L^2(\alpha_0)}$ we deduce from 
\eqref{eq:uf-bd}
that $\|u^{(f)}(0,\cdot)\|_{L^2(\alpha_0)}$ is bounded uniformly 
(over such $f$), as claimed. 

Next, let $\widehat{\mathbb H}$ denote the 
closure of $\oS$ with respect to the 
$L^2(R_x \mathrm{d}m)$-norm. From \eqref{supg} we know that{,} for 
some finite $C_1$, $C$ and all $g \in \oS$,
$$
\Big|\int_{\rr_T} g \, R_{xx} \, \mathrm{d} m \Big| \le C_1 + C \|g\|^2_{\widehat{\mathbb H}} \,.
$$
Hence, $g\mapsto\int_{\rr_T} g\,R_{xx}\,\mathrm{d}m$  
extends to a bounded linear functional on $\widehat{\mathbb H}$, which 
by the Riesz representation theorem is of the form $g\mapsto(g,\hat{h})_{\widehat{\mathbb H}}$ for some $\hat{h}\in\widehat{\mathbb H}$. 
The identity $\int_{\rr_T} g\,R_{xx}\,\mathrm{d}m=\int_{\rr_T} g\,\hat{h}\,R_x\,\mathrm{d}m$ for all $g\in\oS$ implies that $m$-a.e. 
if $R_x=0$ then $R_{xx}=0$, with $\hat{h}=\frac{R_{xx}}{R_x}$ (subject to our 
running convention that $\frac{0}{0}=0$). Consequently, 
\eq\label{eq:Rxx-bd}
\int_{\rr_T} \frac{R_{xx}^2}{R_x} \mathrm{d} m = 
\int_{\rr_T} \hat{h}^2 R_x \mathrm{d}m = \|\hat{h}\|^2_{\widehat{\mathbb H}} < \infty \,. 
\en

\smallskip
\noindent\textbf{Step 4.} Following the derivation of \eqref{eq:Rxx-bd} out
of \eqref{supg}, if 
\eq\label{supg-t}
\sup_{g\in\oS} \int_{\rr_T} (2 g\,R_t - g^2\,R_x)\,\mathrm{d}m<\infty
\en
then $\int_{\rr_T} \frac{R_t^2}{R_x}\,\mathrm{d}m<\infty$. Further, plugging in
\eqref{supg-t} the value of $R_t$ from the \abbr{PDE} \eqref{R_cauchy}, 
all the terms of which are in $L^{3/2}(\rr_T)$,
we find that \eqref{supg-t} amounts to 
\eq\label{eq:supg-t-mod}
\sup_{g\in\oS} \int_{\rr_T} 
\big(2 g [A'(R)R_x+\hat{h} A(R)+h \,A(R)]- \,g^2\big) \,R_x\, \mathrm{d}m < \infty \,,
\en
for $\hat{h} \in {\widehat{\mathbb H}}$ of \eqref{eq:Rxx-bd}
and $h \in L^2(R_x \mathrm{d}m)$ of \eqref{htildeh}. 
Point-wise optimizing in \eqref{eq:supg-t-mod}
over the value of $g(t,x)$ at each point of $\rr_T$,
bounds the 
supremum by the finite
$\|A'(R)R_x + \hat{h} A(R) + h A(R)\|^2_{L^2(R_x \mathrm{d} m)}$,
thereby completing the proof. 
\ep

\begin{rmk}
A crucial step in the proof of Prop. \ref{L3prop} (C) consists of 
showing that the continuous $u=w+1$ solving \eqref{PDE_psi_u}, with  
$w\in W^{1,2}_2(\rr_T) \cap W^{1,2}_{3/2}(\rr_T)$, is further bounded 
and bounded away from $0$. In doing so we relied on the results of \cite{Ar},
but we note in passing that with some additional work such positivity can 
be obtained from \cite[Corollary 4.6]{KR3}, bypassing the need for \cite{Ar}. 
\end{rmk}
\section{Proof of Proposition \ref{approx_prop_main}}\label{sec_approx}

Fixing throughout $\iota \in (0,1]$ and $\mu \in \Miota$, 
we start with the \textit{convexity} of the functionals from 
which $J_{\iota,\mu}(\cdot)$ is composed, followed by its use in 
establishing convergence results for ${\widetilde J} (\cdot)$. 
\begin{lemma}\label{conv_lemma}
The functionals 
\eq\label{eq:def-J123}
J^{(1)} (R) =\int_{\rr_T} \frac{R_t^2}{R_x}\,\mathrm{d}m,\quad
J^{(2)}(R) =\int_{\rr_T} \frac{R_{xx}^2}{R_x}\,\mathrm{d}m,\quad
J^{(3)}(R) =\int_{\rr_T} {R}_x^3\,\mathrm{d}m
\en
are convex on the set $\FF=\FF_{3/2}$ of \eqref{eq:thm-13}. 
\end{lemma}

\noindent\textbf{Proof.} 
The convexity of $J^{(3)}(\cdot)$ 
is an immediate consequence of the convexity of $x\mapsto x^3$ on $[0,\infty)$.
Further, recall from Steps 3 and 4 in the proof of Prop. \ref{L3prop} (C), 
that on $\FF$
\begin{align*}
J^{(1)}({R})
&=\sup_{g\in \oS} \int_{\rr_T} \big(2\,g\,{R}_t-g^2 {R}_x\big)\,\mathrm{d}m,\\
J^{(2)}({R})
&=\sup_{g\in\oS} \int_{\rr_T} \big(2\,g\,{R}_{xx}-g^2 {R}_x\big)\,\mathrm{d}m\,,
\end{align*}
so each of these functionals, being a supremum of linear functionals, must 
therefore be convex.  \ep

\begin{lemma}\label{lem-aprox_rate}
Suppose $A(\cdot) \ge \underline{a} >0$ with $b(\cdot)$ and 
$A'(\cdot)$ continuous and bounded.
Let $R=\Rg$ for $\gamma\in \CC$ such that $J_{\iota,\mu}(\gamma)<\infty$ and suppose 
the strictly positive probability densities $R_x^\epsilon \in C^{1,1}(\rr_T)$ 
are such that  
\begin{equation}\label{T123}
\limsup_{\epsilon\downarrow0} 
J^{(\ell)}(R^{\epsilon}) \leq J^{(\ell)}(R), \quad \ell=1,2,3 \,,
\end{equation}
$R^\epsilon \to R$ uniformly on compacts, 
$R^\epsilon_x \to R_x$ in $L^p(\rr_T)$, $p \in [2,3]$
and $m$-a.e. $R^\epsilon_t \to R_t$, $R^\epsilon_{xx} \to R_{xx}$.
If $\gamma^\epsilon = R^\epsilon_x \mathrm{d}x$ are such that 
$\gamma^\epsilon(0) \in \Mstar$, then
\eq\label{rate_bnd}
\lim_{\epsilon\downarrow 0} {\widetilde J} (\gamma^\epsilon) = 
J_{\iota,\mu} (\gamma)\,.
\en    
\end{lemma}

\smallskip
\noindent\textbf{Proof.}
\noindent\textbf{Step 1.} We first show that as $\epsilon\downarrow0$:
\begin{align}
\int_{\rr_T} \Big|\frac{R^\epsilon_t}{(R^\epsilon_x)^{1/2}}-\frac{R_t}{(R_x)^{1/2}}\Big|^2\,\mathrm{d}m &\rightarrow 0, \label{T1claim}\\
\int_{\rr_T} \Big|\frac{R^\epsilon_{xx}}{(R^\epsilon_x)^{1/2}}-\frac{R_{xx}}{(R_x)^{1/2}}\Big|^2\,\mathrm{d}m &\rightarrow 0, \label{T2claim}\\
\int_{\rr_T} \big|(R^\epsilon_x)^{3/2}-(R_x)^{3/2}\big|^2\,\mathrm{d}m &\rightarrow 0. \label{T3claim}
\end{align}
To this end, with $|x_1-x_2|^2 \le |x_1^2-x_2^2|$ whenever $x_1,x_2 \in \rr_+$
and $R^\epsilon_x\rightarrow R_x$ in $L^2(\rr_T)$,
it follows that{,} for $\epsilon\downarrow 0$,
\eq
\int_{\rr_T} \big|(R^\epsilon_x)^{1/2}-(R_x)^{1/2}\big|^2\,\mathrm{d}m\rightarrow 0 \,. \label{T4claim} 
\en
Similarly, combining the inequality 
$|x_1^{3}-x_2^{3}|\leq \frac{3}{2}|x_1^2-x_2^2| \max(x_1,x_2)$ 
(which holds for all $x_1,x_2 \in \rr_+$) with H\"older's inequality we find that
\begin{align*}
\int_{\rr_T} \big|(R^\epsilon_x)^{3/2}-(R_x)^{3/2}\big|^2\,&\mathrm{d}m 
\leq \frac{9}{4} \int_{\rr_T} |R^\epsilon_x-R_x|^2\,\max(R_x^\epsilon,R_x)\mathrm{d}m \quad\quad\quad\quad\quad\quad\quad\quad \\
&\leq \frac{9}{4} \Big(\int_{\rr_T} |R^\epsilon_x-R_x|^3\,\mathrm{d}m\Big)^{2/3} 
\Big(\int_{\rr_T}\max(R_x^\epsilon,R_x)^3\,\mathrm{d}m\Big)^{1/3}.
\end{align*}
By assumption $R^\epsilon_x\rightarrow R_x$ in $L^3(\rr_T)$ and{,} 
due to \eqref{T123} for $\ell=3$, 
the norms $\|R^\epsilon_x\|_{L^3(\rr_T)}$ are uniformly bounded, thereby
yielding \eqref{T3claim}.
\newline
Turning to prove \eqref{T1claim}, we let $c^\epsilon=R^\epsilon_t (R^\epsilon_x)^{-1/2}$
which by \eqref{T123} {with} $\ell=1$ is $L^2(\rr_T)$-bounded. Hence, for any 
sequence $\epsilon_k \downarrow 0$ we have that  
$c^\epsilon \rightarrow c^*$ weakly in $L^2(\rr_T)$ along some sub-sequence 
(by the Banach-Alaoglu theorem, where $c^* \in L^2(\rr_T)$ may depend on the
chosen sub-sequence). By the triangle inequality, we thus get 
from \eqref{T4claim} that along this sub-sequence, for any fixed $\psi\in C_c(\rr_T)$, 
\eq\label{T1conv1}
\int_{\rr_T} (R^\epsilon_x)^{1/2} c^\epsilon\,\psi\,\mathrm{d}m \rightarrow \int_{\rr_T} (R_x)^{1/2} c^*\,\psi\,\mathrm{d}m \,.
\en
Recall that $R_t \in L^{3/2}(\rr_T)$ (since $R \in \FF$), and by assumption 
$m$-a.e. $(R^\epsilon_x)^{1/2} c^\epsilon=R^\epsilon_t \rightarrow R_t$.
Thus, for any fixed $\psi\in C_c(\rr_T)${,} the l.h.s. of \eqref{T1conv1} 
converges to $\int_{\rr_T} R_t \psi \mathrm{d} m$ 
as $\epsilon \downarrow 0$, resulting with  
\eq\label{T1conv2}
\int_{\rr_T} 
(R_x)^{1/2} c^*\,\psi\,\mathrm{d}m = \int_{\rr_T} R_t\,\psi\,\mathrm{d}m \,.
\en 
Further, with $(R_x)^{1/2}\in L^6(\rr_T)$ and $c^*\in L^2(\rr_T)$, 
by H\"older's inequality $(R_x)^{1/2} c^*\in L^{3/2}(\rr_T)${,}
so from \eqref{T1conv2} we conclude that $m$-a.e. 
$c^*=R_t (R_x)^{-1/2}$, independently of the sequence $\epsilon_k$. 
That is, $R^\epsilon_t (R^\epsilon_x)^{-1/2} \rightarrow R_t (R_x)^{-1/2}$ 
weakly in $L^2(\rr_T)$ when $\epsilon \downarrow 0$. Together with 
the $L^2(\rr_T)$-norm-bound
of \eqref{T123} for $\ell=1$, this 
yields the (strong) convergence of \eqref{T1claim}.  
\newline
Finally, the same argument, just with $R^\epsilon_{t}$ replaced by
$R^\epsilon_{xx}\,$, yields \eqref{T2claim}.

\smallskip
\noindent\textbf{Step 2.} To deduce \eqref{rate_bnd} from 
\eqref{T1claim}--\eqref{T4claim}, recall that $J_{\iota,\mu}(\gamma)<\infty$ 
requires $\gamma \in \cA_{\iota,\mu}$ and{,} in particular, that 
$R=\Rg \in C_b(\rr_T)$ satisfies \eqref{eq:Runif-tight} for some
$M=M_\alpha$ finite (and all $\alpha>0$). Thus, our assumption 
that $R^\epsilon \to R$ uniformly on compact sets, combined
with the monotonicity of the distribution 
functions $x \mapsto R^\epsilon(t,x)$, $x \mapsto R(t,x)$
and \eqref{eq:Runif-tight}, yields that $R^\epsilon \to R$ 
uniformly on $\rr_T$, when $\epsilon\downarrow 0$. 
This, and the assumed continuity of $A'$ 
and $b$ on $[0,1]$, show that as $\epsilon \downarrow 0$, 
\eq\label{eq:unif-sigAb}
\sigma(R^\epsilon)\rightarrow \sigma(R),\;\; A'(R^\epsilon)\rightarrow A'(R),\;\; b(R^\epsilon)\rightarrow b(R) \quad \textrm{uniformly on } \rr_T\,.
\en
Moreover, all functions appearing in \eqref{eq:unif-sigAb} are
uniformly bounded on $\rr_T$. 
Putting this together with the uniform positivity of $\sigma$ and 
\eqref{T1claim}--\eqref{T4claim}, we have the following convergences in $L^2(\rr_T)$ 
when $\epsilon\downarrow0$,
\begin{align}
\frac{1}{\sigma(R^\epsilon)}\,\frac{R^\epsilon_t}{(R^\epsilon_x)^{1/2}}&\rightarrow \frac{1}{\sigma(R)}\,\frac{R_t}{(R_x)^{1/2}},\\
\sigma(R^\epsilon)\frac{R^\epsilon_{xx}}{(R^\epsilon_x)^{1/2}}&\rightarrow \sigma(R)\frac{R_{xx}}{(R_x)^{1/2}},\\
\frac{A'(R^\epsilon)}{\sigma(R^\epsilon)}\,(R^\epsilon_x)^{3/2}&\rightarrow \frac{A'(R)}{\sigma(R)}\,(R_x)^{3/2},\\
\frac{b(R^\epsilon)}{\sigma(R^\epsilon)}(R^\epsilon_x)^{1/2}&\rightarrow \frac{b(R)}{\sigma(R)}(R_x)^{1/2}.
\end{align}
By H\"older's inequality the finiteness of 
$J^{(\ell)}(R^\epsilon)$, $\ell=1,2,3$
implies that 
$R^\epsilon = \Rge \in \FF$, 
and by our assumptions, also {$\gamma^\epsilon \in \tcA$.}
Thus, as $\epsilon\downarrow 0$ we have that 
\begin{align*}
{\widetilde J} (\gamma^\epsilon)&=\frac{1}{2}\Big\|\frac{R^\epsilon_t-(A(R^\epsilon)R^\epsilon_x)_x}{\sigma(R^\epsilon)(R^\epsilon_x)^{1/2}}
+\frac{b(R^\epsilon)}{\sigma(R^\epsilon)}(R^\epsilon_x)^{1/2}\Big\|_{L^2(\rr_T)}^2\\
&\longrightarrow \frac{1}{2}\Big\|
\frac{R_t-(A(R)R_x)_x}{\sigma(R)(R_x)^{1/2}}
+\frac{b(R)}{\sigma(R)}(R_x)^{1/2}
\Big\|_{L^2(\rr_T)}^2=J_{\iota,\mu} (\gamma).
\end{align*}
This completes the proof of Lemma \ref{lem-aprox_rate}. \ep

\red{
Under Asmp. \ref{main_ass}(a),(b) and (d), our 
key result provides the conclusions of Prop. \ref{approx_prop_main} 
for $J_{\iota,\mu}(\gamma)<\infty$ 
} 
and some $\gamma^\epsilon$ in the collection $\GGc_\iota$
that we define next (c.f. Definition \ref{new-def-eta} of $\GG$).
\begin{definition}\label{def:GGc}
Let $\GGc_\iota$ be the subset 
of $\{ \gamma \in \CC : {J}_{\iota,\mu}(\gamma) < \infty, \;$ some
$\; \mu \in \Miota \}$ such that 
$R=\Rg \in C^{\infty}_b(\rre)$ 
with $R_x$ strictly positive and 
\eq\label{eq:unif-der-control}
\big|\partial_t^{j} \partial_x^k R (t,x) \big| 
\le c_{j,k} \Rs (x) \qquad \forall (t,x) \in \rre 
\en
for some $\{c_{j,k},$ finite $j,k \in \nn \}$ and 
\eq\label{eq:rs-def}
\Rs (x) := \sup_{t \in [0,T]} \{ \widetilde{R}(t,x)\}\,, \qquad 
\widetilde{R}(t,x) := 1-R(t,|x|)+R(t,-|x|) \,.
\en
\end{definition}
\begin{prop}\label{approx_prop}
Suppose $A \ge \underline{a} >0$ with $b$ and $A'$ continuous and bounded.
If $J_{\iota,\mu}(\gamma)<\infty$ for some $\mu \in \Miota$, then 
there exist $\{ \gamma^\eps \} \subset \GGc_\iota$  
such that 
$\gamma^\eps \rightarrow\gamma$ in $\CC$ 
as $\eps \downarrow 0$,
$\sup_\eps \int |x|^{1+\iota} \mathrm{d} \gamma^\eps (0)$ is finite,
and \eqref{rate_bnd} holds. 
\end{prop}

\smallskip
\noindent\textbf{Proof.} The proof consists of three steps.
In Step 1 we construct  
$\gamma^{\del,\eps} \in \CC$ whose 
smooth \abbr{cdf} paths
$R^{\del,\eps}=\red{R^{(\gamma^{\del,\eps})}} 
\in C_b^\infty(\rre)$ satisfy
\eqref{eq:unif-der-control} and
having strictly positive \abbr{pdf}-s $R^{\del,\eps}_x$. 
Step 2 confirms that $\gamma^{\del,\eps} \to \gamma$ in $\CC$ 
when $(\del,\eps) \to (0,0)$ and that $\gamma^{\del,\eps} \in \cA_{\iota,\mu}$ 
for $\mu=\gamma^{\del,\eps}(0)$ whose $(1+\iota)$-th moments are
bounded, uniformly over $(\del,\eps)$.
Then, relying on Lemma \ref{lem-aprox_rate}, we verify
in {Step 3} that \eqref{rate_bnd} holds for $\del=\del(\eps)$ 
small enough (so{,} in particular, such $\gamma^{\del,\eps} \in \GGc_\iota$). 

\noindent\textbf{Step 1}. Let $\phi \in C^\infty(\rr)$ 
be a
strictly positive probability density with
$\int |x|^k \phi(x) \mathrm{d} x < \infty$ and 
$|\phi^{(k)}(x)| \le c_k \phi(x)$ for some finite 
$c_k$, $k \ge 1$ and all $x \in \rr$ 
(\red{for example, the smoothing near $x=0$ 
of 
$e^{-2|x|}$ provides such $\phi$}).
With
$\phi_\eps (y)=\eps^{-1} \phi (y/\eps)$, for each
$\del,\eps \in (0,1)$ we consider the function
\eq\label{eq:S-de-ep-def}
S^{\del,\eps}(t,x)=
\int_\rr R((t-3\del)_+,y)\,\phi_\eps (x-y) \,\mathrm{d}y
\en
on $\rr_{T+3\del}$. Then, fixing a probability density 
$\psi \in C_c^\infty(\rr)$ supported on $[0,3]$ with 
\eq\label{eq:cpsi-def}
c_\psi=\inf_{s \in [1,2]} \{\psi(s)\} > 0 \,,
\en
we set $\psi_\del(s)=\del^{-1} \psi(s/\del)$ and consider 
\eq\label{eq:R-de-ep-def}
R^{\delta,\epsilon} (t,x)=
\int_0^{T+3\delta} S^{\delta,\epsilon} (s,x)\,\psi_\del(s-t)\,\mathrm{d}s, 
\quad \delta,\epsilon \in (0,1)\,, \quad (t,x) \in \rre\,.
\en
With $\psi(t)\phi(x)$ a probability density on $\rre$, since
each $R(t,\cdot)$ is a \abbr{cdf}, so are $S^{\delta,\epsilon}(t,\cdot)$
and $R^{\del,\eps}(t,\cdot)$. By the strict positivity of $\phi$ 
we have the same for $S^{\delta,\epsilon}_x${,} and thereby 
also for $R^{\del,\eps}_x$. Next, with $\phi(\cdot)$ smooth, 
$S^{\delta,\epsilon}(t,\cdot) \in C_b^\infty(\rr)$, for each 
$t$, $\epsilon$ and $\delta$, hence
$R^{\del,\eps} \in C_b^{\infty}(\rr_T)$ by smoothness of $\psi$.
Moreover, as $\phi(\cdot)$ point-wise controls its derivatives
it follows that, for each $k \ge 1$ and all $(t,x) \in \rre$,
\red{
\[
\big|\partial_x^k S^{\delta,\epsilon} (t,x) \big| \le c_k \epsilon^{-k}
S^{\del,\eps}(t,x) \,.
\] 
Further, with $\int \phi_{\epsilon}(x-y)\mathrm{d}y =1$, 
the same bound holds for $1-S^{\del,\eps}$, resulting with}
\eq\label{eq:unif-der-cont}
\big|\partial_x^k S^{\delta,\epsilon} (t,x) \big| 
\red{\le c_k \epsilon^{-k} 
[(1-S^{\delta,\epsilon}(t,x)) \wedge S^{\delta,\epsilon}(t,x)]} 
\le c_k \epsilon^{-k} \widetilde{S}^{\del,\eps}(t,x) 
\en  
(where $\widetilde{S}^{\del,\eps}$ is related to $S^{\del,\eps}$ analogously to  
\eqref{eq:rs-def}). 
Thus, having $\psi(\cdot)$ smooth yields, by \eqref{eq:unif-der-cont} and 
\eqref{eq:R-de-ep-def}, that 
\eq\label{eq:der-jk-bd}
|\partial^j_t \partial^k_x R^{\delta,\epsilon}(t,x)|
\le \|\, \psi_\delta^{(j)}\, \|_{\infty} \, c_k \epsilon^{-k} V^{\del,\epsilon}(t,x)
\en
for all $j,k \in \nn$ where 
$$
V^{\del,\epsilon} (t,x) = 
\int_0^{3\delta} \widetilde{S}^{\del,\eps}(t+u,x) 
\, \mathrm{d}u \,.
$$
Clearly, for any $\delta>0$ and $(t,x) \in \rre$,
$$
V^{\del,\eps}(t,x) 
\le 
\frac{3 \delta}{c_\psi} 
\sup_{t \in [0,T]} \; \Big\{
\int_0^{3\delta} \widetilde{S}^{\del,\eps} (t+u,x) 
\psi_\delta (u) \mathrm{d} u \, \Big\} 
= \frac{3 \delta}{c_\psi} \Rs^{\del,\epsilon} (x) \,,
$$
for $\Rs(\cdot)$ of \eqref{eq:rs-def}, which together with 
\eqref{eq:der-jk-bd} implies that $R^{\del,\epsilon}$ 
satisfies \eqref{eq:unif-der-control}. 

\noindent\textbf{Step 2.} Since $S^{\del,\eps}(t,x)=S^{\del,\eps}(0,x)$ when
$t \in [0,3\delta]$, upon specializing \eqref{eq:R-de-ep-def} to $t=0$, 
we get the smooth \abbr{CDF}-s 
\begin{equation}\label{eq:theta-delta}
\red{R^{\del,\eps}(0,x) =
}
\ev [R(0,x-Y_\eps)] =: \Theta_\eps (x) \,,
\end{equation}
for $Y_\eps$ of density $\phi_\eps$. Our assumption 
that $\gamma(0) = \mu \in \Miota$  
translates to $\gamma^{\delta,\epsilon}(0) \in \Miota$ 
with uniformly bounded $(1+\iota)$-moments
(since $Y_\eps$ has this property and 
$f \mapsto 
(\int f(x)^{q_0} \mathrm{d}x,\int |x|^{1+\iota} f(x) \mathrm{d}x)$ 
is convex over probability densities). Similarly, from 
\eqref{eq:S-de-ep-def} and \eqref{eq:R-de-ep-def} 
$$
\int_{\rr_T}   |x|^{1+\iota}\,\mathrm{d}\gamma^{\delta,\epsilon}(t)\,\mathrm{d}t 
\le \int_{\rr_T} \ev [|x-Y_\eps|^{1+\iota}]\,\mathrm{d}\gamma(t) \mathrm{d}t  
+ 3\del \int_{\rr} \ev [|x - Y_\eps|^{1+\iota}]\,\mathrm{d}\gamma(0) \,,
$$
with the r.h.s. finite{,} since $\gamma$ satisfies \eqref{eq:eta-mom},
$\gamma(0) \in \Miota${,} and $Y_\eps$ has finite moments of all order.
That is, $\gamma^{\delta,\epsilon}$ satisfy the moment condition 
\eqref{eq:eta-mom} and are thus in $\cA_{\iota,\gamma^{\del,\eps}(0)}$.
In addition, by dominated convergence, 
$R^{\delta,\epsilon} \to R$ uniformly on compacts, so in particular,
$\gamma^{\delta,\epsilon}\rightarrow\gamma$ in $\CC$.

\noindent\textbf{Step 3}. By construction, $m$-a.e. $R^{\del,\eps}_x \to R_x$, 
$R^{\delta,\epsilon}_t \to R_t$
and $R^{\delta,\epsilon}_{xx} \to R_{xx}$, whenever 
$(\delta,\epsilon) \to (0,0)$. Thus, from  
Lemma \ref{lem-aprox_rate} we get \eqref{rate_bnd}
for $\delta=\delta(\epsilon)$ small, once we show that  
\begin{equation}\label{T123m}
\limsup_{\eps \downarrow 0} 
\limsup_{\del \downarrow 0} J^{(\ell)}(R^{\delta,\epsilon}) \leq J^{(\ell)}(R), \quad \ell=1,2,3 
\end{equation}
(indeed, by Egorov's theorem, \eqref{T123m} for $\ell=3$ implies 
that $R^{\delta,\epsilon}_x \to R_x$ in $L^p(\rr_T)$, $p \in [2,3]$).
Turning to prove \eqref{T123m}, recall from
its definition 
that $J_{\iota,\mu}(\gamma)<\infty$ implies $R=\Rg \in \FF$ and 
$$
J^{(\ell)}(R)=\int_0^{3\del} \psi_\del(s) \mathrm{d} s 
\int_\rr J^{(\ell)} (R(\cdot,y+\cdot)) \phi_\eps(y) \mathrm{d} y<\infty\, \quad \ell=1,2,3\,. 
$$
Hence, starting with the functional $J^{(1)}(\cdot)$, 
we get from the definition of $R^{\delta,\epsilon}$, 
upon applying Lemma \ref{conv_lemma} twice, that
\begin{align*}
J^{(1)} (R^{\delta,\epsilon}) \leq \int_0^{3\del} 
J^{(1)} (S^{\delta,\epsilon}(s+\cdot,\cdot)) \psi_\del (s) 
\mathrm{d} s
\leq J^{(1)} (R) 
\,,
\end{align*}
hence \eqref{T123m} holds for $\ell=1$. 
Next, consider the functionals
\eq\label{eq:wh-J-def}
\widehat{J}^{(2)} (F)=\int_\rr \frac{F''(x)^2}{F'(x)} \mathrm{d}x
\,,\qquad 
\widehat{J}^{(3)} (F)=\int_\rr F'(x)^3 \mathrm{d}x\,,
\en
which{,} as in the proof of Lemma \ref{conv_lemma},  
are convex on the set of twice differentiable \abbr{cdf}-s. 
Thus, for $\ell=2,3${,} we get by the same argument as before that 
\begin{align*}
J^{(\ell)} (R^{\delta,\epsilon}) \leq J^{(\ell)} (R) + 
3 \delta \widehat{J}^{(\ell)} (\Theta_\eps) \,.
\end{align*}  
Our choice of $\phi_\eps$ results with
$|\Theta''_\eps| \le c_1 \eps^{-1} \Theta'_\eps$
and uniformly bounded \abbr{pdf} $\Theta'_\eps$. Consequently,
$\widehat{J}^{(\ell)} (\Theta_\eps) < \infty$
and taking $\del \downarrow 0$ establishes \eqref{T123m} for $\ell=2,3$.
\ep

\smallskip
\red{Specializing to $\iota=1$ and}
applying Prop. \ref{approx_prop} to approximate the given
$\gamma \in \CC$ of finite $J_{1,\mu}(\gamma)$ by a 
suitable sequence from $\GGc_1$, the proof 
of Prop. \ref{approx_prop_main} is thus completed 
\red{by considering our next lemma}
(to get $\gamma^{\ell} \in \GG$ that suitably converge 
to any given $\gamma \in \GGc_1$).
\begin{lemma}\label{lem:new-ofer2}
Suppose $\gamma \in \GGc_\iota$ and Asmp. \ref{main_ass}(a),(b) and (d) hold.
Then,
there exist $\mu \in \Miotas$ 
such that $\gamma^\ell =(1-\ell^{-1}) \gamma + \ell^{-1} \mu \in \GG$ with 
${\widetilde J} (\gamma^\ell) \to {\widetilde J} (\gamma)$,
and $\gamma^\ell \to \gamma$ in $\CC$.
\end{lemma}
\noindent
{\bf Proof}. 
Recall $c_{0,1}=c_{0,1}(\gamma)$ from \eqref{eq:unif-der-control}.
We begin by showing that{,} with $\eta=1/(12 c_{0,1})$,
the function $\Rs$ of \eqref{eq:rs-def} is such that 
\eq\label{eq:Rs-y-Rs-x}
|x-y| \le 3 \eta \quad \Longrightarrow \quad \Rs (y) \le 2 \Rs (x) \,.
\en
Indeed, 
since $\Rs(z)=\Rs(-z)$ is continuous and non-increasing in $|z|$,
it suffices to prove that $\Rs(y) \le 2 \Rs(x)$ for $|x| \in [|y|,|y|+3\eta]$ 
and $x y \ge 0$. To this end, note that by 
\eqref{eq:unif-der-control}, for all
$t \in [0,T]$,
$$
\big|R(t,y) - R(t,x)\big| = \Big|\int_{x}^{y} R_x(t,z) \mathrm{d}z\Big| 
\le  c_{0,1} \Big|\int_{x}^{y} \Rs(z) \mathrm{d}z\Big|
 \le \frac{1}{4} \Rs(y) \,.
$$
The preceding implies that{,} for $\widetilde{R}(t,x)$ of \eqref{eq:rs-def}, 
$$
\widetilde{R}(t,y) \le \widetilde{R}(t,x) + \frac{1}{2} \Rs(y) \,,
$$
so taking the maximum over $t$ results with 
$\Rs(y) \le \Rs(x) + \frac{1}{2} \Rs(y)$, as claimed.

We next claim that
\eq\label{eq:eta-mom-sup}
\int_\rr |x|^{\iota} \Rs (x) \mathrm{d} x < \infty \,.
\en 
Indeed, taking $\widehat{\eta}=1/(8 c_{1,0})>0$, we have 
that{,} for $|t-s| \le 2 \widehat{\eta}$, 
$$
|R(t,x)-R(s,x)| \le c_{1,0} |t-s| \Rs(x) \le \frac{1}{4} \Rs(x) \,,
$$
out of which we deduce that 
$\widetilde{R}(s,x) \le \widetilde{R}(t,x) + \frac{1}{2} \Rs(x)$.
Considering $s \in [0,T]$ such that $\widetilde{R}(s,x)=\Rs(x)$, we thus find
that $\widetilde{R}(t,x) \ge \frac{1}{2} \Rs(x)$ throughout a sub-interval 
of length at least $2 \widehat{\eta}$. Consequently, for all $x \ge 0$,
$$
\widehat{\eta} \Rs(x) \le \int_0^T \widetilde{R}(t,x) \mathrm{d} t \,, 
$$
hence by Fubini's theorem,
$$
\frac{\iota+1}{2} \widehat{\eta} 
\int_\rr |x|^{\iota} \Rs(x) \mathrm{d} x 
\le 
(\iota+1) \int_0^\infty x^{\iota} \mathrm{d} x 
\int_0^T \widetilde{R}(t,x) \mathrm{d} t 
=\int_{\rre} |x|^{\iota+1} \mathrm{d} \gamma(t) \mathrm{d} t 
$$
which for $\gamma \in \GGc_\iota$ is finite in view of \eqref{eq:eta-mom}. 

Continuing with the proof of the lemma,
let $\kappa = \int_\rr \Rs(y) \mathrm{d}y$, which by 
\eqref{eq:eta-mom-sup}
is finite.
Then, for $\eta>0$ of \eqref{eq:Rs-y-Rs-x}
and $\psi \in C_c^\infty(\rr)$ supported on $[0,3]$ 
as in Prop. \ref{approx_prop}, we construct the \abbr{pdf}
\eq\label{eq:r-prime-def}
r'(x) := \frac{1}{\kappa} \int_\rr \Rs(y) \psi_{\eta}(x-y) \mathrm{d}y \,,
\en
and{,} for each $\eps>0${,} consider the path $\gamma^\eps \in \CC$ associated with 
the \abbr{cdf}-s
$$
R^\eps(t,x) := (1-\eps) R(t,x) + \eps r(x) \,.
$$ 
Since $R_x$ is strictly positive, so are $R^\eps_x$. Further
$r \in C_b^\infty(\rr)${,} 
and consequently also $R^\eps \in C_b^\infty(\rre)$.  For $\eps \downarrow 0$
we clearly have that $R^\eps \to R$ uniformly on compacts (so 
$\gamma^\eps \to \gamma$ in $\CC$), and $m$-a.e. 
$R_t^\eps \to R_t$, $R_{xx}^\eps \to R_{xx}$. Since  
$r'(\cdot)$ and $R_x(t,\cdot)$ are both uniformly bounded \abbr{pdf}-s,
obviously also $R_x^\eps \to R_x$ in $L^p(\rre)$ for all $p \in [2,3]$. 
Thus, in view of Lemmas \ref{conv_lemma} and \ref{lem-aprox_rate}, we have 
that $\widetilde{J} (\gamma^\eps) \to \widetilde{J} (\gamma)$ 
provided we show that $\widehat{J}^{(\ell)}(r)$, $\ell=2,3$ of 
\eqref{eq:wh-J-def} are finite. The finiteness of $\widehat{J}^{(3)}(r)$ is 
trivial, for $r'$ is a bounded \abbr{pdf}, whereas 
$\widehat{J}^{(2)}(r) < \infty$ due to the boundedness 
of $|r''(x)|/r'(x)$. To see the latter, note that for 
$c_\psi>0$ as in \eqref{eq:cpsi-def}, by 
\eqref{eq:Rs-y-Rs-x}
and \eqref{eq:r-prime-def}, 
\eq\label{eq:r-prime-Rs}
\frac{c_\psi}{2\eta} \Rs(x) \le
\frac{c_\psi}{\eta} \inf_{y \in [x-3\eta,x]} \{ \Rs(y) \} \le
\kappa r'(x) \le \sup_{y \in[x-3\eta,x]} \{ \Rs(y) \} 
\le 2 \Rs(x)  \,.
\en
Similarly, $\kappa |r''(x)| \le 2 \|\psi'\|_1 \eta^{-1} \Rs(x)$ which together
with the l.h.s. of \eqref{eq:r-prime-Rs} implies that 
$|r''(x)|/r'(x)$ is uniformly bounded (by $4 \|\psi'\|_1/c_\psi$). 

Next, combining the r.h.s. of \eqref{eq:r-prime-Rs} with 
\eqref{eq:eta-mom-sup}{,} we deduce that the $\iota$-th moment of 
$\mu = r'(x) \mathrm{d} x$
is finite.
It thus remains only to show that the continuous function 
$$
h(R^\eps) = \frac{R^\eps_t-(A(R^\eps)R^\eps_x)_x}{A(R^\eps)\,R^\eps_x} 
$$
is further uniformly bounded and globally Lipschitz in $x$.
With $A$ bounded below, $A'$, $R^\eps_x$ bounded 
above, and $R_t^\eps=(1-\eps) R_t$, the boundedness of 
$h(R^\eps)$ follows from that of $(|R_t|+|R_{xx}|+|r''|)/r'$. To this end,
we have just shown the uniform boundedness of
$|r''(x)|/r'(x)$ and recall {from} \eqref{eq:unif-der-control} 
that $|R_t(t,x)|+|R_{xx}(t,x)| \le (c_{1,0}+c_{0,2}) \Rs(x)$, 
which by the l.h.s. of \eqref{eq:r-prime-Rs} is 
further bounded by $C r'(x)$, for some $C=C(\gamma)$ finite.
As for showing that $h(R^\eps)$ is globally Lipschitz continuous in $x$, 
note that
\begin{align*}
[h(R^\eps)]_x &= \frac{(1-\eps) R_{tx} - (A(R^\eps) R^\eps_x)_{xx}
- h(R^\eps) (A(R^\eps) R^\eps_x)_x}{A(R^\eps) R^\eps_x} 
\,.
\end{align*}
Recall that $A(\cdot)$ is bounded below, $h(R^\eps)$, $R_x$, $r'$
and $(
|R_{xx}|+|r''|)/r'$ are bounded above, 
and $A'$ is globally Lipschitz. We thus have that 
$|[h(R^\eps)]_x|
$ is uniformly bounded, provided that 
$(|R_{tx}|+|R_{xxx}|
+|r'''|)/r'$ is uniformly bounded.
The latter holds by the l.h.s. of \eqref{eq:r-prime-Rs}{,} since
from \eqref{eq:unif-der-cont} we have that 
$|R_{tx}|+|R_{xxx}|
\le C \Rs$ for some ${C=}C(\gamma)$ finite,
and by the same reasoning as above, we also get that
$\kappa |r'''| \le 2 \|\psi''\|_1 \eta^{-2} \Rs$.
\ep

\section{Proof of Proposition \ref{local_ubd_exp_tight}}


We suppose throughout this section 
that parts (a)--(c) of Asmp. \ref{main_ass} hold 
\red{for some $\ieta \in (0,1]$ and use the simplified notations 
$I(\cdot)$ and $\cA$ for $I_{\ieta,\rho_0}(\cdot)$ of \eqref{eq:tilde-I-def} 
and $\cA_{\ieta,\rho_0}$ of \eqref{eq:eta-mom}, respectively.}
In this setting we establish a local large deviations 
upper bound for $\rho^N$ being near $\gamma$, starting with 
$\gamma\in \cA$
(where
$I(\gamma)=\sup_{g\in\oS} \,[\Phi_\gamma (g) - (g,g)_\gamma]$).
\begin{prop}\label{contcenter} 
For each $\gamma\in \cA$ and any $g \in \oS$ 
we have 
the 
upper bound
\eq
\lim_{\delta\downarrow0}\limsup_{N\rightarrow\infty}\frac{1}{N}\log \pp(\rho^N\in B(\gamma,\delta)) \leq (g,g)_\gamma - \Phi_\gamma (g) \;.
\en
\end{prop}
\noindent\textbf{Proof}. 
Fixing $g\in\oS$, for each $t \in[0,T]$ 
and $\xi\in \CC${,} we set
\eqs
H^g(\xi)(t)=(\xi,{\RR}^\xi g)(t)\,{,}
\ens
with ${\RR}^\xi g = g_t+b(\Rxi)g_x+A(\Rxi)g_{xx}$ of \eqref{eq:rgamma-def}. Then,
applying It\^o's formula for the 
real-valued stochastic processes $Z_N^g(t):=(\rho^N,g)(t)$, one finds that 
\eq\label{N_mckean_vlasov}
Z_N^g(t)-Z_N^g(0) = \int_0^t H^g (\rho^N)(s) \mathrm{d} s + M_N^g(t) \,, 
\en
{with} the continuous martingale 
\eq\label{eq:MNg}
M_N^g(t):=\frac{1}{N} \sum_{i=1}^N \int_0^t 
\sigma(F_{\rho^N(s)}(X_i(s))) g_x(s,X_i(s))\,\mathrm{d}W_i(s)\,.
\en
Its quadratic variation is 
$\langle M_N^g \rangle (t) = \frac{1}{N}\int_0^t V^g(\rho^N)(s)\,\mathrm{d}s$,
{with} 
\eq\label{eq:Vg}
V^g(\xi)(t):=2 \big(\xi, A(\Rxi) g_x^2 \big)(t) \,.
\en
Hence, by  
the martingale representation theorem (see \cite[Theorem 3.4.2]{ks}),  
\eqs
M^g_N(t)=\frac{1}{\sqrt{N}} \int_0^t \sqrt{V^g(\rho^N)(s)}\,\mathrm{d}\beta_N(s) 
\ens
for some one-dimensional standard Brownian motion $\beta_N$. 
Next, fixing $\gamma\in \cA$, let
\eq
\overline{H}^g(\gamma)(t)=\int_0^t H^g(\gamma)(s) \mathrm{d}s \,,\quad
\overline{V}^g(\gamma)(t)=\int_0^t V^g(\gamma)(s) \mathrm{d}s \,,
\en
and recall that{,} by \eqref{eq:phi-def} and \eqref{eq:ip-def},
\eq\label{eq:Ig-Phig}
(g,g)_\gamma - \Phi_\gamma (g) = \frac{1}{2} \overline{V}^g(\gamma)(T)  
+ (\gamma,g) (0) + \overline{H}^g(\gamma)(T)  - (\gamma,g)(T) \,.
\en
We thus proceed to 
compare $Z_N^g(\cdot)$ with the process
\eq\label{eq:YNg-def}
Y_N^g(t)=(\gamma,g)(0) + \overline{H}^g(\gamma)(t)+M^{g,\gamma}_N(t),\quad t\in[0,T],
\en
having the martingale part 
\eqs
M^{g,\gamma}_N(t) = \frac{1}{\sqrt{N}} \int_0^t \sqrt{V^g(\gamma)(s)}\,
\mathrm{d}\beta_N(s),\quad t\in[0,T] \,.
\ens 

\smallskip
\noindent \textbf{Step 1.} We first show that{,}
on the event $\rho^N\in B(\gamma,\delta)$, 
$$
\|Z_N^g-Y_N^g\|_\infty := \sup_{t \in [0,T]} |Z_N^g(t) - Y_N^g(t)| \le \epsilon \,,
$$ 
up to a probability which is negligible at our large deviations 
exponential scale (in the limit $N \to \infty$ followed by $\delta \downarrow 0$). 
To this end, fixing $\rho\in B(\gamma,\delta)$ we note that  
\begin{align*}
\int_0^T \big|H^g(\rho)(s)-H^g(\gamma)(s)\big|\,\mathrm{d}s
&\leq \int_0^T \big|(\rho,{\RR}^{\rho} g)(s) - (\rho,{\RR}^{\gamma} g)(s)
\big|\,\mathrm{d}s\\
&\quad +\int_0^T \big|(\rho,{\RR}^{\gamma} g)(s)
-(\gamma,{\RR}^{\gamma} g) (s) \big|\,\mathrm{d}s
\end{align*}
and setting $\wt{R} = \Rro$, $R=\Rg$, 
we further have
\begin{align*}
\int_0^T \big|V^g(\rho)(s)-V^g(\gamma)(s)\big|\,\mathrm{d}s 
&\leq 2 \int_0^T \big|(\rho,A(\wt{R})g_x^2)(s)-(\rho,A(R) g_x^2)(s)\big|\,\mathrm{d}s\\
&\quad +2\int_0^T \big|(\rho,A(R) g_x^2)(s)-(\gamma,A(R) g_x^2)(s)\big|\,\mathrm{d}s.
\end{align*}
Since $\gamma \in \cA${,} we know that $R \in C_b(\rr_T)$
and consequently so are ${\RR}^{\gamma} g$ and $A(R) g_x^2$, from which we
deduce that the second term in both upper bounds tends to zero 
as $\delta \downarrow 0$, uniformly in $\rho\in B(\gamma,\delta)$. 
Now, recall that $b$ and $A$ are Lipschitz functions (under 
Asmp. \ref{main_ass}(a)-(b)). Hence, 
to get the same uniform convergence for the first term
in both upper bounds, it suffices to show that  
\eq\label{R_comp_est}
\lim_{\delta\downarrow0}\sup_{\rho\in B(\gamma,\delta)}\int_0^T
(\rho,|\wt{R}-R|)(s)\,\mathrm{d}s=0\,.
\en
Moreover, fixing $\delta>0$, by definition of the metric $d(\cdot,\cdot)$ on
$\CC$,  
for any $\rho\in B(\gamma,\delta)$ and $(s,x)\in\rr_T$ one has that
\eq\label{eq:dl-r}
R(s,x-\delta)-\delta\leq \wt{R}(s,x)\leq R(s,x+\delta)+\delta 
\en
(see \eqref{eq:dL-sup-def}), and hence, also
\eqs
|\wt{R}(s,x)-R(s,x)|\leq R(s,x+\delta)-R(s,x-\delta)+\delta\,.
\ens
Recall that $R(s,x+\delta)-R(s,x-\delta) = \gamma(s,(x-\delta,x+\delta])$,
resulting by Fubini's theorem and yet another application of the preceding bound with
\begin{align*}
(\rho,|\wt{R}-R|)(s) &\le
\delta + \int_{\rr} (R(s,x+\delta)-R(s,x-\delta)) \rho(s,\mathrm{d}x) \\
&= \delta + \int_{\rr} \rho(s,[y-\delta,y+\delta)) \gamma(s,\mathrm{d}y) \\
&\le 3\delta + \int_{\rr} (R(s,y+3\delta)-R(s,y-3\delta)) \gamma(s,\mathrm{d} y) \,.
\end{align*}
Integrating both sides over $s \in [0,T]$, we see that 
\eqref{R_comp_est} is a consequence of 
\eqs
\lim_{\delta\downarrow0} \int_{\rr_T} (R(s,y+3\delta)-R(s,y-3\delta))\,
\gamma(s,\mathrm{d}y)\,\mathrm{d}s=0 \,,
\ens
which in turn follows from the dominated convergence theorem and
the fact that $R \in C_b(\rr_T)$.
All in all, we have shown that 
\begin{align}
&\lim_{\delta\downarrow0}\sup_{\rho\in B(\gamma,\delta)} \int_0^T \big|H^g(\rho)(s)-
H^g(\gamma)(s)\big|\,\mathrm{d}s=0,\label{iden1}\\
&\lim_{\delta\downarrow0}\sup_{\rho\in B(\gamma,\delta)} \int_0^T \big|V^g(\rho)(s)-V^g(\gamma)(s)\big|\,\mathrm{d}s=0.\label{iden2}
\end{align}
Recall that $(\rho^N,g)(0) \to (\rho_0,g(0,\cdot)) =(\gamma,g)(0)$ by 
Asmp. \ref{main_ass}(c) and the definition of $\cA$. Hence,
comparing \eqref{N_mckean_vlasov} 
with \eqref{eq:YNg-def}{,}
it immediately follows from \eqref{iden1} that{,} for any fixed $\epsilon>0$,
\begin{align*}
 \lim_{\delta\downarrow0} & \limsup_{N\rightarrow\infty}
\frac{1}{N}\log\pp\big(\rho^N\in B(\gamma,\delta),\|Z^g_N-Y^g_N\|_\infty >2\epsilon\big)\\
&\leq \lim_{\delta\downarrow0}\limsup_{N\rightarrow\infty}
\frac{1}{N}\log\pp\big(\rho^N\in B(\gamma,\delta),\|M^g_N-M^{g,\gamma}_N\|_\infty>\epsilon\big). 
\end{align*}

\noindent
Recall \cite[Theorem 8.5.7]{ok} that $M^g_N(\cdot)-M^{g,\gamma}_N(\cdot)$
has the law of time-changed standard Brownian motion $\beta(\tau(\cdot))$, for 
$\tau(t)=\langle M^g_N-M^{g,\gamma}_N\rangle(t)$
(the quadratic variation process of the martingale $M^g_N-M^{g,\gamma}_N$). 
Moreover, on the event $\{\rho^N\in B(\gamma,\delta)\}$, we have that
$$
\tau(T)=\frac{1}{N}\int_0^T \Big(\sqrt{V^g(\rho^N)(s)}-\sqrt{V^g(\gamma)(s)}\Big)^2\,\mathrm{d}s\leq\frac{1}{N} \kappa_\gamma (\delta) \,,
$$
with $\kappa_\gamma (\delta) \to 0$ as $\delta\downarrow 0$ (due to 
the inequality $(\sqrt{x_1}-\sqrt{x_2})^2\leq|x_1-x_2|$ for
$x_1,x_2 \ge 0$ and \eqref{iden2}). Thus, from
Bernstein's inequality for Brownian motion, 
\begin{align*}
\lim_{\delta\downarrow0} & \limsup_{N\rightarrow\infty}
\frac{1}{N}\log\pp\big(\rho^N\in B(\gamma,\delta),\|M^g_N-M^{g,\gamma}_N\|_\infty 
>\epsilon \big)\\
&\leq \lim_{\delta\downarrow0}\limsup_{N\rightarrow\infty}\frac{1}{N}\log\pp\big(\sup_{t\in[0,\kappa_\gamma (\delta)/N]} \,|\beta(t)|\, >\epsilon\big)
\leq - \lim_{\delta \downarrow 0} \Big\{\frac{\epsilon^2}{2 \kappa_\gamma (\delta)} \Big\} = -\infty\,.
\end{align*}

\noindent \textbf{Step 2.} 
Recall that{,} for any $\alpha_1,\alpha_2 \in M_1(\rr)$,  
\eq\label{eq:dBL}
d_{BL}(\alpha_1,\alpha_2)=\sup_{\|f\|_\infty+\|f\|_{\textrm{Lip}} \le 1}
\{ |(\alpha_1,f) - (\alpha_2,f)| \} \le 2 d_L(\alpha_1,\alpha_2) 
\en
(see \cite[Corollary 11.6.5]{du}). Hence, by \eqref{eq:dL-sup-def} 
there exists 
$r_g:\,(0,\infty)\rightarrow(0,\infty)$ such that $\lim_{\delta\downarrow0} 
r_g(\delta)=0${,} 
and
\eqs
\rho\in B(\gamma,\delta) \quad \Longrightarrow 
\quad \|(\rho,g)-(\gamma,g)\|_\infty
\leq r_g(\delta).
\ens
Recalling that $Z_N^g(t)=(\rho^N,g)(t)$, we thus have from Step 1 that{,} for any $\epsilon>0$, 
\begin{align*}
\lim_{\delta\downarrow0}&\limsup_{N\rightarrow\infty}\frac{1}{N}\log\pp(\rho^N\in B(\gamma,\delta))\\
&=\lim_{\delta\downarrow0}\limsup_{N\rightarrow\infty}\frac{1}{N}\log
\pp\big(\rho^N\in B(\gamma,\delta),\|Z^g_N-Y^g_N\|_\infty \leq\epsilon\big)\\
&\leq\lim_{\delta\downarrow0}\limsup_{N\rightarrow\infty}\frac{1}{N}\log
\pp\big(\|(\rho^N,g)-(\gamma,g)\|_\infty\leq r_g(\delta),\;\|Z^g_N-Y^g_N\|_\infty \leq\epsilon\big)\\
&\leq\limsup_{N\rightarrow\infty}\frac{1}{N}\log
\pp\big(\|(\gamma,g)-Y^g_N\|_\infty \leq 2\epsilon\big)\,.
\end{align*}
Therefore, it suffices to prove the simpler local large deviations upper bound 
\begin{align}\label{eq:simple-ubd}
\lim_{\epsilon \downarrow 0} \limsup_{N\rightarrow\infty}\frac{1}{N}\log
\pp\big(\|(\gamma,g)-Y^g_N\|_\infty \leq \epsilon\big) \le - I^g((\gamma,g)) \,,
\end{align}
provided that (see \eqref{eq:Ig-Phig})
\eq\label{eq:bd-Ig}
I^g((\gamma,g)) \ge (\gamma,g)(T)
- (\gamma,g) (0) - \overline{H}^g(\gamma)(T) 
- \frac{1}{2} \overline{V}^g(\gamma) (T)  \,.
\en

\smallskip
\noindent\textbf{Step 3.} We establish \eqref{eq:simple-ubd} as a 
consequence of the \abbr{LDP} holding for $\{Y^g_N\}$ in the space 
$C([0,T],\rr)$ with a good rate function $I^g(\cdot)$. Indeed,
note that by the time-change formalism for It\^o integrals 
(see e.g. \cite[Theorem 8.5.7]{ok}),
the process $Y^g_N$ can be obtained as $\Psi(N^{-1/2} \tilde{\beta})$ for 
a one-dimensional standard Brownian motion $\tilde{\beta}(t)$, $t \ge 0$ 
and the deterministic operator
$$
\Psi:\;C([0,\infty),\rr)\rightarrow C([0,T],\rr),\quad (\Psi h)(t)=(\gamma,g)(0)+
\overline{H}^g(\gamma)(t)+h(\overline{V}^g(\gamma)(t)). 
$$ 
Clearly, $\Psi$ is continuous with respect to uniform convergence on compacts in
$C([0,\infty),\rr)$, hence by Schilder's theorem (see \cite[Theorem 5.2.3]{dz}),
and the contraction principle (see \cite[Theorem 4.2.1]{dz}), the sequence
$\{Y^g_N\}$ satisfies the \abbr{LDP} in $C([0,T],\rr)$ with the good rate function
\eq\label{eq:Ig-def}
I^g(f)=\inf_{\{h:\,\Psi(h)=f\}} \, \frac{1}{2} \int_0^{S} \Big(\frac{dh}{du}\Big)^2\,\mathrm{d}u\,,
\en
where $S=\overline{V}^g(\gamma)(T)$ and the infimum is over all absolutely continuous 
functions $h$ on $[0,S]$, starting at $h(0)=0$, with Radon-Nikodym derivative 
$\frac{dh}{du} \in L^2([0,S])$. In particular, since
\eqs
0 \le \frac{1}{2} \int_0^S \Big(\frac{dh}{du} -1\Big)^2 \mathrm{d} u 
= \frac{1}{2} \int_0^S \Big(\frac{dh}{du}\Big)^2 \mathrm{d} u - h(S) + \frac{1}{2} S \,,
\ens
it follows from the requirement $(\Psi h)(T)=(\gamma,g)(T)$ that 
\eqs 
I^g((\gamma,g)) \ge h(S) - \frac{1}{2} S = (\gamma,g)(T) - (\gamma,g)(0) 
- \overline{H}^g(\gamma)(T) - \frac{1}{2} S \,,
\ens
which is precisely our claim \eqref{eq:bd-Ig}.
\ep 

\smallskip
We proceed with the local large deviations upper bound for 
paths $\gamma\notin \cA$.
\begin{prop}\label{A_UBD}
If $\gamma\notin \cA${,} then 
\eq\label{infinite_rate}
\lim_{\delta\downarrow0}\limsup_{N\rightarrow\infty}\frac{1}{N}\log\,\pp(\rho^N\in B(\gamma,\delta))=-\infty \,.
\en
\end{prop}

\smallskip
\noindent\textbf{Proof.}
Fixing $N$, let $\qq^{(b)}$ denote the law of the solution of the \abbr{SDS} \eqref{sde} 
with $\qq=\qq^{(0)}$ corresponding to the solution of \eqref{sde} in 
case $b \equiv 0$. Recall that 
\eq\label{eq:Qb-Girsanov}
\frac{\mathrm{d}\qq^{(b)}}{\mathrm{d}\qq}=\exp\big(M^b_{N}(T)-\frac{1}{2}\langle M^b_{N}\rangle(T)\big)
\en
(see \cite[Theorem 3.5.1]{ks}), {with} the continuous martingale
\eqs
M^b_{N}(t)=\sum_{i=1}^N \int_0^t \frac{b(F_{\rho^N(s)}(X_i(s)))}{\sigma(F_{\rho^N(s)}(X_i(s)))}\,\mathrm{d}W_i(s),\quad t\in[0,T] \,,
\ens
whose quadratic variation is
$\langle M^b_{N} \rangle (t) = N \int_0^t U_N^2 (s) \mathrm{d}s$,
with $U_N(s)
$ uniformly bounded by $\|b/\sigma\|_\infty$. Hence, setting
$\kappa=\frac{T}{2} \|b/\sigma\|^2_\infty$, by the Cauchy-Schwarz inequality,
\begin{align}
\qq^{(b)}(\rho^N\in B(\gamma,\delta))
&=\qq\big[e^{M^b_{N}(T)-\frac{1}{2}\langle M^b_{N}\rangle(T)} 
\mathbf{1}_{\{\rho^N\in B(\gamma,\delta)\}}\big]\nonumber \\
&\leq e^{\kappa N} 
\qq\big[e^{2 M^b_{N}(T)-2 \langle M^b_{N}\rangle (T)}\big]^{1/2}\qq\big(\rho^N\in B(\gamma,\delta)\big)^{1/2}\nonumber \\
&= e^{\kappa N}\,\qq^{(2b)}[1]^{1/2} \qq(\rho^N\in B(\gamma,\delta))^{1/2}\,.
\label{eq:drift-removal}
\end{align}
Consequently, 
it suffices to establish 
\eqref{infinite_rate} when $b \equiv 0$ in order to have the same conclusion 
for any other choice of $b(\cdot)$. With $\cA$ independent of such choice, 
we proceed throughout the proof with $b\equiv 0$ (without loss of generality).
In addition, since $\rho^N \to \rho_0$, the bound \eqref{infinite_rate} trivially 
holds when $\gamma(0) \ne \rho_0$. Assuming hereafter that $\gamma(0)=\rho_0$, we
distinguish the three reasons for $\gamma \notin \cA$. We start 
with {\bf case (a)} where $\Rg \notin C_b(\rr_T)$,
proceed to {\bf case (b)}
in which $\Rg \in C_b(\rr_T)$ but
$\gamma$ fails to satisfy the moment condition \eqref{eq:eta-mom} for 
$\iota=\ieta$ of Asmp. \ref{main_ass}(c), and finish with {\bf case (c)}, 
where $\Rg \in C_b(\rr_T)$ and \eqref{eq:eta-mom} holds, but
$t \mapsto (\gamma,g)(t)$ is not absolutely continuous for 
some $g \in \oS$.

\smallskip
\noindent{\bf Case (a).} By assumption $\gamma \in \CC$ and 
$\gamma(0)$ has a density, so 
if $\Rg \notin C_b(\rr_T)$ then necessarily 
$\gamma(s)(\{y\})=3\epsilon$ for some $s \in(0,T]$, $y\in\rr$ 
and $\epsilon>0$.
Fixing $N$ and $0<\delta<\epsilon$, it follows from \eqref{eq:dL-sup-def},
the definition of $d_L(\cdot,\cdot)${,} and the union bound that{,} for 
$m=\lceil N \epsilon \rceil$,
\begin{align} 
\pp(\rho^N\in B(\gamma,\delta)) &\leq \pp\big(\rho^N(s)([y-\delta,y+\delta])\geq\epsilon\big)=\pp\Big(\sum_{i=1}^N \mathbf{1}_{\{|X_i(s)-y|\leq\delta\}}\geq N \epsilon \Big)
\nonumber \\
&\leq\binom{N}{m}\;
\sup_{\{u_j,y_j\}} \; \pp\big(|Z_{j}^{(u_j)}(s)-y_j| \le \delta\,, \quad
\forall j \le m \big)\,.
\label{eq:sup-ui}
\end{align}
Here the supremum is over non-random $\{y_1,\ldots,y_m\}$ (into which we incorporated
the initial conditions $X_i(0)$) and all $[\inf\,\sigma,\sup\,\sigma]$-valued
processes $\{u_1,u_2,\ldots,u_m\}$ adapted to the filtration 
${\HH}_t$ generated by $\{W_i(r), r \in [0,t],1 \le i \le N\}$, while 
\eq\label{eq:Zj-def}
Z_j^{(u_j)}(s)=\int_0^s u_j(r)\,\mathrm{d}W_j(r)\,.
\en
Each It\^o integral in \eqref{eq:Zj-def} is the 
$L^2$-limit of some ${\HH}_t$-adapted stochastic integrals of processes
which are piece-wise constant in time. Therefore, we can and shall assume hereafter that 
$\{u_j\}$ are simple processes, which are constant on each of the 
time intervals $[0,t_1),\dots,[t_{k-1},t_k)$ for some partition 
$0=t_0<t_1<\dots<t_{k}=s$ and $k\in\nn$ (possibly dependent on $N$). 
The probability we maximize in \eqref{eq:sup-ui} is then  
\eq\label{prob_simple_process}
\ev
\big[\prod_{j=1}^m G_\delta (\sum_{\ell=0}^{k-1} u_{j}(t_\ell) \Delta W_{j}(t_\ell),y_j)\big]\,,
\en 
for $\Delta W_i(t_\ell)=W_i(t_{\ell+1})-W_i(t_\ell)$ and 
$G_\delta(z,v)=\mathbf{1}_{|z-v|\le \delta}$. Such expectation, 
when conditioned upon
$$ 
\Gamma_r = \{u_j(t_\ell), \Delta W_i(t_{\ell}), \ell \le k-2\} 
\cup \{u_j(t_{k-1}), \Delta W_{i}(t_{k-1}), i,j \ne r \}
$$
becomes
$$
\prod_{j \ne r}  G_\delta(Z_{j}^{(u_j)}(s),y_j)
\ev\Big[ G_\delta(u_{r}(t_{k-1}) \Delta W_{r}(t_{k-1})
+Z^{(u_r)}_{r}(t_{k-1}),y_r) \, \big| \,\Gamma_r \Big].
$$
The value of $\Delta W_{r}(t_{k-1})$ is independent of everything else,
and clearly all that matters from $\Gamma_r$ to the choice of 
$u_{r}(t_{k-1})$ that maximizes \eqref{prob_simple_process}
is the value of $Z_{r}^{(u_r)}(t_{k-1})$. We thus 
conclude that the optimal $u_r(t_{k-1})$ may be assumed 
measurable 
with respect to the $\sigma$-algebra generated by $W_r(t)$, 
$t\in[0,t_{k-1}]$ and
$u_r(t_\ell), \ell \le k-2$. Substituting optimal $u_j(t_{k-1})$ of this type,
for $1 \le j \le m$, 
changes the function $G_\delta$ considered in \eqref{prob_simple_process},
but retains its form, namely we need thereafter to maximize
\eqs
\ev\Big[\prod_{j=1}^m G'_\delta (Z_j^{(u_j)}(t_{k-1}),y_j)\Big] \,,
\ens
for some (new) function $G'_\delta(z,v)$. The previous argument still applies,
so proceeding by backward induction, from $t_{k-1}$ to $t_{k-2},\ldots,t_0$,
we conclude that it suffices to take the supremum in \eqref{eq:sup-ui} 
only over $\{u_j, j \le m\}$ such that each process $u_j$ is adapted to 
the filtration generated by $W_j$. The bound of \eqref{eq:sup-ui} then becomes
\eq\label{steering_diff}
\pp(\rho^N\in B(\gamma,\delta)) \leq \binom{N}{m}
\big[ \sup_{u_1,y_1} \pp(|Z_1^{(u_1)}(s)-y_1| \le \delta) \big]^m \,,
\en
where the supremum is now over non-random $y_1 \in \rr$ and 
all $[\inf\,\sigma,\sup\,\sigma]$-valued
processes $u_1$ adapted to the filtration generated by $W_1$. Since
$m/N \ge \epsilon$ is bounded away from zero, we get \eqref{infinite_rate} from
\eqs
\limsup_{\delta\downarrow 0} \sup_{u_1,y_1} \pp(|Z_1^{(u_1)}(s)-y_1| \le \delta) = 0 \,,
\ens
which is an immediate consequence of 
the stronger result in \cite[Theorem 1]{McN}. 

\smallskip
\noindent\textbf{Case (b).} We continue with $b \equiv 0$, $\gamma(0,\cdot)=\rho_0${,}
and $\Rg \in C_b(\rr_T)$, whereas \eqref{eq:eta-mom} fails, namely
\eq\label{inf_moment}
\int_0^T (\gamma(t),|x|^{1+\ieta})\,\mathrm{d}t=\infty\,
\en
for $\ieta \in (0,1]$ of Asmp. \ref{main_ass}(c). Let
$0 \le f_K \uparrow f_\infty$ be 
infinitely differentiable functions such that  
\begin{align*}
&f_K(x)=|x|^{1+\ieta}\;\;\mathrm{on}\;\;[-K,-1]\cup[1,K],\quad f_K(x)\leq 1\;\;\mathrm{on}\;\;[-1,1],\\
|f_K'|^2 &
\leq 8 f_K,\;\; \|f_K\|_\infty\leq 2K^{1+\ieta},\;\; \|f'_K\|_\infty\leq 2 K^{\ieta},\;\; \|f''_K\|_\infty\leq 2 
\end{align*}
(we construct $f_K$ by smoothing the 
function $\big(K \wedge (|x| \vee 1) \big)^{1+\ieta}$ 
around the points $\{\pm 1,\pm K\}$). Next, for $Z^{f_K}_N(t):=(\rho^N(t),f_K)$ 
consider the stopping times
\eqs
\tau^K_{N} (r) :=\inf\big\{ t\ge 0 : \; Z_N^{f_K}(t) \geq 2 r \big\} \,.
\ens
Since $f_K(x) \le 2 |x|^{1+\ieta} + 1$ we have from Asmp. \ref{main_ass}(c) that 
for some $C_\ieta$ finite
\eq\label{eq:C-eta}
\sup_{K, N \in\nn} \{Z_N^{f_K}(0)\} \leq 2
\sup_{N\in\nn} (\rho^N(0),|x|^{1+\ieta}) + 1 \le C_\ieta \,.
\en
With $\|f_K''\|_\infty \le 2$ we get upon applying It\^o's formula for $Z_N^{f_K}$
that{,} for the continuous martingale $M^{f_K}_N(t)$ of \eqref{eq:MNg}
and any $r \ge r_0 := C_\ieta + 2 \|A\|_\infty T$,
\eq\label{eq:bd-tauK-tail}
\pp(\tau^K_{N}(r)\leq T)\leq\pp\big( M^{f_K}_N (\tau^K_N (r) \wedge T) \ge r\big)\,.
\en
Further,
since $|f'_K|^2\leq 8f_K$, we have for $V(\cdot)$ of \eqref{eq:Vg} that 
\eqs
\langle M_N^{f_K} \rangle (t) =\frac{1}{N} \int_0^t V^{f_K}(\rho^N)(s) \mathrm{d} s
\le \frac{16 \|A\|_\infty}{N}  \int_0^t Z_N^{f_K}(s) \mathrm{d} s \,.
\ens
In particular, by the 
definition of $\tau^K_N(r)$,  
\eq\label{eq:bd-MfK-tauK}
\langle M_N^{f_K}\rangle (\tau^K_N(r)\wedge T) \le \kappa \frac{r}{N}
\en
for $\kappa = 128 \|A\|_\infty T$. Appealing 
to the martingale representation theorem,  
\eqs
M_N^{f_K}(\tau^K_N(r) \wedge T) \stackrel{d}{=}  
\beta\Big(\langle M_N^{f_K}\rangle (\tau^K_N(r)\wedge T) \Big) \,,
\ens
for some standard Brownian motion $\beta$. Hence, by \eqref{eq:bd-tauK-tail}
and \eqref{eq:bd-MfK-tauK}, for any $r \ge r_0$,
\eq\label{eq:bd-tau-k-r}
\pp(\tau_{N}^K (r) \leq T)
\leq \pp\Big(\sup_{t\in [0,\kappa r/N]} \{ \beta(t) \} \geq r \Big)
\leq 2 \exp\Big(-\frac{N r}{2 \kappa}\Big)\,.
\en
Since $d_L$ is a metric for the weak convergence in $M_1(\rr)$ 
and $f_K \in C_b(\rr)$, by dominated convergence the functionals
$G_K(\xi) := \int_0^T (\xi(t),f_K) \mathrm{d} t$ on $\CC$ are
continuous with respect to the distance $d(\cdot,\cdot)$ 
of \eqref{eq:dL-sup-def}.
Consequently, for any $K \in \nn${,} there exists $\delta_K>0$ such that 
$$
\rho \in B(\gamma,\delta_K) \quad \Longrightarrow \quad 
G_K(\rho) \ge \frac{1}{2} G_K(\gamma) \,.
$$ 
Further, if $G_K(\rho^N) \ge 2 T r${,} then necessarily $\tau_N^K(r) \le T$. 
Thus, setting $r_K=\frac{1}{4T} G_K(\gamma)$, we have{,}
for any $K$ and $\delta \le \delta_K${,} that 
$$
\pp(\rho^N \in B(\gamma,\delta)) \le \pp(G_K(\rho^N) \ge 2 T r_K) \le 
\pp(\tau_N^K(r_K) \le T) \,,
$$
which{,} in view of \eqref{eq:bd-tau-k-r}, yields the bound 
\eqs
\lim_{\delta\downarrow 0}\limsup_{N\rightarrow\infty}
\frac{1}{N}\log \pp\big(\rho^N\in B(\gamma,\delta)\big)
\le - \frac{G_K(\gamma)}{8 \kappa T}
\ens
provided $G_K(\gamma) \ge 4Tr_0$. Our assumption \eqref{inf_moment} and the 
fact that $f_K(x)=|x|^{1+\ieta}$ for all $|x| \in [1,K]$ imply that 
$G_K(\gamma) \to \infty$ as $K \to \infty$, thereby establishing 
\eqref{infinite_rate}.

\smallskip
\noindent\textbf{Case (c).} Steps 1 and 2 of the proof of Prop. \ref{contcenter} 
require only that $\gamma(0)=\rho_0$ and $\Rg \in C_b(\rr_T)$,
both of which hold here. Hence, in view of the derivation leading to \eqref{eq:simple-ubd},
we get \eqref{infinite_rate} as soon as
\eq\label{eq:Ig-inf}
\sup_{g \in \oS} I^g((\gamma,g)) = \infty \,,
\en
for $I^g(\cdot)$ of \eqref{eq:Ig-def}. Further, 
$I^g((\gamma,g))$ is finite only if the identity
\eq\label{eq:contra-abs-cont}
(\gamma,g)(t)=(\Psi h)(t)=(\gamma,g)(0)+\overline{H}^g(\gamma)(t) + 
h\big(\overline{V}^g(\gamma)(t)\big)
\en
holds for some $h$ absolutely continuous. By assumption, there exists
$g\in\oS$ for which the l.h.s. 
of \eqref{eq:contra-abs-cont} is not absolutely continuous on $[0,T]$.
In contrast, both $\overline{H}^g(\gamma)(\cdot)$ and
$\overline{V}^g(\gamma)(\cdot)$ are absolutely continuous for any $g$
and $\gamma$, hence so is the r.h.s. of \eqref{eq:contra-abs-cont} for any absolutely 
continuous $h$, resulting with \eqref{eq:Ig-inf}.
\ep

\smallskip
We conclude this section by showing
the exponential tightness of $\{\rho^N, N\in\nn\}$. 
\begin{prop}\label{exp_tightness_prop} 
Under Asmp. \ref{main_ass}(a)-(c), 
the sequence $\{\rho^N\}$ is exponentially tight on $\CC$.
That is, for any finite $M$ there exists a compact $K_M\subset \CC$ 
for which
\eq
\limsup_{N\rightarrow\infty}\frac{1}{N}\log \pp(\rho^N\notin K_M)\leq -M\,.
\en
\end{prop}

\smallskip
\noindent\textbf{Proof.} 
In view of \eqref{eq:dBL} it suffices to confirm 
the criterion for exponential tightness given in \cite[Lemma A.2]{dza}. 
Specifically, this amounts to showing as
{Step 1}
that $\{\rho^N(t), N\in\nn\}$ is exponentially tight on 
$(M_1(\rr),d_L)$ for each fixed $t \in [0,T]$ (rational), and then proving 
as {Step 2} that{,} for any fixed $\epsilon>0${,} one has
\eq\label{modofcont}
\lim_{\kappa\downarrow0}\sup_{N\in\nn}\frac{1}{N}\log \pp\Big(\sup_{0\leq s,t\leq T,|t-s|\leq\kappa} d_L(\rho^N(s),\rho^N(t))>\delta\Big)=-\infty.
\en 

\smallskip
\noindent\textbf{Step 1.} By Prokhorov's theorem, the set
\eqs
\{\alpha\in M_1(\rr):\;(\alpha,\phi)\leq C\} \,,
\ens 
is pre-compact in $(M_1(\rr),d_L)$ for $\phi(x)=|x|$ and any $C$ finite. 
Hence, to prove our first assertion, it suffices to show that for some
$C=C(M,T)<\infty$ 
\eq\label{exptight2}
\limsup_{N\rightarrow\infty}\frac{1}{N}\log\sup_{t \in [0,T]} 
\pp((\rho^N(t),\phi)>2C)\leq -M\,.
\en
To this end, from \eqref{eq:C-eta} we have that $\sup_N (\rho^N(0),\phi) \le C_0$ finite,
hence with 
\eq\label{eq:ZX-def}
Z_i(t) = \int_0^t \sigma(F_{\rho^N(r)}(X_i(r)))\,\mathrm{d}W_i(r) \,,
\en
it follows by Markov's inequality, that for any $C \ge C_0 + T \|b\|_\infty$
\begin{align}
\pp((\rho^N(t),\phi)>2C) &= \pp\big(\sum_{i=1}^N |X_i(t)| > 2 C N \big) \nonumber \\
&\leq \pp\big(\sum_{i=1}^N |Z_i(t)|>C N \big)
\leq e^{-C N}\, \sup_{\{u_j\}} \, \ev \Big[\prod_{j=1}^N e^{|Z_j^{(u_j)}(t)|}\Big]\,.
\label{eq:bdZX-byZj}  
\end{align}
The supremum here is over the same collection of simple adapted processes $\{u_j\}$ 
we considered in \eqref{prob_simple_process} and with $Z_j^{(u_j)}(t)$ 
as in \eqref{eq:Zj-def}.
By the same argument we have used en-route to \eqref{steering_diff}, it suffices to
consider the situation where each process $u_j$ is adapted to the filtration 
generated by the Brownian motion $W_j$. Consequently, the preceding 
upper bound simplifies to
\eq\label{eq:rho-n-bd}
\pp((\rho^N(t),\phi)>2C) 
\leq e^{-C N}\, \sup_{u_1} \, \ev \big[ e^{|Z_1^{(u_1)}(t)|} \big]^N \,,
\en
where the supremum is now over all $[\inf\,\sigma,\sup\,\sigma]$-valued
processes $u_1$ adapted to the filtration generated by $W_1$. Viewing 
such process $Z_1^{(u_1)}(t)$ as a time changed standard Brownian motion 
$\beta$, results with 
\eqs
\sup_{u_1}\, \ev \big[ e^{|Z_1^{(u_1)}(t)|} \big]
\leq \ev \big[ \exp( \sup_{r\in [0,2 t \|A\|_\infty]} |\beta(r)| \, ) \big] < \infty\,,
\ens 
which together with \eqref{eq:rho-n-bd} yields \eqref{exptight2}. 

\smallskip  
\noindent\textbf{Step 2.} 
First note that for 
$d_{BL}(\cdot,\cdot)$ of
 \eqref{eq:dBL}, $Z_i(t)$ of \eqref{eq:ZX-def}
and any $s,t\in [0,T]$
\begin{align}\label{eq:dl-dbl}
\frac{1}{2} d_L(\rho^N(s)& , \rho^N(t))^2 \leq d_{BL}(\rho^N(s),\rho^N(t))\\
&\leq \frac{1}{N}\sum_{i=1}^N |X_i(t)-X_i(s)| \leq
\frac{1}{N} \sum_{i=1}^N |Z_i(t)-Z_i(s)| + |t-s| \|b\|_\infty 
\nonumber
\end{align}
(see \cite[proof of Theorem 11.3.3]{du} for the left-most inequality). 
Thus, we have \eqref{modofcont}, as soon as we show that{,} for
$$
\mathop{\sf osc} 
(Z;\kappa,T) = \sup_{0\leq s,t\leq T,|t-s|\leq\kappa} \; \{|Z(t)-Z(s)|\} 
$$
and any fixed $\lambda, \delta > 0$,
\eq\label{modofcontclaim}
\lim_{\kappa\downarrow0}\sup_{N\in\nn}\frac{1}{N}\log \pp
\big(\sum_{i=1}^N \mathop{\sf osc} (Z_i;\kappa,T) > \delta N \big) \le 
- \lambda \delta \,.
\en
As for the proof of \eqref{modofcontclaim}, 
similarly to the derivation of \eqref{eq:bdZX-byZj} and 
\eqref{eq:rho-n-bd}, by Markov's inequality we have that  
\begin{align}
\pp (\sum_{i=1}^N \mathop{\sf osc} (Z_i;\kappa,T) > \delta N)
&\le e^{-\lambda \delta N} \sup_{\{u_j\}} \, \ev 
\Big[ \prod_{j=1}^N e^{\lambda \mathop{\sf osc}(Z_j^{(u_j)};\kappa,T)} \Big] 
\nonumber \\
&\le e^{-\lambda \delta N} \Big( \sup_{u_1}
\, \ev \big[e^{\lambda \mathop{\sf osc}(Z_1^{(u_1)};\kappa,T)} \big] \Big)^N \,, 
\label{eq:osc-bd}
\end{align}
where 
the latter supremum is over all $[\inf \sigma, \sup \sigma]$-valued processes
$u_1$ adapted to the filtration generated by $W_1$.
Further, viewing the continous martingale $Z_1^{(u_1)}(t)$ as a time-changed
standard Brownian motion $\beta$, it follows 
by the same reasoning we applied in Step 1 that 
\eq\label{eq:bd-osc-Zu1}
 \sup_{u_1} \, \ev \big[e^{\lambda \mathop{\sf osc}(Z_1^{(u_1)};\kappa,T)} \big] 
\le \ev\big[ e^{\lambda \mathop{\sf osc}(\beta;2 \|A\|_\infty \kappa, 2 \|A\|_\infty T)}
\big] \,.
\en 
Now, since
\eqs
\mathop{\sf osc}(\beta;2 \|A\|_\infty \kappa, 2 \|A\|_\infty T) \le
 2\sup_{r \in [0, 2 \|A\|_\infty T]} \, \{|\beta(r)|\} \,,
\ens
which has finite exponential moments of all order, 
by sample path continuity of $\beta$ and dominated convergence{,} 
the r.h.s. of \eqref{eq:bd-osc-Zu1} decays to one as $\kappa \downarrow 0$.
Thus, combining \eqref{eq:bd-osc-Zu1} and
\eqref{eq:osc-bd} results with \eqref{modofcontclaim}, thereby completing the proof of 
the proposition.
\ep

\section{Proof of Proposition \ref{local_lbd}}\label{sec_tilt}

We work throughout 
under parts (a)--(c) of Asmp. \ref{main_ass} with $\ieta=1$. 

\noindent
{\bf Proof of Part (a)}. Recall from \cite[Section 3]{BFK} that 
in any solution of \eqref{rank_sde}
the ordered particles $X_{(1)}(t) \le X_{(2)}(t) \le \cdots\le X_{(N)}(t)$ 
 satisfy the \abbr{SDS}
\begin{equation}\label{ranks_proc}
\mathrm{d}X_{(j)}(t)=b_j\,\mathrm{d}t+\sigma_j\,\mathrm{d}\beta_j(t)+\frac{1}{2}\,\mathrm{d}\Lambda_{j-1}(t)
-\frac{1}{2}\,\mathrm{d}\Lambda_{j}(t),\quad j=1,\ldots,N\,,
\end{equation}
for independent standard Brownian motions $\{\beta_j\}$ where 
$\Lambda_{0}(t) = \Lambda_N(t) \equiv 0$ 
and $\Lambda_{j}(t)$ for {$j=1,\,\ldots,\,N-1$}  
denotes the local time at zero accumulated 
by the $\rr_+$-valued path $X_{(j+1)}(\cdot)-X_{(j)}(\cdot)$ by time $t$.
It is also shown in \cite[Section 3]{BFK} 
that strong existence and uniqueness holds for the \abbr{SDS} \eqref{ranks_proc}. 
We claim that w.p.1 
\eq
t \mapsto 
\Delta_N (t): = \frac{1}{N} \sum_{j=1}^N \big|X_{(j)}(t)-\widetilde{X}_{(j)}(t)\big| \,,
\en
is non-increasing, for any two strong solutions $X$ and $\widetilde{X}$ 
of \eqref{ranks_proc} driven by the same Brownian motions $\{\beta_j\}$
(extending \cite[inequality (15)]{JR} to arbitrary initial conditions).
Indeed, since $t \mapsto (X_{(j)}(t)-\widetilde{X}_{(j)}(t))$ is of 
bounded variation, its local time process at zero vanishes. Thus,
setting $S_j(t):=\mathrm{sgn}\big(X_{(j)}(t)-\widetilde{X}_{(j)}(t)\big)${,}
we have by Tanaka's formula, followed by summation by parts, that
\begin{align}
\mathrm{d}\Delta_N(t) 
&=\frac{1}{2N}\sum_{j=1}^N S_j(t)\,
\big(\mathrm{d}\Lambda_{j-1}(t)-\mathrm{d}\Lambda_{j}(t)-\mathrm{d}\widetilde{\Lambda}_{j-1}(t) +\mathrm{d}\widetilde{\Lambda}_{j}(t)\big) \nonumber \\
&= 
\frac{1}{2N}\sum_{j=2}^N \big(S_j(t)-S_{j-1}(t)\big)\,
\big(\mathrm{d}\Lambda_{j-1}(t)-\mathrm{d}\widetilde{\Lambda}_{j-1}(t)\big)\,.
\label{eq:Delta_N-bd}
\end{align}  
Since $X_{(j)}(t)=X_{(j-1)}(t)$ at times of increase of $\Lambda_{j-1}(t)$ 
and $\widetilde{X}_{(j)}(t)=\widetilde{X}_{(j-1)}(t)$ at times of increase
of $\widetilde{\Lambda}_{j-1}(t)$, it is easy to check that{,} for 
{$j=2,\,\ldots,\,N$} and all $t \ge 0$,
$$
\big(S_j(t)-S_{j-1}(t)\big)    \mathrm{d}\Lambda_{j-1}(t) \le 0 \le 
\big(S_j(t) - S_{j-1}(t) \big) \mathrm{d}\widetilde{\Lambda}_{j-1}(t)\,.
$$
Hence, by \eqref{eq:Delta_N-bd} we have, as claimed, that $\mathrm{d}\Delta_N(t) \le 0$.
Next, with $\rho^N$ and $\widetilde{\rho}^N$ denoting the 
paths of empirical measures of $X$ and $\widetilde{X}$, respectively, 
we see that for any $t \ge 0$,
\eq\label{eq:wass-bd}
d_{BL}(\rho^N(t),\widetilde{\rho}^N(t)) \le \Delta_N(t) \le \Delta_N (0) =  
W_1(\rho^N(0),\widetilde{\rho}^N(0)) \,,
\en
where $W_1(\cdot,\cdot)$ stands for the $L_1$-Wasserstein distance 
on $M_1(\rr)$:
$$
W_1(\alpha_1,\alpha_2):=\inf\{ \, \ev [\,|Y_1-Y_2|\,] \, : Y_1 \sim \alpha_1, 
Y_2 \sim \alpha_2 \} \,.
$$
Further, as $\frac{1}{2} d_L(\alpha_1,\alpha_2)^2 \le
d_{BL}(\alpha_1,\alpha_2)$ (see \eqref{eq:dl-dbl}), it follows from \eqref{eq:wass-bd}
that 
$$
\frac{1}{2} d(\rho^N,\widetilde{\rho}^N)^2 \le W_1(\rho^N(0),\widetilde{\rho}^N(0)) \,.
$$
Fixing $\ell,\eps$ we let $\rho^{N,\ell,\eps}(0)$ 
be the empirical measure of $\widetilde{X}_{(j)}(0)=F^{-1}(u_{j,N})$,
{with} $u_{j,N}=j/(N+1)$ and $F:=F_{\gamma^{\ell,\eps}(0)}$ 
a continuous \abbr{cdf}. By the preceding, we get 
\eqref{eq:coupling} upon showing that 
\eq\label{eq:bd-delta-N-0}
\limsup_{\eps \to 0}\limsup_{\ell \to \infty} \limsup_{N \to \infty}
W_1(\rho^{N}(0),\rho^{N,\ell,\eps}(0)) 
= 0 \,.
\en
We then complete the proof of part (a)
by observing that{,} 
\begin{align} 
\frac{N}{N+1}( |x|,  \rho^{N,\ell,\eps}(0) ) &= 
\frac{1}{N+1} \sum_{j=1}^N |F^{-1}(u_{j,N})| \nonumber \\
&\le \int_0^1 |F^{-1} (u)| \mathrm{d} u = 
( |x|, \gamma^{\ell,\eps}(0)) \,,
\label{eq:bd-int-rho-ell-eps}
\end{align}
which is finite since $\gamma^{\ell,\eps}(0) \in \Mstar$, 
whereas for $N \to \infty$ and any fixed $x \in \rr$, 
$$
\rho^{N,\ell,\eps}(0) \big( (-\infty,x] \big)= \frac{1}{N} 
\Big\lfloor (N+1) F(x) \Big\rfloor \to F(x) \,.
$$
That is, as claimed our $\{\rho^{N,\ell,\eps}(0)\}$ satisfy 
Asmp.~\ref{main_ass}(c) (with $\ieta=0$). In particular, 
\eq\label{eq:first-bd}
\lim_{N \to \infty} d_{BL}(\rho^{N,\ell,\eps}(0),\gamma^{\ell,\eps}(0)) = 0 \,.
\en
Turning to prove 
\eqref{eq:bd-delta-N-0}{,} recall that{,} for any $\kappa \ge 1${,}  
\eq\label{eq:was-bd}
W_1(\alpha_1,\alpha_2) \le \kappa d_{BL}(\alpha_1,\alpha_2) 
+ 3 (|x| {\bf 1}_{|x| > \kappa}, \alpha_1+\alpha_2) 
\en
(for example, this follows from \cite[Theorem 11.8.2]{du}).
From Asmp.~\ref{main_ass}(c) we know that 
$d_{BL}(\rho^N(0),\gamma(0)) \to 0$ when $N \to \infty$, whereas 
$d_{BL}(\gamma(0),\gamma^{\ell,\eps}(0)) \to 0$
when $\ell \to \infty$ followed by $\eps \to 0$
(see \eqref{eq:d-eta-conv}). Combining these facts with 
\eqref{eq:first-bd}, \eqref{eq:was-bd}{,} and the triangle inequality for 
$d_{BL}(\cdot,\cdot)$, we get \eqref{eq:bd-delta-N-0} by showing that 
\begin{align}\label{eq:ui-rho-N}
\lim_{\kappa \to \infty} \limsup_{N \to \infty} \;
(|x| {\bf 1}_{|x| > \kappa}, \rho^N(0)) &= 0 \,, \\
\lim_{\kappa \to \infty} \limsup_{\eps \to 0} \limsup_{\ell \to \infty} 
\limsup_{N \to \infty} \;
(|x| {\bf 1}_{|x| > \kappa}, \rho^{N,\ell,\eps} (0)) &= 0 \,.
\label{eq:ui-rho-N-ell-eps}
\end{align}
With $\sup_N (|x|^2,\rho^N(0))$ finite (see Asmp.~\ref{main_ass}(c), $\ieta=1$), 
we obviously have \eqref{eq:ui-rho-N}. As for \eqref{eq:ui-rho-N-ell-eps}, from 
\eqref{eq:first-bd} and $\gamma^{\ell,\eps}(0)$ having a density, we
deduce that as $N \to \infty$
$$
(|x| {\bf 1}_{|x| \le \kappa},\rho^{N,\ell,\eps}(0)) \to 
(|x| {\bf 1}_{|x| \le \kappa},\gamma^{\ell,\eps}(0)) \,.
$$
Further, when $\ell \to \infty$, $\eps \to 0${,} and then $\kappa \to \infty$, we
have by \eqref{eq:iota-ubdd} that
$$
(|x| {\bf 1}_{|x| > \kappa}, \gamma^{\ell,\eps} (0))
\le \frac{1}{\ell} (|x|,\widehat{\mu}^\eps) +
\frac{1}{\kappa}(|x|^2,\gamma^\eps(0)) \to 0 \,,
$$ 
which together with \eqref{eq:bd-int-rho-ell-eps} yields 
\eqref{eq:ui-rho-N-ell-eps}.
\ep
 
\bigskip
\noindent
{\bf Proof of Part (b).} From part (a) we know 
that Asmp.~\ref{main_ass}(c) holds for $\ieta=0$, 
$\rho^{N,\ell,\eps}(0)${,} and $\rho_0^{\ell,\eps}=\gamma^{\ell,\eps}(0)$,
so we simplify notation by dropping hereafter the 
superscripts $\ell,\eps$. That is, we fix 
$\gamma \in \GG$ with $\widetilde{J}(\gamma)<\infty$ 
(see Definition \ref{new-def-eta} and \eqref{eq:J-def}, 
respectively), having $R=\Rg \in 
C_b^{\infty}(\rre) \bigcap \FF_{3/2}$ (for $\FF_q$ of \eqref{eq:thm-13}), 
starting at $R(0,x)=F_{\rho_0}(x)${,} with $R_x$ strictly positive 
on $\rre$ {and} such that 
\eq\label{eq:h-def}
h = \frac{R_t-(A(R)R_x)_x}{A(R)\,R_x} \in C_b(\rre) 
\en
with $x \mapsto h(t,x)$ uniformly 
Lipschitz continuous 
on $\rr_T$. The functional on $\CC$
\eq
\widehat{J}(\rho) = \frac{1}{2} \int_0^T (\rho,U(\rho)^2)(s) \mathrm{d}s \,,
\en 
{with} the continuous in time, bounded function on $\rr_T$
$$
U(\rho):= \frac{h A(\Rro) + b(\Rro)}{\sigma(\Rro)}\,,
$$ 
is then such that  
${\widetilde J}(\gamma)=\widehat{J}(\gamma) < \infty$.
\newline
To prove the local large deviations lower bound 
\eqref{eq:ldp-loc-lbd} we introduce in {Step 1} a suitable change 
of measure to $\qq_1$ for which the relevant 
Radon-Nikodym derivative is at least $e^{-N (\widehat{J}(\gamma)+\epsilon)}$
on the event $\{\rho^N \in B(\gamma,\delta)\}$ (when $N \to\infty$ followed
by $\delta \downarrow 0$), then verify in {Step 2}
the relevant \abbr{lln} \eqref{lln_claim} for 
$\rho^N$ under $\qq_1$. This step relies on proving 
in Lemma \ref{no_atoms_lemma} that any limit point $\widehat{\gamma}$ of $\rho^N$
has no atoms, from which we deduce that the corresponding path of \abbr{cdf}-s 
$w=\Rwg$
solves the weak form \eqref{eq:lln-weak2-form} 
of the porus medium equation \eqref{PDEtilt}. The tilt $h$ in the latter equation 
is given by \eqref{eq:h-def}, so $\Rg$ is one such solution,
and the analysis of Lemma \ref{PDE_uniqueness} guarantees its uniqueness.

\smallskip
\noindent\textbf{Step 1.} Since $x \mapsto h(t,x)$ is uniformly 
bounded, Lipschitz function on $\rr_T$,
there exist{s}{,} for any $q \in \rr${,} a probability measure $\qq_q$ under which{,} for 
$i=1,\,\ldots,\,N$,
\begin{align}
X_i(t)=X_i(0)
+ \int_0^t \psi_q (\rho^N)(s,X_i(s)) \,\mathrm{d}s 
+\int_0^t \sigma(F_{\rho^N(s)}(X_i(s)))\,\mathrm{d} W^{q}_i (s),
\label{Q_dynamics}
\end{align}
with $\psi_q(\rho)=b(\Rro) - q U(\rho) \sigma(\Rro)$
and $\{W^q_i\}$ independent standard Brownian motions. 
Further, $\pp=\qq_0$ and{,} similarly to \eqref{eq:Qb-Girsanov} 
(which has $h\equiv 0$, i.e. $\psi_q=(1-q) b(\Rro)$),
by Girsanov's theorem  
\eqs
\frac{\mathrm{d}\qq_q}{\mathrm{d}\pp}=
\exp\big( -q M_N(T)-\frac{q^2}{2} \langle M_N \rangle (T)\big),
\ens
{with} the continuous martingales 
\eqs
M_N(t)=\sum_{i=1}^N \int_0^t U(\rho^N)(s,X_i(s))
\,\mathrm{d} W^{0}_i(s)
\ens
{satisfying} $\langle M_N \rangle (T) = 2 N \widehat{J}(\rho^N)$
(see \cite[Theorem 3.5.1]{ks}).
By the triangle inequality,
for any $\rho \in B(\gamma,\delta)$,  
\begin{align}
2 
\big| \widehat{J}(\rho) - \widehat{J}(\gamma)
\big| &\leq \sup_{\rho\in B(\gamma,\delta)}
\Big| \int_0^T (\rho,U(\rho)^2-U(\gamma)^2) (s) \mathrm{d}s \Big|
 \nonumber \\
&+\sup_{\rho\in B(\gamma,\delta)}\Big|
\int_0^T (\rho,U(\gamma)^2)(s) \mathrm{d} s - \int_0^T (\gamma,U(\gamma)^2)(s)
 \mathrm{d} s \Big|
\label{eq:til-J-bd}
\end{align}
Since $b$ and $\sigma$ are bounded, 
Lipschitz functions, 
and $h$, $\sigma^{-1}$ are uniformly bounded, it is easy 
to check that 
$U(\rho)^2$ is a
bounded, Lipschitz function of $\Rro$. Consequently, 
$|U(\rho)^2-U(\gamma)^2| \le C |\Rro-\Rg|$ for 
some finite $C$ and all $(t,x) \in \rr_T$. Thus, from
\eqref{R_comp_est} it follows that the first term 
on the r.h.s. of \eqref{eq:til-J-bd} converges 
to $0$ as $\delta\downarrow 0$. Further, 
with $h$ and $\Rg$ in $C_b(\rr_T)$, also $U(\gamma)^2 \in C_b(\rr_T)$.
Hence, by dominated convergence, the functional
$\xi \mapsto \int_0^T (\xi,U(\gamma)^2)(s) \mathrm{d} s$  
is continuous on $\CC${,} and the second term on the r.h.s.
of \eqref{eq:til-J-bd} also converges to $0$ as $\delta\downarrow 0$.
We conclude that{,} for any fixed $\epsilon>0$ and all $\delta<\delta_0(\epsilon)$,
\eq
\rho^N \in B(\gamma,\delta) \quad \Longrightarrow \quad
\big|\frac{1}{2} \langle M_N \rangle (T)  - N \widehat{J} (\gamma)\big|
\leq \epsilon N \,.
\en 
For such $\delta<\delta_0(\epsilon)${,} any
$p=q/(q-1)>1$, $q>1${,} and all $N \in \nn$,
we thus have by H\"older's inequality that 
\begin{align*} 
\qq_1(\rho^N\in B(\gamma,\delta))=\pp\Big[e^{-M_N(T)-\frac{1}{2} 
\langle M_N\rangle (T)}\mathbf{1}_{\{\rho^N\in B(\gamma,\delta)\}}\Big]&\\
\leq\pp \big[e^{-q M_N(T) - \frac{q}{2} \langle M_N \rangle (T)}
\mathbf{1}_{\{\rho^N\in B(\gamma,\delta)\}}
\big]^{1/q}&\,\pp\big(\rho^N\in B(\gamma,\delta)\big)^{1/p}\\
\le e^{(q-1)N(\widehat{J}(\gamma)+\epsilon)}
\;\pp \big[e^{-q M_N(T)-\frac{q^2}{2} \langle M_N \rangle (T)} \big]^{1/q} 
&\;\pp \big(\rho^N\in B(\gamma,\delta)\big)^{1/p}\\
= e^{(q-1)N(\widehat{J}(\gamma)+\epsilon)} \; \big[ \qq_q(1) \big]^{1/q}
& \;\pp \big(\rho^N\in B(\gamma,\delta)\big)^{1/p}\,.
\end{align*}
Consider $\frac{p}{N}\log(\cdot)$ of both sides, taking first
$N \to \infty$ followed by $\delta \downarrow 0$, to find that  
$$
\lim_{\delta\downarrow0}\liminf_{N\rightarrow\infty}\frac{1}{N}\log \pp(\rho^N\in B(\gamma,\delta))\geq -q (\widehat{J} (\gamma) + \epsilon) \,,
$$
for any $\epsilon>0$, $q>1$ (hence 
\eqref{eq:ldp-loc-lbd} holds), provided that 
\eq
\lim_{N\rightarrow\infty}\qq_1\big(\rho^N\in B(\gamma,\delta)\big)=1 \label{lln_claim}
\,.
\en

\smallskip
\noindent\textbf{Step 2.} To prove \eqref{lln_claim}, recall that $\qq_1$ corresponds
to $\{X_i(t)\}$ of \eqref{Q_dynamics} for drift $\psi_1(\rho)= - h A(\Rro)$ 
and{,} 
for each $N\in\nn$, let $Q_{N}$ denote the law of $\rho^N$ under $\qq_1$. 
One then deduces the uniform tightness, and hence pre-compactness, of 
the collection $\{Q_{N}\}$ in the space of 
probability measures on $\CC$.
Indeed, Steps 1 and 2 of the proof of \cite[Theorem 1.1]{sh} 
rely only on a general compactness criterion for subsets of $\CC$ 
(taken from \cite[Lemma 1.3]{ga}), so can be carried out mutatis mutandis,
appealing here to boundedness of the drift and diffusion coefficients {in} \eqref{Q_dynamics} 
(in place of boundedness of the corresponding coefficients in the dynamics treated in \cite[Theorem 1.1]{sh}). Moreover, since $\rho^N(0) \to \rho_0=\gamma(0,\cdot)$,
the computations there which involve the initial conditions, can be omitted here.
To prove \eqref{lln_claim}, it thus suffices to show that the atomic measure
$\delta_\gamma$ is the only possible limit point of the sequence $\{Q_N\}$. 
Alternatively, passing to the relevant sub-sequence and utilizing the 
Skorokhod representation theorem in the form of \cite[Theorem 3.5.1]{du2}, we 
can and shall assume that the variables $\rho^N$, $N\in\nn$ are defined on the 
same probability space and converge almost surely in $\CC$,
when $N \to \infty${,} to some limiting variable $\widehat{\gamma}$. Thus, the task of 
proving \eqref{lln_claim} amounts to showing that $\widehat{\gamma}=\gamma$ w.p.1.
\newline
To this end, fixing $g \in \oS$ and
replacing hereafter 
$b(\Rxi)$ by $-h A(\Rxi)$ in the definition of 
${\RR}^\xi g$, note that the l.h.s.
of the identity \eqref{N_mckean_vlasov} converges w.p.1 as $N \to \infty$,
to $(\widehat{\gamma},g)(t)-(\rho_0,g(0))$, for all $t\in[0,T]$. 
We claim that w.p.1 the r.h.s. of \eqref{N_mckean_vlasov} converges to 
$\int_0^t H^g(\widehat{\gamma}) (s) \mathrm{d} s$ and consequently for all
$t \in [0,T]$,
\eq\label{eq:lln-weak-form}
(\widehat{\gamma},g)(t)-(\rho_0,g(0)) =
\int_0^t (\widehat{\gamma}, g_t + (g_{xx} - h g_x) A(\Rwg) ) (s)
\,\mathrm{d}s\,.
\en
Indeed, recall that $\overline{V}^g(\rho^N)(T)$ is uniformly 
bounded by $C=2 T \|A g_x^2\|_\infty$ finite,
so $\langle M^g_N \rangle (T) \le C/N$ and w.p.1
the continuous martingale $M_N^g(t) \to 0$ uniformly over
$[0,T]$ (for example, combine the Burkholder-Davis-Gundy inequalities, 
see \cite[Theorem 3.3.28]{ks}, and the Borel-Cantelli lemma). 
Further, we show in Lemma \ref{no_atoms_lemma} that w.p.1. 
$\Rwg \in C_b(\rr_T)$, from 
which one deduces as in the derivation of \eqref{R_comp_est}, that w.p.1
\eqs
\lim_{N \to \infty} \int_0^T (\rho^N,|F_{\rho^N}- \Rwg|)(s)\, \mathrm{d} s = 0 \,.
\ens
Since $\Rwg \in C_b(\rr_T)${,} also 
${\RR}^{\widehat{\gamma}} g \in C_b(\rr_T)$. Thus, 
with $h g_x$ and $g_{xx}$ bounded, the preceding convergence to zero implies,
as in the derivation of \eqref{iden1}, that w.p.1
\eqs
\lim_{N\rightarrow\infty} \int_0^T
\big|H^g(\rho^N)(s)- H^g(\widehat{\gamma})(s)\big|\,\mathrm{d}s = 0 \,,
\ens
thereby completing the proof of \eqref{eq:lln-weak-form}. Now,
setting $w=\Rwg$, and 
having \eqref{eq:lln-weak-form} hold w.p.1 for all 
$g$ in a countable dense subset of $\oS$, 
we deduce that \eqref{eq:lln-weak-form}
holds for all $g \in \oS$,
which amounts after integration by parts over $\rr$ to
\eq\label{eq:lln-weak2-form}
\int_{\rr} (f w) (t,x) \mathrm{d} x 
- \int_{\rr} (f w) (0,x) \mathrm{d} x =
\int_{\rr_t} \big[ f_t w +  f_{xx} \Sigma(w)  - (h f)_x \Sigma(w) \big] \, \mathrm{d} m \,,
\en
holding for $w(0,x)=F_{\rho_0}(x)$ and all $f=g_x \in \DxS$. Here 
$\Sigma(w)=\int_0^w A(r) \mathrm{d} r$
and $h_x$ is well-defined ($m$-a.e.), since $x \mapsto h(t,x)$ is a Lipschitz function.
Note that{,} for $w \in C^{1,2}(\rr_T)$, further integration by parts
(to eliminate all derivatives of $f$) confirms that this is equivalent 
to $w$ solving the porous medium equation \eqref{PDEtilt}
(with initial condition $F_{\rho_0}$).  
In view of \eqref{eq:h-def} 
one such solution is $w=\Rg$, hence w.p.1
$\widehat{\gamma}=\gamma$ provided we establish the uniqueness of 
such generalized solution.
In conclusion, all that remains for establishing Prop.
\ref{local_lbd} is to prove the following two lemmas.
\begin{lemma}\label{no_atoms_lemma}
Consider the unique weak solution of \eqref{Q_dynamics} 
for $\psi_1(\rho)=-h A(\Rro)$, and suppose that 
$\widehat{\gamma}$ is an a.s. limit point in 
$\CC$ of $\rho^N$. Then, w.p.1 
the probability measures $\widehat{\gamma}(t,\cdot)$ 
have no atoms for all $t \in [0,T]$.
\end{lemma}
\begin{lemma}\label{PDE_uniqueness} 
Suppose $h \in C_b(\rre)$ with $x \mapsto h(t,x)$ uniformly Lipschitz 
on $\rre$, for which \eqref{PDEtilt} has a bounded classical solution 
$R=u \in C^{1,2}(\rr_T)$ with $u(0,x)=F_{\rho_0}(x)$. 
It is then the only solution $w \in C_b(\rr_T)$ of 
\eqref{eq:lln-weak2-form} with such initial conditions, 
for which $w(t,\cdot)$ are \abbr{cdf}-s 
of some path in $\CC$.
\end{lemma}

\noindent\textbf{Proof of Lemma \ref{no_atoms_lemma}.} 
When proving Prop. \ref{A_UBD}, we first removed the drift
$b(\Rro)$ thanks to the bound \eqref{eq:drift-removal}, which 
applies for any uniformly bounded drift. So, using the same argument
for the measure $\qq_1$ with bounded drift $-h A(\Rro)$, we conclude
as in case (a) of the proof of Prop. \ref{A_UBD} that{,} 
for any $\epsilon>0$ and some $\delta=\delta(\epsilon) \in (0,\epsilon)$,  
\eq\label{delta_property1}
\limsup_{N\rightarrow\infty}\frac{1}{N}\log
\sup_{s \in [0,T], y \in \rr} \qq_1 \big(\rho^N(s)([y-3\delta,y+3\delta]) \ge \epsilon) 
\le -1 \,.
\en  
Similarly, Step 1 and Step 2 of the proof of 
Prop. \ref{exp_tightness_prop} apply whenever the drift is uniformly 
bounded, hence{,} for any $\delta,\epsilon>0$, some 
$M=M(\epsilon)$ finite{,} and $\kappa(\delta)>0$, 
\begin{align}
\label{delta_property2}
\limsup_{N\rightarrow\infty}\frac{1}{N}\log & \sup_{s \in [0,T]} \qq_1
\big(\rho^N(s)([-M,M]) \le 1- \epsilon \big) \le -1 \,, \\
\limsup_{N \to \infty} \frac{1}{N} \log & \qq_1 \big(
\sup_{0\leq s,t\leq T,|t-s|\leq\kappa} d_L(\rho^N(s),\rho^N(t))>\delta
\big) \le -1
\label{delta_property3}   
\end{align}     
(see \eqref{exptight2} and \eqref{modofcont}, respectively). 
Fixing $\delta(\epsilon)>0$ and $M(\epsilon)$ finite, then 
$\kappa(\delta)>0$, consider \eqref{delta_property1} and 
\eqref{delta_property2} at all $\{s_i, i \le m\}$ on a finite $\kappa$-net 
in $[0,T]$ and all $\{y_j, j \le \ell\}$ in a finite 
$\delta$-net within $[-(M+2\delta),(M+2\delta)]$ 
to conclude by \eqref{delta_property3} and the Borel-Cantelli lemma that w.p.1 for all $N$ large enough,
\begin{align*}
& \sup_{1 \le i \le m} \, \sup_{y\in\rr} \, 
\{\rho^N(s_i)([y-2\delta,y+2\delta]) \} \le \epsilon \,,\\
& \sup_{t \in [0,T]} \, \inf_{1 \le i \le m} \, \{ d_L(\rho^N(s_i),\rho^N(t)) \} \le \delta \,.
\end{align*}
Consequently, as in \eqref{eq:dl-r},
$$
\sup_{t \in [0,T]}\, \sup_{y \in \rr} \,\{
\rho^N(t)([y-\delta,y+\delta]) \} \le 3\epsilon \,.
$$
Combining this with the Portmanteau theorem, we infer that, for every $\epsilon>0$, 
there exists a $\delta>0$ such that for any limit point $\widehat{\gamma}$
of $\rho^N$, w.p.1
\begin{align*}
\sup_{t\in[0,T]}\sup_{y\in\rr} \;\widehat{\gamma}(t)\big((y-\delta,y+\delta)\big)
&\leq\sup_{t\in[0,T]}\sup_{y\in\rr}\;\liminf_{N\rightarrow\infty}\;\rho^N(t)
\big((y-\delta,y+\delta)\big)\\
&\leq \limsup_{N\rightarrow\infty}\sup_{t\in[0,T]}\sup_{y\in\rr}\;\{\rho^N(t)
\big((y-\delta,y+\delta)\big)\}\leq3\epsilon\,.
\end{align*}
This shows that w.p.1 the path of measures 
$\widehat{\gamma}$ has no atoms of mass at least 
$3 \epsilon$, and taking $\epsilon \downarrow 0$ finishes the proof of 
the lemma. \ep

\smallskip
\noindent\textbf{Proof of Lemma \ref{PDE_uniqueness}.} First{,} note 
that if \eqref{eq:lln-weak2-form} holds for $w(t,\cdot)$ which are the 
\abbr{cdf}-s of some path in $\CC$, and all
$f \in \DxS${,} then it further holds for all
$\widehat{f} \in \oS$. Indeed, for any 
such $\widehat{f}$ there exist $f_k \in \DxS$ 
coinciding with $\widehat{f}$ on $[0,T] \times [-k,\infty)$ such that
both $\sup_k \|f_k - \widehat{f} \|_{W_1^{1,2}(\rr_T)}$ 
and $\sup_{k,t} \|f_k(t,\cdot)-\widehat{f}(t,\cdot)\|_{L^1(\rr)}$ 
are finite.
Thus, with $x \mapsto h(x,t)$ uniformly bounded, globally Lipschitz on $\rr_T$ and 
$\Sigma(w) \le \|A\|_\infty w$, it follows from 
\eqref{eq:Runif-tight} that the value each side of \eqref{eq:lln-weak2-form} 
takes for $f_k$, converges as $k \to \infty$ to its value 
for $\widehat{f}$, thereby extending the scope of 
\eqref{eq:lln-weak2-form} to all of $\oS$. The 
latter identity involves only $(f,f_t,f_x,f_{xx})$, 
hence holds for any $f \in C_c^{1,2}(\rr_T)$. 
Further, setting $\mathbb{K}_{r}=[0,t] \times [-r,r]$,
the identity \eqref{eq:lln-weak2-form} applies 
for any $f \in C^{1,2}(\mathbb{K}_{r})$ such that 
$f(s,\pm r)=0$ for all $s \in [0,t]$, provided 
one adds to its l.h.s. the boundary term  
$$
\int_0^t \big[(\Sigma(w) f_x) (s,r^-) - (\Sigma(w) f_x)(s,-r^-)\big] \mathrm{d} s 
$$
(which takes into account the jump-discontinuity 
of $f_x$ at $\partial \mathbb{K}_r$). That is, 
$w(t,x)$ is a generalized solution of the 
(CP) problem as in \cite[Definition 1.1]{dk},
except for replacing the term $f_x b(w)$ there by $(h f)_x \Sigma(w)$. 
The uniqueness of such non-negative
$w \in C_b(\rr_T)$, starting at the non-negative 
$w(0,x)=u(0,x) \in C_b(\rr)${,} thus follows by 
adapting the proof 
of \cite[Theorem 4.2 (1)]{dk} to handle $h \not\equiv 1$. 

To this end, we modify hereafter 
\eqref{eq:lln-weak2-form} as above and 
prove the analog of \cite[(4.2)]{dk}.
That is, we call $w$ a sub-solution if the r.h.s. of \eqref{eq:lln-weak2-form} 
is greater or equal its (modified) l.h.s. for every 
non-negative $f \in C^{1,2}(\mathbb{K}_r)$, any $r>0${,} and all $t \in (0,T]$,
while $w$ is a super-solution when the corresponding l.h.s. is greater
or equal the r.h.s. for any such $f,r,t$. It then suffices to show that for some 
$c=c(h)$ finite, 
any super-solution $w$, all $t\in[0,T]$, $\ell \in (0,\infty)${,} and $[0,1]$-valued 
$\omega \in C_c^\infty([-\ell,\ell])$,
\eq\label{uniq_refined}
\int_\rr (u(t,x)-w(t,x)) \omega(x)\,\mathrm{d}x\leq c\int_\rr
(u(0,x)-w(0,x))_+ \,\mathrm{d}x\,,
\en
where for any sub-solution $w$, the same inequality holds 
with $u-w$ replaced by $w-u$. Indeed, by definition of a super-solution 
(or sub-solution), the r.h.s. of \eqref{uniq_refined} is zero. Hence, 
choosing $\omega(x)$ as smooth approximations
of $\mathbf{1}_{u(t,x)>w(t,x)}$ on $(-\ell,\ell)$ and sending $\ell \to \infty$, we 
deduce that $m$-a.e. $u \le w$ for super-solutions and $u \ge w$ for 
sub-solutions, yielding the stated uniqueness of the solution $w$. 
Fixing $t,\ell${,} and $\omega$, we prove \eqref{uniq_refined} 
for a given super-solution $w$ (exchanging the roles of $u$ and $w$ then yields 
the proof for sub-solutions). As in the proof of 
\cite[(4.2)]{dk}, for $f=f^{(n,r)}$, $n, r \ge \ell+1$
which solve suitable linear parabolic first boundary value problems,
we bound the difference between the left most terms of \eqref{eq:lln-weak2-form} 
for $u$ and $w$. Taking $n \to \infty$ followed by $r \to \infty$ then 
yields \eqref{uniq_refined}. 

Specifically, 
note that $\Sigma(u)-\Sigma(w) = (u-w) A_\star$ where 
\begin{align*}
A_\star = \int_0^1 A(\tau u+(1-\tau)w)\,\mathrm{d}\tau
\end{align*} 
is in $C_b(\rr_T)$ and $A_\star \ge \underline{a}>0$. Hence, 
there exist uniformly bounded smooth functions $(A_{n},B_{n},C_{n})$, 
such that 
$A_{n} \downarrow A_\star$ and $B_{n}\downarrow B_\star = -h A_\star \in C_b(\rr_T)$, 
uniformly on $\rr_T$, whereas $m$-a.e. $C_{n} \to C_\star = h_x A_\star$.
Now{,} consider{,} for each $n$ and $r \ge \ell+1$,
the unique classical solution $f^{(n,r)} \in C^{1,2}(\mathbb{K}_r^o)$ of 
the first boundary value problem
\begin{align}
\LL_n f := f_t+A_{n} f_{xx}+B_{n} f_x - C_{n} f &=0
\quad \;\; \mathrm{on}\;\; (0,t)\times(-r,r),\nonumber \\
f(t,x)&=\omega(x)\,,\quad x\in[-r,r],
\label{linear_PDE}
\\
f(s,-r)=f(s,r)&=0\,,\qquad\;\; s\in[0,t]\,.\nonumber
\end{align}
With $A_n$ bounded away from zero, the existence and 
uniqueness of such solution for \eqref{linear_PDE} is 
well-known (see e.g. \cite[chapter IV, Theorem 5.2]{lsu}),
and further, with $\omega \ge 0$ also $f^{(n,r)} \ge 0$.
Next, setting $g^{\pm}_{n,r}(s) = [f^{(n,r)}_x (\Sigma(u)-\Sigma(w))](s,\pm r^-)$, we
use the test function $f=f^{(n,r)}$ in our modified
\eqref{eq:lln-weak2-form}, to bound the l.h.s. of \eqref{uniq_refined} by
\begin{align}\label{eq:ubd-nr}
& \int_{-r}^r (u(0,x)-w(0,x))_+ f^{(n,r)}(0,x) \mathrm{d} x -
\int_0^t [g^+_{n,r}(s) - g^-_{n,r}(s)]\,\mathrm{d} s \nonumber \\
& + 
\int_{\mathbb{K}_r^o} (u-w) \big[f_{xx}^{(n,r)} (A_\star - A_n) + f_x^{(n,r)} (B_\star-B_n) - f^{(n,r)} 
(C_\star -C_n) \big] \, \mathrm{d} m \,.
\end{align}
Recall the uniform boundedness of $u-w$ and the uniform convergence to zero of 
$A_\star - A_n$ and $B_\star - B_n$. 
Thus, similarly to the derivation of \cite[(4.12)]{dk}, 
the last term of \eqref{eq:ubd-nr} goes to zero when $n \to \infty$, provided
$\sup_{s,x,n,r} f^{(n,r)}(s,x) \le c$ finite and both 
$e_n := \|f^{(n,r)}_{xx}\|_{L^2(\mathbb{K}_r)}$ and 
$\chi_n := \|f^{(n,r)}_x\|_{\infty}$ are uniformly bounded in $n$. Taking then 
$r \to \infty$ yields \eqref{uniq_refined}, if in addition  
$\sup_{s,n} |g^{\pm}_{n,r}(s)| \to 0$ as $r \to \infty$. 
\newline
Turning to prove the latter four estimates, 
assume first that $m$-a.e. $h_x \ge 0$, hence $C_n \ge 0$ for all $n$ 
and $(s,x) \in \rr_T$. Then, with $A_n \geq \underline{a} >0$, 
the maximum principle applies to the parabolic equation \eqref{linear_PDE} 
(see e.g. \cite[Sec. 2.1, Theorem 1]{fr}, 
where our time direction is reversed compared to the setting there),
resulting with $0 \le f^{(n)}(s,x) \le \sup_x \omega(x) \le 1$ 
for all $n,r$ (as in \cite[Lemma 4.1, (i)]{dk}). Similarly, 
taking $\kappa > \sup_n \|A_n + |B_n| - C_n\|_\infty$ and
$\v^{(\pm)}=f^{(n,r)}-\phi^{(\pm)}$ for $\phi^{(\pm)}=e^{\ell \pm x+\kappa(t-s)}$, 
one has that $\LL_n \v^{(\pm)} \ge 0${,} while 
$\v^{(\pm)} (s,\pm r)=-\phi^{(\pm)}(s, \pm r) \le 0$
and $\v^{(\pm)}(t,x) = \omega(x) - \phi^{(\pm)}(t,x) \le 0$. 
Hence{,} by the maximum principle, both $\v^{(+)} \le 0$ and $\v^{(-)} \le 0$, 
yielding the bound $f^{(n,r)}(s,x) \le e^{\ell-|x|}$ (as in 
\cite[Lemma 4.1, (ii)]{dk}). Equipped with this bound, 
we follow the proof of \cite[Lemma 4.1, (iii)]{dk}. Specifically, taking
$\xi^{(\pm)}=y(r) e^{\pm \kappa x}$ for $y(r)=e^{\ell-(\kappa+1)(r-1)}$
and constant $\kappa \ge 1$ which makes $\kappa^2 A_n - \kappa |B_n| - C_n$
non-negative on $\rr_T$ for all $n$, results with  
$\LL_n \v^{(\pm),\pm} \ge 0$ for $\v^{(\pm),\pm} = \xi^{(\pm)} \pm f$ 
and $f=f^{(n,r)}$. Thus, the four functions $\v^{(\pm),\pm}$ 
satisfy the maximum principle on each of the two components of 
$\mathbb{K}_r \setminus \mathbb{K}_{r-1}$. Since $f(s,\pm r)=f(t,x)=0$ 
when $|x| \ge r-1$, while $|f(s,\pm (r-1))| \le e^{\ell+1-r}$, 
by our choice of $\xi^{(+)}$, the maximum of $\v^{(+),\pm}$ on 
the positive component of $\mathbb{K}_r \setminus \mathbb{K}_{r-1}$
is attained at $x=r$, where $\v^{(+),\pm}(s,r)=\xi^{(+)}(r)$ is 
constant. Hence, $\v^{(+),\pm}_x(s,r^-) \ge 0$ yielding that 
$|f_x(s,r^-)| \le \kappa e^\kappa e^{\ell+1-r}$. Similarly,
the maximum of $\v^{(-),\pm}$ on the negative
component of $\mathbb{K}_r \setminus \mathbb{K}_{r-1}$ is attained
at $x=-r$, where $\v^{(-),\pm}(s,-r) = \xi^{(-)}(-r)$ is constant. 
Hence, $\v^{(-)}_x (s,-r^-) \le 0$, so 
$|f_x(s,-r^{-})| \le \kappa e^\kappa e^{\ell+1-r}${,} and
$\sup_{s,n} |g_n^{\pm}(s)| \to 0$ when $r \to \infty$.
Having uniform 
ellipticity and $(A_n,B_n,C_n)$ uniformly bounded, the 
uniform bound on $\chi_n$ 
follows by applying \cite[chapter III, Theorem 11.1]{lsu} in our setting.
Finally, to bound $e_n$ we 
multiply the linear \abbr{PDE} \eqref{linear_PDE} by 
$f_{xx}$ and integrate over $\mathbb{K}_{r}$
as in the proof of \cite[Lemma 4.1, (v)]{dk}. 
Following the derivation after \cite[(4.10)]{dk}, since
$f^{(n,r)}_t(s,\pm r)=0$, integration by parts 
of the term $f_t f_{xx}$ results with  
\eq\label{eq:dk-lemm4.1-v}
\int_{\mathbb{K}_{r}^o} A_n (f^{(n,r)}_{xx})^2 \mathrm{d} m 
\le \frac{1}{2} \int_{-r}^{r} \omega'(x)^2 \mathrm{d} x 
+ \int_{\mathbb{K}_{r}^o} f^{(n,r)}_{xx} (C_n f^{(n,r)} - B_n f^{(n,r)}_x) \mathrm{d} m \,.
\en
Further, $\mathbb{K}_{r}$ is compact and $C_n f^{(n,r)}$, $B_n f^{(n,r)}_x$
uniformly bounded. Hence, by Cauchy-Schwarz inequality{,} the {rightmost} term 
of \eqref{eq:dk-lemm4.1-v} is bounded by $\kappa_2 e_{n}$ 
for some $\kappa_2=\kappa_2(r,h)$ finite. Further, the l.h.s. of 
\eqref{eq:dk-lemm4.1-v} is at least $\underline{a} e_{n}^2$, 
whereas the first term on its r.h.s. is some $\kappa_1=\kappa_1(\omega)$
finite. Consequently, $\underline{a} e_{n}^2 \le \kappa_1 + \kappa_2 e_{n}$,
yielding the desired uniform bound on $e_{n}$.
%
This completes the proof in case $h_x \geq 0$. More generally, setting 
$c=e^{\upsilon T}$ for
$\upsilon >  \| (C_\star)_- \|_\infty$, the function 
$\widetilde{f}^{(n,r)}=e^{\upsilon (s-t)} f^{(n,r)}$ satisfies 
\eqref{linear_PDE} with $\widetilde{C}_n = C_n + \upsilon \ge 0$.
Thus, by the preceding 
$\sup f^{(n,r)} \le c \sup \widetilde{f}^{(n,r)} \le c$, 
$\sup_{n} \|g^{\pm}_{n,r}\|_\infty \le c \sup_{n} \|\widetilde{g}^{\pm}_{n,r}\|_\infty
\to 0$ as $r \to \infty$, and 
$e_n \le c \|\widetilde{f}^{(n,r)}_{xx}\|_{L^2(\mathbb{K}_r)}$, 
$\chi_n \le c \|\widetilde{f}^{(n,r)}_x\|_{\infty}$ are both uniformly 
bounded in $n$, as claimed. 
\ep

\bibliographystyle{plain}

\bigskip\bigskip\bigskip

\end{document}